\RequirePackage{etoolbox}
\csdef{input@path}{{style/}{graphics/}}
\documentclass[aop,MSNbibl,seceqn,dvips]{arximspdf}
\usepackage{mathbh}
\usepackage{graphicx}

\doi{10.1214/13-AOP838} %kopijuoti is PTS
\volume{43}
\issue{2}
\pubyear{2015}
\firstpage{468}
\lastpage{527}

\makeatletter
\def\smallint{\int}
\def\colonn{\dvtx}

\newcommand{\angler}{\rangle}
\newcommand{\anglel}{\langle}
\newcommand{\rrVert}{\Vert}
\newcommand{\llVert}{\Vert}
\newcommand{\eqref}[1]{(\ref{#1})}
\newproclaim{example}{Example}[section]

\newproclaim{remark}{Remark}[section]
\newtheorem{lemma}{Lemma}[section]
\newtheorem{theorem}[lemma]{Theorem}
\newproclaim{definition}[lemma]{Definition}
\newtheorem{prop}[lemma]{Proposition}
\newtheorem{corollary}[lemma]{Corollary}
\newproclaim{assumption}[lemma]{Assumption}

\newcommand{\C}{\mathrm{C}}
\newcommand{\Z}{\mathbb{Z}}
\newcommand{\eps}{\varepsilon}
\newcommand{\dist}{\operatorname{dist}}

\newcommand{\N}{{\mathbb{N}}}
\newcommand{\R}{{\mathbb{R}}}
\newcommand{\E}{{\mathbb{E}}}
\renewcommand{\P}{{\mathbb{P}}}
\newcommand{\1}{{\mathbh{1}}}

\def\CL{\mathcal{L}}
\def\CO{\mathcal{O}}
\makeatother

\begin{document}
\begin{frontmatter}

\title{Loss of regularity for Kolmogorov equations}
\runtitle{Loss of regularity for Kolmogorov equations}
%``Numerical solutions
%of stochastic differential equations
%with nonglobally Lipschitz continuous coefficients'' and by the research
%project ``Numerical approximation
%of stochastic differential equations
%with nonglobally Lipschitz
%continuous coefficients'' both
%funded by
%the German Research Foundation.}

\begin{aug}
\author[A]{\fnms{Martin} \snm{Hairer}\thanksref{T1}},
\author[B]{\fnms{Martin} \snm{Hutzenthaler}\corref{}\thanksref{T2}\ead[label=e1]{hutzenth@math.uni-frankfurt.de}}
\and
\author[C]{\fnms{Arnulf} \snm{Jentzen}\thanksref{T3}}
\thankstext{T1}{Supported by EPSRC, the Royal Society and by the Leverhulme
Trust.}
\thankstext{T2}{Supported by the research project ``Numerical
approximation of stochastic differential equations with non-globally
Lipschitz continuous
coeffcients''.}
\thankstext{T3}{Supported by the research project ``Numerical
solutions of stochastic differential equations with non-globally
Lipschitz continuous
coeffcients''.}
\runauthor{M. Hairer, M. Hutzenthaler and A. Jentzen}
\affiliation{University of Warwick, Goethe University Frankfurt, and
ETH Zurich and~Princeton~University}
\address[A]{M. Hairer\\
Mathematics Department\\
University of Warwick\\
Coventry, CV4 7AL\\
United Kingdom} %adresu isvedimo komanda gale!
\address[B]{M. Hutzenthaler\\
Faculty of Mathematics\\
Institute for Mathematics\\
Goethe University Frankfurt\\
60325 Frankfurt\\
Germany}
\address[C]{A. Jentzen\\
Seminar f\"ur Angewandte Mathematik\\
ETH Zurich\\
8092 Z\"urich\\
Switzerland\\
and\\
Program in Applied and Computational Mathematics\\
Princeton University\\
Princeton, New Jersey 08544-1000\\
USA}
\end{aug}

% HISTORY:
\received{\smonth{10} \syear{2012}}
\revised{\smonth{1} \syear{2013}}

% ABSTRACT
%
\begin{abstract}
The celebrated
H\"ormander condition
is a sufficient
(and nearly necessary)
condition
for a second-order
linear Kolmogorov
partial differential equation
(PDE)
with smooth coefficients
to be hypoelliptic.
As a consequence, the solutions of
Kolmogorov PDEs
are smooth
at all positive times
if the coefficients
of the PDE are smooth
and satisfy
H\"ormander's condition
even if the initial function
is only continuous but not
differentiable.
First-order linear
Kolmogorov PDEs
with smooth coefficients
do not have
this smoothing effect
but at least preserve
regularity in the sense
that solutions are smooth
if their initial functions are smooth.
In this article, we consider
the intermediate regime of
nonhypoelliptic
second-order Kolmogorov PDEs
with
smooth coefficients.
% Our
% surprising
% result
%
The main observation of this article
is that
there exist counterexamples
to regularity preservation
in that case.
More precisely, we give an
example of a second-order
linear Kolmogorov PDE
with globally bounded and
smooth
coefficients
and a smooth initial function with
compact support
such that the unique
globally bounded viscosity
solution
of the PDE
is not even
locally H\"older continuous.
From the perspective of probability theory, the
existence of this example PDE has the consequence that
%observe
%In other words,
there exists a
stochastic differential equation
(SDE)
with globally bounded and
smooth coefficients
and a smooth function
with compact support
which is mapped by
the corresponding transition semigroup
to a function which is not locally H\"older continuous.
In other words, degenerate noise can have
a roughening effect.
A further implication of this
loss of regularity phenomenon
is
that numerical approximations
may converge
% slower than
% with
without any arbitrarily
small polynomial rate of convergence
to the true solution
of the SDE.
% More precisely,
% we show that there exists
% an SDE
% with globally bounded and
% smooth coefficients
% to which the
% standard
% Euler approximations
% converge
% in the strong and numerically
% weak sense without
% any arbitrarily small
% rate of convergence.
%
More precisely,
we prove for an
example SDE
with globally bounded and
smooth coefficients that
the standard
Euler approximations
converge
to the exact solution of the SDE
in the strong and numerically
weak sense, but at a rate that is slower then any power law.
% @@@@.
\end{abstract}

% KEYWORDS
\pdfkeywords{35B65,
Kolmogorov equation,
loss of regularity,
roughening effect,
smoothing effect,
hypoellipticity,
Hormander condition,
viscosity solution,
degenerate noise,
nonglobally Lipschitz continuous}

\end{frontmatter}

%s1 #&#
\section{Introduction and main results}\label{sec1}

The key observation of this
article is to reveal the
phenomenon
of \textit{loss of regularity}
in Kolmogorov
partial
differential equations (PDEs).
This observation has a
direct consequence
on the literature
on
\textit{regularity analysis
of linear PDEs},
on the literature
on \textit{regularity analysis
of stochastic differential
equations} (SDEs)
and
on the literature
on \textit{numerical approximations
of SDEs}.
We will illustrate the implications for
each field separately.

\textit{Regularity analysis of linear partial differential
equations.}
For some
$
d, m \in\mathbb{N}
$,
let
$
\mu\colonn
\R^d
\to
\R^d
$
and
$
\sigma
\colonn
\R^d \to
\R^{ d \times m }
$
be smooth
functions such that there exists
a real number $\rho> 0$
such that
$
\langle
x
,
\mu(x)
\rangle
\leq
\rho
( 1 + \|x\|^2 )
$
and\break 
% \mbox{and}
% \operatorname{tr} \Big(
% \sigma(x)
% \left[ \sigma(x) \right]^{*}++
% \Big)
$
\llVert\sigma(x) \rrVert_{L(\R^m,\R^d)}^2
\leq
\rho
( 1 + \|x\|^2 )
$
for all $ x \in\R^d $. (Here and below, we write
$
\anglel \cdot, \cdot\angler
$
and
$\llVert \cdot\rrVert
$
for the Euclidean scalar product and norm on $\R^n$.)
Let furthermore
$
\varphi\colonn\R^d \to
\R
$
be a globally bounded
and continuous function and consider the
second-order PDE
%
%e1.1 #&#
\begin{eqnarray}
\label{eqKolmogorovequationPDE} %\label{eqKolmogorovequationPDE}
\frac{ \partial}{ \partial t } u(t,x) &=& \frac{1}{2} \sum
_{ i, j = 1 }^{ d } % \left(
\sum
_{ k = 1 }^{ m } \sigma_{ i, k }(x) \cdot
\sigma_{ j, k }(x) % \right)
\cdot\frac{ \partial^2 }{ \partial x_i\, \partial x_j
} u(t,x)
\nonumber
\\[-8pt]
\\[-8pt]
\nonumber
&&{} + \sum
_{ i = 1 }^{ d } \mu_i(x) \cdot
\frac{ \partial}{ \partial x_i } u(t,x),\qquad %
u(0,x) = \varphi(x)
\end{eqnarray}
for
$
(t,x) \in(0,\infty) \times\R^d
$.
%(see, e.g.,~Theorem 3.1 in Pardoux \citationand\ Peng~
%
%Kolmogorov showed in his celebrated paper~\cite{Kolmogorov1931}
%that the density functions of the marginal distributions of a diffusion
%satisfy a partial differential equation (PDE),
%see Display (125) in~\cite{Kolmogorov1931}.
%Using integration by parts, this PDE transfers to the following
%PDE for the transition semigroup.
%
The
PDE~\eqref{eqKolmogorovequationPDE}
is referred to
as
\emph{Kolmogorov equation}
in the literature
(see, e.g.,
%Da Prato \& Zabczyk~\cite{dz02b},
Cerrai~\cite{Cerrai2001},
Da Prato~\cite{DaPrato2004},
R\"{o}ckner~\cite{Roeckner1999} and
R\"ockner and Sobol~\cite{RoecknerSobol2006};
it is also referred to as
\emph{Kolmogorov backward equation}
or
\emph{Kolmogorov PDE},
see, e.g.,
Da Prato and Zabczyk~\cite{dz02b},
{\O}ksendal~\cite{Oeksendal2000}).
It has a strong
link to probability theory
and appeared first
(in
a slightly different form;
see
display~(125)
in~\cite{Kolmogorov1931})
in Kolmogorov's celebrated
paper~\cite{Kolmogorov1931}.
%in which he laid the foundations of modern %probability theory.
Corollary~\ref{corKolmogorovviscosity}
in Section~\ref{secKolmogorovequations}
below
implies that
the
PDE~\eqref{eqKolmogorovequationPDE}
admits a unique globally
bounded viscosity solution.
More precisely,
Corollary~\ref{corKolmogorovviscosity}
proves
that there exists
a unique globally bounded continuous function
$u\colonn[0,\infty)\times\R^d\to\R$
such that
$
u|_{(0,\infty)\times\R^d}
$
is a viscosity solution
of
\eqref{eqKolmogorovequationPDE}
and
such that
$
u(0,x) = \varphi(x)
$
for all
$
x \in\R^d
$.
In this article,
we are interested to know
whether solutions $ u $ of
the
PDE~\eqref{eqKolmogorovequationPDE}
\emph{preserve regularity}
in the sense
that
$
u|_{ (0,\infty) \times\mathbb{R}^d }
$
is smooth if the
initial function
$
u(0,\cdot) =
\varphi( \cdot)
$
is smooth.
In particular, we will
answer the question whether
smoothness and global
boundedness of the initial
function
$ \varphi\colonn\R^d \to\R$
implies the existence of
a \emph{classical solution}
of the
PDE~\eqref{eqKolmogorovequationPDE}.

In the case of first-order
Kolmogorov PDEs with smooth coefficients,
that is,
$ \sigma\equiv0 $
in
\eqref{eqKolmogorovequationPDE},
regularity preservation
of solutions of
\eqref{eqKolmogorovequationPDE}
is well known.
More precisely, if
$ \sigma(x) = 0 $
for all $ x \in\R^d $
and if the initial function
$
\varphi\colonn\R^d
\to\R^d
$
in~\eqref{eqKolmogorovequationPDE}
is smooth,
then it is well known
that there exists a unique
smooth
classical solution of
\eqref{eqKolmogorovequationPDE}.
In this sense,
the
PDE~\eqref{eqKolmogorovequationPDE}
is
\emph{regularity preserving}
in the purely first-order case
$
\sigma\equiv0
$.
In the second-order case
$ \sigma\not\equiv0 $,
the situation may be even better
in the sense that the PDE~\eqref{eqKolmogorovequationPDE}
often has a \emph{smoothing effect}.
More precisely,
if
%the differential operator
%associated with
the
PDE~\eqref{eqKolmogorovequationPDE}
is
\emph{hypoelliptic},
then by definition solutions
$ u $ of the
PDE~\eqref{eqKolmogorovequationPDE}
are smooth in the sense
that
$
u|_{(0,\infty)\times\R^d}
$
is infinitely often differentiable
even if the initial function
$ u(0, \cdot) = \varphi( \cdot) $
is only continuous
%and globally bounded
but not differentiable.
%(see also Proposition~\ref{propKolmogorov}
%below).
In the seminal
paper~\cite{Hoermander1967},
H\"ormander
gave a sufficient
(and also nearly necessary;
see the discussion before
Theorem~1.1
in \cite{Hoermander1967}
and Section~2
in Hairer~\cite{Hairer2011})
condition for
%differential
%operator associated with
\eqref{eqKolmogorovequationPDE}
to be hypoelliptic;
see Theorem 1.1
in~\cite{Hoermander1967}.
To formulate
H\"ormander's condition,
set
$
\sigma_0(x)
=
\mu(x) -
\frac{1}{2}
\sum_{ k = 1 }^{ m }
\sigma'_k(x)
\sigma_k(x)
$
for all
$
x \in\R^d
$.
Then
%the
%PDE~\eqref{eqKolmogorovequationPDE}
%satisfies
the
\textit{H\"ormander condition}
is fulfilled if
%
%e1.2 #&#
\begin{eqnarray}
\label{eqhoermandercondition} &&\operatorname{span} \bigl\{ \sigma
_{ i_0 }(x), [
\sigma_{ i_0 }, \sigma_{ i_1 } ](x), \bigl[ [
\sigma_{ i_0 }, \sigma_{ i_1 } ], \sigma_{ i_2 }
\bigr](x), \ldots\in\R^d \colonn
\nonumber
\\[-8pt]
\\[-8pt]
\nonumber
&& \hspace*{ 86pt}i_0, i_1,
i_2, \ldots\in\{ 0, 1, \ldots, m \}, i_0 \neq0 \bigr\}
=\R^d
\end{eqnarray}
for all $ x \in\R^d $
where
$
[f,g]$
denotes the Lie bracket
of two smooth vector fields
$ f, g \colonn\R^d \to\R^d $.
Consequently,
if H\"ormander's
condition~\eqref{eqhoermandercondition}
is satisfied, then
the PDE~\eqref{eqKolmogorovequationPDE}
admits a unique globally bounded
smooth classical solution even if the
initial function
$
\varphi\colonn\R^d \to\R
$
is assumed to be continuous
and globally bounded only.
Clearly, there are many
cases
where
the H\"ormander
condition~\eqref{eqhoermandercondition}
fails to be fulfilled
and
where \eqref{eqKolmogorovequationPDE}
is not hypoelliptic,
for example, if $ \sigma\equiv0 $.
Next, we point out
that if all derivatives
of the drift coefficient
$ \mu$,
of the diffusion coefficient
$ \sigma$ and of
the initial function
$ \varphi$
are globally bounded
($ \mu$ and $ \sigma$
are then, in particular, globally
Lipschitz continuous), then
smoothness of the solution of the
PDE~\eqref{eqKolmogorovequationPDE}
is known even in the
nonhypoelliptic case
(see, e.g., Theorem~4.32
in Krylov~\cite{Krylov1999}
for twice differentiability
of the solution;
infinitely often differentiability
of the solution follows analogously
as in the proof of Theorem~4.32
in Krylov~\cite{Krylov1999}).
Obviously, there are
many cases where
$ \mu$ and $ \sigma$ are
not both globally Lipschitz continuous,
for example, when $ \mu$ is a polynomial
with a degree greater or equal $ 2 $
(see, e.g., Section~4 in~\cite{HutzenthalerJentzen2014Memoires}
for a list of examples).
To the best of our knowledge,
regularity of solutions of the
PDE~\eqref{eqKolmogorovequationPDE}
is in general unknown
in the
nonhypoelliptic case
if $ \sigma\not\equiv0 $
and if $ \mu$ and $ \sigma$
are not both globally Lipschitz
continuous.

In this article, we address the question
whether second-order linear PDEs
with smooth coefficients
of the form
\eqref{eqKolmogorovequationPDE}
at least preserve regularity in the nonhypoelliptic case.
The following
Theorem~\ref{thmPDEintroduction}
answers this
question to the negative.
More precisely, the key observation
of this article is to reveal the
phenomenon
of
\textit{loss of regularity}
in the sense that
the solution $ u $
of the PDE~\eqref{eqKolmogorovequationPDE}
starting with a
smooth compactly supported function
$
u(0,\cdot)
\in
\C^\infty_{\mathrm{cpt}}(\R^d,\R)
$
may turn into a nondifferentiable function
$
u(t,\cdot)
\notin
\C^1(\R^d,\R)
$
for every positive time
$ t \in(0,\infty) $.
In analogy to the well-known ``smoothing effect''
in the hypoelliptic case,
we will say
in the case of loss of regularity
that the PDE~\eqref{eqKolmogorovequationPDE}
has a \emph{roughening effect}.
Here is a simple two-dimensional
example with polynomial $\mu$
and linear $\sigma$
which has this roughening effect.
In the special case $d=2,m=1$ and
$
\mu(x)
=
(
x_1 \cdot x_2,
- x_1^2
)
$
and
$ \sigma(x)
= (0,x_2)$
for all
$
x = (x_1, x_2)
\in\mathbb{R}^2
$,
the PDE~\eqref{eqKolmogorovequationPDE}
reads as
%
%e1.3 #&#
\begin{equation}
\label{eqPDEbsp1introduction} \frac{ \partial}{ \partial t }
u(t,x) = \frac{ x_2^2
}{ 2 } \frac{ \partial^2
}{ \partial x_2^2
}u(t,x)
+ x_1 x_2 \frac{ \partial}{ \partial x_1
}u(t,x) - x_1^2
\frac{ \partial}{ \partial x_2 } u(t,x)
\end{equation}
for
$
(t,x) \in(0,\infty)
\times\R^2
$.
%Note that $\langle x,\mu(x)\rangle=0$ for all $x\in\R^2$ in this
%example.
%Thus the solution process of the associated ordinary differential
%equation
%stays on the circle centered at
%$
% (0,0) \in\R^2
%$
%going through the starting point.
%Corollary~\ref{corKolmogorovviscosity},
Theorem~\ref{thmcount1}
and
Corollary~\ref{corKolmogorovviscosity}
below
imply
that
there exists an infinitely
often
differentiable
function
$
\varphi
\in
\C_{ \mathrm{cpt} }^{ \infty
}(\R^d,\R)
$
with compact support
such
that
the unique
globally
bounded viscosity solution
$
u \colonn[0,\infty)
\times\R^2
\to\R
$
to \eqref{eqPDEbsp1introduction}
with
$
u(0, \cdot) = \varphi( \cdot)
$
has the property
that
$
u|_{ (0,\infty) \times\R^d }
$
is not differentiable and
not locally Lipschitz continuous.
%the function
%%$
%% \R^2 \ni
%% x \mapsto\E\big[ X^x(t) \big]
%% \in\R^2
%%$
%$u(t,\cdot)$
%is continuous but
%not differentiable
%and
%not locally Lipschitz continuous
%for every $ t \in(0,\infty) $.
In particular,
we thereby disprove
the existence of a
globally bounded classical
solution of the
PDE~\eqref{eqPDEbsp1introduction}
with
$ u(0, \cdot) = \varphi( \cdot) $.
Note that the drift coefficient $ \mu$
of the
PDE~\eqref{eqPDEbsp1introduction}
grows superlinearly.
One could wonder whether the roughening effect
of example~\eqref{eqPDEbsp1introduction} is due to this superlinear
growth of $\mu$. To exclude this possibility,
we prove for an example PDE with
{globally bounded and smooth coefficients} that
there exists a smooth initial function with compact support
such that the solution $ u $
is not even locally H\"older
continuous;
see Theorem~\ref{thmPDEintroduction}
below.
In particular,
Theorem~\ref{thmPDEintroduction}
implies that, in general,
the PDE~\eqref{eqKolmogorovequationPDE}
does not have a classical solution
even if the coefficients
and the initial function
are \textit{globally bounded}
and \textit{infinitely often differentiable}.
%%%%%%%%%%%%%%%%%%%%%%%%%%%%%%%%%%%%%%%%%%%%%%%%%%%%%%%%%%%%%%%%

%th1.1 #&#
\begin{theorem}[(Disprove
of the existence of classical solutions
of the Kolmogorov PDE
with smooth and globally
bounded coefficients and initial
function)]
\label{thmPDEintroduction}
There exists
a natural number $ d \in\mathbb{N} $,
a globally bounded
and infinitely often differentiable
function
$
\mu\colonn\mathbb{R}^d
\rightarrow\mathbb{R}^d
$,
a symmetric nonnegative
matrix
$
A =
( A_{i, j}
)_{
i, j \in\{ 1, 2, \ldots, d \}
}
\in\mathbb{R}^{ d \times d }
$
and
an infinitely
often differentiable function
$
\varphi
\in\C^{ \infty}_{ \mathrm{cpt} }( \R^d, \R)
$
with compact support
such that
there exists no
globally bounded
classical solution of
the PDE
%
%e1.4 #&#
\begin{eqnarray}\label{eqthmPDEintroduction}
\frac{ \partial}{ \partial t } u(t,x) & =& \sum_{ i, j = 1 }^{ d }
A_{ i, j } \cdot\frac{
\partial^2
}{
\partial x_i\, \partial x_j
} u(t,x) + \sum
_{ i = 1 }^{ d } \mu_i(x) \cdot
\frac{ \partial}{ \partial x_i } u(t,x),
\nonumber
\\[-8pt]
\\[-8pt]
\nonumber
u(0,x)&=&\varphi(x)
\end{eqnarray}
for
$
(t,x) \in(0,\infty) \times\R^d
$.
In addition,
there exists a unique
globally bounded viscosity solution
$
u \colonn[0,\infty) \times\R^d \to\R
$
of \eqref{eqthmPDEintroduction}
and
this function
% $ u $
fails to be
locally H\"older continuous.
\end{theorem}

Theorem~\ref{thmPDEintroduction}
follows immediately from
Corollary~\ref{corKolmogorovviscosity} in Section~\ref{secKolmogorovequations}
and from
Theorem~\ref{corex2bfinalcorollary}
in Section~\ref{secex2}.
More precisely,
Corollary~\ref{corKolmogorovviscosity}
and
Theorem~\ref{corex2bfinalcorollary}
imply that
there exists an infinitely differentiable
function
$
\varphi\in\C_{\mathrm{cpt}}^\infty(\R^3,\R)
$
with compact support
such that
the unique globally
bounded
viscosity solution $u\colonn[0,\infty)\times\R^3\to\R$
of the PDE
%
%e1.5 #&#
\begin{equation}
\label{eqsimplifiedversion} \frac{ \partial}{ \partial t } u(t,x)
= \frac{ \partial^2 }{ \partial x_2^2 } u(t,x) + \cos\bigl(
x_3 \exp\bigl( x_2^3 \bigr) \bigr) \cdot
\frac{ \partial}{ \partial x_1 } u(t,x)
\end{equation}
with initial condition
$ u( 0, x ) = \varphi( x ) $
for
$
(t,x) =
(t, x_1, x_2, x_3)
\in(0,\infty) \times
\R^3
$
is not locally H\"older continuous.
In particular,
the PDE~\eqref{eqsimplifiedversion}
with $ u(0,\cdot) = \varphi(\cdot) $
has no globally
bounded classical solution.
The PDE~\eqref{eqsimplifiedversion}
has a globally bounded and
highly oscillating drift coefficient
and a constant diffusion coefficient
and serves as a counterexample
to regularity preservation
for Kolmogorov PDEs.
An SDE with a globally
bounded
and highly oscillating diffusion
coefficient and a vanishing drift
coefficient has been presented
in Li and Scheutzow~\cite{LiScheutzow2011}
as a counterexample for strong
completeness of SDEs.
Another interesting observation
is that
%It is interesting to observe
the PDE~\eqref{eqsimplifiedversion}
without the second-order term
on the right-hand side of
\eqref{eqsimplifiedversion}
preserves regularity and
has a smooth classical solution
and that
the PDE~\eqref{eqsimplifiedversion}
without the first-order term
on the right-hand side of~\eqref{eqsimplifiedversion}
also preserves regularity
and
has a smooth classical
solution.
Thus,
the roughening effect
of the PDE~\eqref{eqsimplifiedversion}
is a consequence of the interplay
between the first-order and the
second-order term in
\eqref{eqsimplifiedversion}.
%
%Roughly speaking, the
%second-order term
%in \eqref{eqsimplifiedversion}
%results in unbounded
%support
%$
% | u(t,x) |
%$
%for large values
%results in an unbounded
%support for large
%increases
%brings the $ x_2 $-
%@
%unbounded
%support of the
%Wiener measure allows for large
%values of the variable
%$
% x_2 \in\R
%$
%which in turn results in
%increasingly strong oscillations
%due to
%the exponential term in
%the argument of the cosine as
%$ x_2 \to\infty
%$.
%The cube
%on the right-hand side of~\eqref{eqsimplifiedversionPDE}
%is chosen to ensure that $\E\big[\exp([W_2(t)]^3)\big]=\infty$ for all
%$t\in(0,\infty)$
%so that the oscillations of the cosine are sufficiently strong
%to roughen the solution
%function.
%We also emphasize
%that the roughening effect
%revealed here is not a
%special property
%of a particular initial
%condition
%$
% \varphi\colonn\R^3 \to\R
%$
%but a property of the
%PDE~\eqref{eqsimplifiedversion}.
%
We add that
Theorem~\ref{thmirr2}
in Section~\ref{secex2}
is
a stronger version
of Theorem~\ref{thmPDEintroduction}
in which the roughening effect
appears on every arbitrarily
small open subset
of the state space;
see
Section~\ref{secex2}
and also
Theorem~\ref{thmirr2introduction}
below
for more details.
Note that in both counterexamples
to regularity preservation
[PDE~\eqref{eqsimplifiedversion}
and PDE~\eqref{eqPDEbsp1introduction}]
it does not hold that
all derivatives of $ \mu$
and $ \sigma$
are globally bounded.
Indeed,
in both counterexamples
the drift coefficient
$ \mu$ is not globally Lipschitz
continuous.
As observed above, regularity preservation
is known if all derivatives of $ \mu$
and $ \sigma$ are globally bounded.
Moreover, note that the coefficients
in our counterexample
PDE~\eqref{eqsimplifiedversion}
are analytic functions
and
that the initial
function
$
\varphi\colonn\R^d \to\R
$
may be chosen to
be analytic
(see
Theorem~\ref{corex2bfinalcorollary}
for details).
We emphasize that this
does not contradict
the
classical Cauchy--Kovalevskaya
theorem
(e.g.,~Theorem 4.6.2
in Evans~\cite{Evans2010})
proving existence, uniqueness and
analyticity of solutions of PDEs
with analytic coefficients
as the Cauchy--Kovalevskaya
theorem
applies to
\eqref{eqthmPDEintroduction}
in the case $ A = 0 $ only.
Moreover,
we would like to
point out that
Theorem~\ref{thmPDEintroduction}
does not contradict to
Theorems 7.1.3, 7.1.4 and 7.1.7 in Evans~\cite{Evans2010},
which show the existence
of a unique classical solution of
\eqref{eqthmPDEintroduction}
if $ A $ is strictly positive
[note that~$ A $ in
\eqref{eqsimplifiedversion}
is nonnegative but not strictly positive].

Theorem~\ref{thmPDEintroduction}
shows that a general
existence
theorem for globally bounded classical solutions
of the
PDE~\eqref{eqKolmogorovequationPDE}
cannot be established.
However,
it is possible to
ensure the
existence
of a viscosity solution
of the
PDE~\eqref{eqKolmogorovequationPDE}
under rather general assumptions on
the coefficients.
More precisely,
one of our main results,
Theorem~\ref{thmKolmogorovviscosity} below,
establishes
the existence of a within a certain class
unique viscosity solution
for every second-order linear Kolmogorov PDE
whose coefficients are locally Lipschitz continuous
and satisfy the Lyapunov-type inequality~\eqref{eqLyapunovcondition2}.
To the best of our knowledge,
this is the first result in the literature
proving existence and uniqueness of solutions
of the Kolmogorov PDE~\eqref{eqKolmogorovequationPDE}
in the above generality;
see also the discussion after
Theorem~\ref{thmKolmogorovviscosity}
for a short review of existence and uniqueness results
for Kolmogorov PDEs.
A crucial result on the route to
Theorem~\ref{thmKolmogorovviscosity}
is
the uniqueness result of Corollary~\ref{coruniqueness2}
for viscosity solutions
of degenerate parabolic
second-order linear
PDEs.

The roughening effect
of the
PDE~\eqref{eqKolmogorovequationPDE}
revealed in this first paragraph
of this \hyperref[sec1]{Introduction}
has a direct consequence
on the literature on regularity
analysis of SDEs.
This is the subject of the
next paragraph.

\textit{Regularity analysis of stochastic differential
equations.}
For the rest of this \hyperref[sec1]{Introduction},
we use the following notation.
Let
$ ( \Omega, \mathcal{F}, \P ) $
be an arbitrary probability space with
a normal filtration
$
( \mathcal{F}_t )_{ t \in[0,\infty) }
$
which supports
a standard
$
( \mathcal{F}_t
)_{ t \in[0,\infty) }
$-Brownian motion
$
W \colonn[0,\infty) \times
\Omega\to\R^m
$
with continuous sample paths.
It is a classical result that
the above
assumptions on $\mu$ and $\sigma$
ensure the existence
of a family
$
X^x = (X^x_1, \ldots, X^x_d)
\colonn[0,\infty)
\times\Omega
\rightarrow\R^d
$,
$ x \in\R^d $,
of
up to indistinguishability unique
solution processes
(see, e.g.,
Theorem~3.1.1 in
\cite{PrevotRoeckner2007})
with continuous sample
paths
of the SDE
%
%e1.6 #&#
\begin{equation}
\label{eqSDEintroduction} d X^x(t) = \mu\bigl( X^x(t) \bigr) \,dt +
\sigma\bigl( X^x(t) \bigr) \,dW(t)
\end{equation}
for $ t \in(0,\infty) $
and $ x \in\R^d $
and with
$
X^x(0) = x
$
for all $ x \in\R^d $
(see, e.g., Theorem~1 in
Krylov~\cite{Krylov1990}).
Here, the function
$
\mu\colonn
\R^d \to\R^d
$
is the infinitesimal mean and
the function
$
\sigma\cdot
\sigma^{*}
\colonn
\R^d \to\R^{d\times d}
$
is the infinitesimal covariance
matrix of the
SDE~\eqref{eqSDEintroduction}.
It is also well known that
the coercivity assumption
on $ \mu$ and
the linear growth
bound on $ \sigma$
additionally imply
moment bounds
$
\sup_{
x \in\{y\in\R^d\colonn
\| y \| \leq p\}
}
\E[
\sup_{ t \in[0,p] }
\| X^x(t) \|^p
] < \infty
$
for all
$ p \in[0,\infty) $
for the solution processes
of the
SDE~\eqref{eqSDEintroduction}.
The \textit{transition semigroup}
$
P_t \colonn
\C_b(\R^d,\R) \to\C_b(\R^d,\R)
$,
$ t \in[0,\infty) $
of the SDE~\eqref{eqSDEintroduction}
is defined by
$
( P_t \varphi)(x)
:=
\E[
\varphi(X^x(t))
]
$
for all
$ t \in[0,\infty) $,
$ x \in\R^d $ and all
$ \varphi\in\C_b(\R^d,\R) $
where $ \C_b( \R^d, \R) $
is as usual the space of
globally bounded and continuous functions from $\R^d$ to $\R$.
Note for every
$
\varphi\in\C_b( \R^d, \R)
$
that
the function
$
\R^d \ni x \mapsto
\E[
\varphi(X^x(t))
]
\in\R
$
is continuous
(see, e.g., Theorem~1.7
in
Krylov~\cite{Krylov1999})
and hence,
the semigroup
$
( P_t )_{ t \in[0,\infty) }
$
is well defined.
Observe also that
the function
$
\R^d \ni x \mapsto
\E[
\varphi(X^x(t))
]
\in\R
$
is continuous
for every
$
\varphi\in\C_b( \R^d, \R)
$
although the
SDE~\eqref{eqSDEintroduction}
is, in general,
not strongly complete;
see
Li and Scheutzow~\cite{LiScheutzow2011}
and see, for example,
also Elworthy~\cite{Elworthy1978},
Kunita~\cite{Kunita1990}
and Fang,
Imkeller and Zhang~\cite{FangImkellerZhang2007}
for further
results on strong completeness
of SDEs.

Theorem 1.1 in
H\"ormander~\cite{Hoermander1967}
and Proposition~\ref{propKolmogorov}
below
%in Subsection~\ref{ssecSolutionsofKolmogorovequationsinthe%distributionalsense}
imply that if
the H\"ormander condition~\eqref{eqhoermandercondition}
is fulfilled,
then
the semigroup is smoothing
in the sense
that
$
P_t (
\C_{ b }( \R^d, \R)
)
\subseteq
\C^{ \infty}_{ b }( \R^d, \R)
$
for all $t\in(0,\infty)$.
To the best of our knowledge,
it remained an open question
in the nonhypoelliptic case
whether SDEs
with infinitely often
differentiable coefficients
such as
\eqref{eqSDEintroduction}
in general
preserve regularity
in the sense that
$
P_t (
\C^{ \infty}_{ b }( \R^d, \R)
)
\subseteq
\C^{ \infty}_{ b }( \R^d, \R)
$
for all $t\in(0,\infty)$.
This article answers this
question to the negative.
More precisely,
the following theorem
reveals
%a key observation
%of this article is to reveal the
%phenomenon
%of
%in the sense
that
smooth functions
with compact support
may be mapped to nonsmooth functions
by the transition
semigroup
of the
SDE~\eqref{eqSDEintroduction}.
In analogy to the well-known ``smoothing effect''
of many SDEs,
we will say that the semigroup
has a \emph{roughening effect} in that case.
Here is a simple two-dimensional example SDE with polynomial drift coefficient
and linear diffusion coefficient
which has this roughening effect.
In the special case $d=2$, $m=1$ and
$\mu(x)
= (
x_1 \cdot x_2,
- x_1^2
)$
and
$
\sigma(x)
=
(0, x_2)
$
for all
$
x = (x_1, x_2)
\in\mathbb{R}^2
$, the SDE~\eqref{eqSDEintroduction} reads as
%
%e1.7 #&#
\begin{eqnarray}
\label{eqSDEbsp1introduction} d X_1^x(t) &=& X^x_1(t)
\cdot X_2^x(t) \,dt,
\nonumber
\\[-8pt]
\\[-8pt]
\nonumber
d X_2^x(t) &= &- X^x_1(t)^2
\,dt + X_2^x(t) \,dW(t)
\end{eqnarray}
for $ t \in(0,\infty) $
and $ x \in\mathbb{R}^2 $.
Observe that
\eqref{eqPDEbsp1introduction}
is the
Kolmogorov
PDE
of
\eqref{eqSDEbsp1introduction};
see
Corollary~\ref{corKolmogorovviscosity}
for details.
Moreover, note that
$
\langle x,\mu(x)\rangle=0
$
for all
$
x \in\R^2
$
in this example.
Thus, the solution process of the associated ordinary differential equation
stays on the circle centered at $(0,0)\in\R^2$ going through the
starting point.
Theorem~\ref{thmcount1}
in Section~\ref{secex1}
shows for the SDE~\eqref{eqSDEbsp1introduction}
that
there exists an infinitely
often
differentiable
function
$
\varphi
\in
\C_{ \mathrm{cpt} }^{ \infty
}(\R^d,\R)
$
with compact support
such
for every $ t \in(0,\infty) $
the functions
$
\R^2
\ni
x
\mapsto
\E[
\varphi(X^x(t))
] \in\R
$
and
$
\R^2
\ni
x
\mapsto
\E[
X^x(t)
] \in\R^2
$
are continuous
but not differentiable
and not locally Lipschitz
continuous.
For every
$ t \in(0,\infty)$,
we hence have
the roughening effect
$
P_t (
\C^{ \infty}_{ \mathrm{cpt} }( \R^d, \R)
)
\nsubseteq
\C^{ 1 }( \R^d, \R)
$
in the case of the
SDE~\eqref{eqSDEbsp1introduction}.
The drift coefficient $ \mu$
of the
SDE~\eqref{eqSDEbsp1introduction}
grows superlinearly.
As above,
the superlinear growth of $ \mu$
is not necessary for the transition
semigroup of the SDE
to be roughening.
This is subject of the next
main result of this article.

%%%%%%%%%%%%%%%%%%%%%%%%%%%%%%%%%%%%%%%%%%%%%%%%%%%%%

%th1.2 #&#
\begin{theorem}[(A counterexample
to regularity preservation with
degenerate additive noise)]
\label{thmirr2introduction}
There exists
a natural number
$ d \in\N$,
a globally bounded
and infinitely often differentiable
function
$
\mu\colonn\R^d
\rightarrow\R^d
$
and
a constant
function
$
\sigma\colonn\R^d
\to\R^{ d \times d }
$, that is
$
\sigma(x) = \sigma(0)
$
for all $ x \in\R^d $,
with the following properties.
For every
$ t \in(0,\infty) $
the function
$
\R^d
\ni
x \mapsto
\mathbb{E} [
X^x(t)
]
\in\R^d
$
is continuous but nowhere locally
H\"{o}lder continuous
and for every nonempty open
set $ O \subset\R^d $
there exists an infinitely
often differentiable function
$
\varphi
\in
\C^{ \infty}_{ \mathrm{cpt} }( \R^d, \R)
$
with compact support
such that
the function
$
O
\ni
x \mapsto\break 
\mathbb{E} [
\varphi( X^x(t) )
]
\in\R
$
is continuous
but not locally
H\"{o}lder
continuous.
In particular,
for every
$
t \in(0,\infty)
$
we have
$
P_t (
\C^{ \infty}_{ \mathrm{cpt} }( \R^d, \R)
)
\nsubseteq
\bigcup_{ \alpha\in(0,\infty) }
\C^{ \alpha}( \R^d, \R)
$.
\end{theorem}

Theorem~\ref{thmirr2introduction}
follows immediately from
Theorem~\ref{thmirr2}
in Section~\ref{secex2}.
The roughening effect of
some SDEs with smooth
coefficients
revealed through
example~\eqref{eqSDEbsp1introduction}
and
Theorem~\ref{thmirr2introduction}
above,
has a direct consequence
on the literature
on numerical approximations
of SDEs.
This is the subject of the next paragraph.

\textit{Numerical approximations
of stochastic differential equations.}
Starting with Maruyama's adaptation
of Euler's method to SDEs in 1955 (see~\cite{m55}),
an extensive literature
on
the numerical approximation
of solutions of
SDEs
has been published
in the last six decades;
see, for example, the books and
overview articles
\cite{kp92,kps94,m95,g02,BurrageEtAl2004,mt04,mr08,jk09d,KloedenNeuenkirch2013}
for extensive lists of references.
A key objective
in this field
of research is
to prove convergence
of suitable
numerical
approximation processes
to the solution process
of the SDE and
to establish a rate of convergence
for the considered
approximation scheme
in the strong,
in the almost sure or
in the numerically weak sense.

Almost sure convergence
rates of many numerical
schemes such as the
standard Euler method
or the higher-order Milstein method
are well known
for the
SDE~\eqref{eqSDEintroduction}
and even for a much larger
class of nonlinear SDEs;
see
Gy\"{o}ngy~\cite{g98b}
and Jentzen,
Kloeden and Neuenkirch~\cite{jkn09a}.
Many applications,
however,
require the
numerical approximation
of moments or other functionals
of the solution process,
for instance,
the expected pay-off of an
option in computational finance;
see, for example,
Glasserman~\cite{g04} for details.
For this reason,
%particular interest is
applications are particularly
interested in
strong and numerically weak
convergence rates.
The vast majority of research results
establishing
strong and
numerically weak
convergence rates
assume that the coefficients
of the SDE are
globally Lipschitz
continuous
or at least
that they satisfy the
\textit{global
monotonicity condition}
that there exists a real number
$ \rho\in\R$ such that
$
\langle x - y, \mu(x) - \mu(y)
\rangle
+
\frac{ 1 }{ 2 }
\sum_{ k = 1 }^{ m }
\| \sigma_k(x) - \sigma_k(y)
\|^2
\leq
\rho\| x - y \|^2
$
for all $ x, y \in\R^d $
(see, e.g.,
Theorem~2.4 in Hu~\cite{h96},
Theorem~5.3 in
Higham, Mao and Stuart~\cite{hms02},
Schurz~\cite{Schurz2006},
Theorems 2 and 3 in
Higham and Kloeden~\cite{hk07},
Theorem~6.3 in Mao and Szpruch~\cite{MaoSzpruch2013Rate},
Theorem 1.1 in Hutzenthaler, Jentzen and Kloeden~\cite
{HutzenthalerJentzenKloeden2012},
Theorem~3.2 in
Wang and Gan~\cite{WangGan2013}).
%For instance,
%for the Euler
%method the
%strong convergence
%rate $ \frac{ 1 }{ 2 } $
%and the numerically
%weak convergence rate $ 1 $
%are established
%in the case of globally
%Lipschitz continuous coefficients
%of the SDE.
%Other numerical approximation schemes
%such as the Milstein
%scheme or the implicit Euler method
%convergence with the same
%or with higher convergence rates
%under these assumptions
Strong
and numerically
weak convergence
rates
without assuming
global monotonicity
are established in
%Theorem 2.1 in
Gy{\"o}ngy and R{\'a}sonyi~\cite{GyoengyRasonyi2011}
in the case of a class of scalar SDEs
with globally H\"{o}lder continuous
coefficients,
in
D\"orsek~\cite{Doersek2012}
in the case of the two-dimensional
stochastic Navier--Stokes equations
%in Theorem 1.1
and in Dereich,
Neuenkirch and Szpruch~\cite{DereichNeuenkirchSzpruch2012},
\mbox{Alfonsi}~\cite{Alfonsi2012},
Neuenkirch and Szpruch~\cite{NeuenkirchSzpruch2012}
in the case of a class of scalar SDEs (including, e.g.,
the Cox--Ingersoll--Ross process) that
can be transformed in a suitable sense to
SDEs that satisfy the global monotonicity assumption.
The global monotonicity assumption
is a serious restriction on the coefficients
of the SDE and excludes
many interesting SDEs in the literature
(e.g., stochastic Lorenz equations,
stochastic Duffing--van der Pol oscillators
and the stochastic SIR model;
see Section~4 in~\cite{HutzenthalerJentzen2014Memoires}
for details and further examples).
It remains an
open problem to establish
strong and numerically
weak convergence rates
in the general setting of the
SDE~\eqref{eqSDEintroduction}.

In this article, we establish
in the setting \eqref{eqSDEintroduction}
the existence of an SDE
with \textit{globally
bounded} and
\textit{infinitely often differentiable
coefficients}
for which the Euler approximations
converge in the
strong and in the numerically
weak sense
without
% slower than
any arbitrarily small
polynomial rate of
convergence.
More precisely,
our main result for the
literature on the numerical
approximation
of SDEs is the following theorem.
%%
%%%%%%%%%%%%%%%%%%%%%%%%%%%%%%%%%

%th1.3 #&#
\begin{theorem}[(A counterexample
to the rate of convergence
in the numerical approximation
of nonlinear SDEs with
additive noise)]
\label{thmnonrateintroduction}
Let
$ T \in(0,\infty) $
and $ x_0 \in\R^4 $
be arbitrary.
Then there exists
a globally bounded
and infinitely often differentiable
function
$
\mu\colonn\mathbb{R}^4
\rightarrow\mathbb{R}^4
$
and
a symmetric nonnegative matrix
$
B \in\R^{ 4 \times4 }
$
such that
the stochastic process
$
X \colonn[0,T] \times
\Omega\to\R^4
$
with continuous sample
paths
satisfying
$
X(t)
=
x_0
+
\int_0^t \mu( X(s) ) \,ds
+
B W(t)
$
for all
$ t \in[0,T] $
and its Euler--Maruyama
approximations
$
Y^N \colonn
\{ 0, \frac{ T }{ N },
\frac{ 2 T }{ N },
\ldots, T \}
\times\Omega\to
\R^4
$,
$ N \in\N$,
satisfying
$ Y^N(0) = x_0 $
and\vadjust{\goodbreak}
$
Y^N( \frac{ ( n + 1 ) T }{ N }
)
=
Y^N(
\frac{ n T }{ N }
)
+\break 
\mu(
Y^N(
\frac{ n T }{ N }
)
)
\frac{ T }{ N }
+
B
(
W_{ (n + 1) T / N }
- W_{ n T / N }
)
$
for all
$ n \in\{ 0, 1, \ldots, N - 1 \} $,
$ N \in\N$,
fulfill that
%
%e1.8 #&#
\begin{eqnarray}
\label{eqlimitnumerik} &&\lim_{ N \to\infty} \bigl( N^{ \alpha}
\cdot\E
\bigl[ \bigl\| X(T) - Y^N(T) \bigr\| \bigr] \bigr) \nonumber\\
&&\qquad= \lim_{ N \to\infty}
\bigl( N^{ \alpha} \cdot\bigl\llVert\E\bigl[ X(T) \bigr] - \E
\bigl[
Y^N(T) \bigr] \bigr\rrVert\bigr)
\\
&&\qquad = \cases{0, & \quad$\alpha= 0,$
\vspace*{2pt}
\cr
\infty,&\quad  $\alpha> 0,$}\nonumber
\end{eqnarray}
for all $ \alpha\in[0,\infty) $.
In particular, for every
$ \alpha\in(0,\infty) $
there exists no real number
$
c_{ \alpha}
\in(0,\infty)
$
such that
$
\|
\E[
X(T)
]
-
\E[
Y^N(T)
]
\|
\leq
c_{ \alpha}
\cdot N^{ - \alpha}
$
for all $ N \in\N$.
\end{theorem}

Theorem~\ref{thmnonrateintroduction}
follows immediately from
Theorem~\ref{thmnonrate}
in Section~\ref{secex4}.
In the deterministic case $\sigma\equiv0$,
it is well known that
the Euler approximations
converge to the solution
process
of \eqref{eqSDEintroduction}
with the rate $ 1 $.
In the stochastic case $\sigma\not\equiv0$,
this rate of convergence
can often not be achieved.
In particular,
Clark and Cameron~\cite{ClarkCameron1980}
% were the first to
proved for an SDE
in the setting
of \eqref{eqSDEintroduction}
that
a class of
Euler-type schemes
cannot, in general, converge
strongly with a higher-order than $ \frac{ 1 }{ 2 } $.
Since then, there have been many results
on lower bounds of strong and numerically weak approximation
errors for numerical approximation
schemes of SDEs;
see, for example,
\cite{Ruemelin1982,CambanisHu1996,hmr00a,hmr00b,dg01,m02,mr07a,mr08,hjk11,Kruse2011}
and the references therein.
Now
the
%surprising
observation of
Theorem~\ref{thmnonrateintroduction}
is
that
%also the convergence order $\frac{1}{2}$
%(which is the standard convergence order in the
%globally Lipschitz case)
%cannot always be achieved with the Euler method.
%Even more,
there exist
SDEs with
smooth and globally
bounded coefficients
for which the standard Euler approximations
converge in the strong and
numerically weak sense
\textit{without
% slower than
any arbitrarily small
polynomial rate of convergence}.
To the best of our
knowledge,
Theorem~\ref{thmnonrateintroduction}
is the first result in the literature
in which it has been
established that
Euler's method
%for SDEs
%the Euler approximations
converges to the solution
of an SDE with smooth
coefficients in the strong
and numerical weak sense
% slower than
without
any arbitrarily
small polynomial rate of convergence.
Clearly, this lack of
a rate of convergence
is not
a special property of the
Euler scheme and
holds for other schemes
such as the Milstein scheme, too.
It is based on
%a consequence of
the fact that
our counterexample SDE
for Theorem~\ref{thmnonrateintroduction}
[see \eqref{eqex3SDE}]
%to which
%we show that Euler's method
%converges in the strong
%and numerically weak sense
%sl....wer than any arbitrarily small
%polynomial rate of convergence
suffers under the
\textit{roughening
effect} revealed in
Theorems~\ref{thmPDEintroduction}
and \ref{thmirr2introduction}
(see Corollary~\ref{corex4initial}
and Theorem~\ref{thmnonrate}
in Section~\ref{secex4}
for details).

Comparing
Theorem~\ref{thmnonrate}
with
Theorem~2.4
in Gy\"{o}ngy~\cite{g98b}
reveals the remarkable difference
that the
Euler approximations for some SDEs
have almost sure convergence
rate
$
\frac{ 1 }{ 2 }-$ but no strong and no numerically
weak rate of convergence.
More formally,
Theorem~2.4 in
\cite{g98b} shows
in the setting of
Theorem~\ref{thmnonrateintroduction}
that there exist
finite random variables
$
C_{ \varepsilon} \colonn\Omega
\to[0,\infty)
$,
$
\varepsilon\in(0,\frac{1}{2})
$,
such that
$
\| X(T) - Y^N(T) \|
\leq
C_{ \varepsilon}
\cdot
N^{ ( \varepsilon- { 1 }/{ 2 } )
}
$, $\P$-a.s. for
all $ N \in\N$
and all
$ \varepsilon\in(0,\frac{1}{2}) $.
Taking expectation then
results in
$
\E[
\| X(T) - Y^N(T) \|
]
\leq
\E[ C_{ \varepsilon} ]
\cdot
N^{ ( \varepsilon- { 1 }/{ 2 } )
}
$
for all $ N \in\N$
and all
$ \varepsilon\in(0,\frac{1}{2}) $
and from
Theorem~\ref{thmnonrateintroduction}
we hence get that
the error constants
have infinite expectations, that is,
$
\E[
C_{ \varepsilon}
] = \infty
$
for all
$
\varepsilon\in(0,\frac{1}{2})
$.
In addition, we refer to Theorem~2.3
in Milstein and Tretyakov~\cite{mt05}
for a weak convergence
result restricted to certain subevents of the
probability space.
Finally, we emphasize
that Monte Carlo simulations
confirm
the slow strong
and numerically weak convergence
phenomenon of Euler's
method revealed
in Theorem~\ref{thmnonrateintroduction}.
For details, the reader is referred to
Figure~\ref{fnonrate}
in Section~\ref{secex4}
below.

%The remainder of this
%article is organized
%as follows. $ \dots$

%s2 #&#
\section{Counterexamples to regularity preservation with linear multiplicative
noise}
\label{secex1}

In this section,
we
establish the phenomenon of
loss of regularity
of the simple example SDE~\eqref{eqSDEbsp1introduction}
with polynomial drift coefficient and linear diffusion coefficient.
For this, we consider
the following setting.
Let
$
( \Omega, \mathcal{F},
\mathbb{P} )
$
be a probability space
with a normal filtration
$
( \mathcal{F}_t )_{
t \in[0,\infty)
}
$,
let
$
W \colonn[0,\infty)
\times\Omega\rightarrow
\mathbb{R}
$
be a one-dimensional standard
$
( \mathcal{F}_t )_{
t \in[0,\infty)
}
$-Brownian motion with
continuous sample paths
and let
%$
% \mu, \sigma\colonn\mathbb{R}^2
% \rightarrow\mathbb{R}^2
%$
%be two functions
%given by
% \mu(x)
% = \left(
% \begin{array}{cc}
% x_1 \cdot x_2 \\
% - \left( x_1 \right)^2
% \end{array}
% \right)
%
% \mbox{and}
%
% \sigma(x)
% = \left(
% \begin{array}{cc}
% 0 \\
% x_2
% \end{array}
% \right)
%for all
%$
% x = (x_1, x_2)
% \in\mathbb{R}^2
%$.
%Note that the condition
%%\begin{equation}
%%\label{eqcoercivity}
% $\left< x, \mu(x) \right> = 0$
%%\end{equation}
%for all $ x \in\mathbb{R}^2 $
%ensures
%the existence of a family
%$
% X^x = (X^x_1, X^x_2)
% \colonn
% [0,\infty)
% \times\Omega
% \rightarrow\mathbb{R}^2
%$,
%$ x \in\mathbb{R}^2 $,
%of up to indistinguishability
%unique adapted stochastic
%processes with continuous
%sample paths
%satisfying
%$
% X^x(t) =
% x
% +
% \int_0^t
% \mu( X^x(s) )
% ds
% +
% \int_0^t
% \sigma( X^x(s) )
% dW(s)
%$
%for all
%$ t \in[0,\infty) $
%$ \mathbb{P} $-a.s.\ and all
%$ x \in\mathbb{R}^2 $;
%see, e.g.,
%Corollary 2.6 in Gy\"ongy \citationand\ Krylov~
%The stochastic processes
$
X^x =
( X^x_1, X^x_2 )
\colonn[0,\infty)
\times\Omega\rightarrow\mathbb{R}^2
$,
$ x \in\mathbb{R}^2 $,
be the up to indistinguishability
unique solution processes
with continuous sample paths
of the SDE
% adapted stochastic processes
% with continuous sample paths
% satisfying
%solution processes
%of the SDE
%
%e2.1 #&#
\begin{eqnarray}
\label{eqSDEbsp1} d X_1^x(t) &=& X^x_1(t)
\cdot X_2^x(t) \,dt,
\nonumber
\\[-8pt]
\\[-8pt]
\nonumber
d X_2^x(t) &=& - \bigl( X^x_1(t)
\bigr)^2 \,dt + X_2^x(t) \,dW(t)
\end{eqnarray}
for $ t \in(0,\infty) $
and $ x \in\mathbb{R}^2 $
satisfying
$ X^x(0) = x $
for all $ x \in\R^2 $.
Corollary~2.6 in Gy\"ongy
and Krylov~\cite{GyoengyKrylov1996}
ensures that the processes
$
X^x \colonn[0,\infty) \times\Omega
\to\R^2
$, $ x \in\R^2 $,
do indeed exist.
%@@ Ok this way?
The following
Theorem~\ref{thmcount1}
shows that
the semigroup associated with the SDE~\eqref{eqSDEbsp1}
loses regularity in the
sense that there exists an
infinitely often differentiable function
with compact support, which is mapped
to a nonsmooth function
by the semigroup.

\begin{theorem}[(A counterexample
to regularity preservation with
linear multiplicative noise)]
\label{thmcount1}
Let
$
X^x
\colonn
[0,\infty)
\times\Omega
\rightarrow\mathbb{R}^2
$,
$ x \in\mathbb{R}^2 $,
be solution processes
of the SDE~\eqref{eqSDEbsp1}
with continuous sample paths and
with
$ X^x(0) = x $
for all $ x \in\R^2 $.
Then
$
\sup_{
x\in\{ y \in\R^2 \colonn\| y \| \leq p \}
}
\E[
\sup_{
t \in[0,p]
}
\|
X^{ x }(t)
\|^p
]
< \infty
$
for all $ p \in[0,\infty) $
and there exists
an infinitely often differentiable
function
$
\varphi
\in\C^{ \infty}_{ \mathrm{cpt} }( \R^2, \R)
$
with compact support
such that
for every $ t, p \in(0,\infty) $
the mappings
$
\mathbb{R}^2
\ni x
\mapsto
\mathbb{E} [
X^x(t)
]
\in\mathbb{R}^2
$,
$
\mathbb{R}^2
\ni x
\mapsto
\mathbb{E} [
\varphi( X^x(t) )
]
\in\mathbb{R}
$
and
$
\mathbb{R}^2
\ni x
\mapsto
X^x(t)
\in
L^p( \Omega; \mathbb{R}^2 )
$
are continuous but
not locally Lipschitz continuous
and not differentiable.
\end{theorem}

The proof
of
Theorem~\ref{thmcount1}
is deferred
to the end of this section.
The proof of
Theorem~\ref{thmcount1}
uses the following simple
lemma.
\begin{lemma}[(Restricted
exponential integrals
of a geometric Brownian motion)]
\label{leWintegrate1}
Let
$
( \Omega, \mathcal{F},
\mathbb{P} )
$
be a probability space and
let
$
W \colonn[0,\infty)
\times\Omega
\rightarrow\mathbb{R}
$
be a one-dimensional
standard Brownian
motion with continuous sample paths.
Then
%
%e2.2 #&#
\begin{equation}
\mathbb{E} \biggl[ \mathbh{1}_{
\{
a \leq
e^{ W(t) }
\leq b
\}
} \exp\biggl( c \cdot\int
_0^t e^{ W(s) } \,ds \biggr) \biggr] =
\infty
\end{equation}
for all
$
t, a, b, c \in(0,\infty)
$
with $ a < b $.
\end{lemma}
\begin{pf}%{Proof
%of Lemma~\ref{leWintegrate1}}
Independence of
$ W(t) $
from
$
(
W(s) - \frac{ s }{ t } W(t)
)_{
s \in[0,t]
}
$
for all
$ t \in(0,\infty) $
implies
%
%e2.3 #&#
\begin{eqnarray}
\label{eqestimatelemeWintegrate1} && \mathbb{E} \biggl[ \mathbh{1}_{
\{
a \leq
e^{ W(t) }
\leq b
\}
} \exp\biggl( c
\cdot\int_0^t e^{ W(s) } \,ds \biggr)
\biggr]\nonumber\\ %\\ & =
% \mathbb{E} \left[
% \mathbh{1}_{
% \left\{
% a \leq
% e^{ W(t) }
% \leq b
% \right\}
% }
% \exp\left(
% c \cdot
% \int_0^t
% e^{
% \left( W(s) - \frac{ s }{ t } W(t)
% \right)
% }
%
% e^{
% \frac{ s }{ t } W(t)
% }
% ds
% \right)
% \right]
% \\ &
&&\qquad\geq
\mathbb{E} \biggl[ \mathbh{1}_{
\{
a \leq
e^{ W(t) }
\leq b
\}
} \exp\biggl( c \cdot\int
_0^t e^{
( W(s) - ({ s }/{ t }) W(t)
)
} a^{
{ s }/{ t }
} \,ds
\biggr) \biggr]
\nonumber
\\
&&\qquad \geq% \mathbb{E} \left[
% \mathbh{1}_{
% \left\{
% a \leq
% e^{ W(t) }
% \leq b
% \right\}
% }
% \exp\left(
% c \cdot
% \min( a, 1 )
% \cdot
% \int_0^t
% e^{
% \left( W(s) - \frac{ s }{ t } W(t)
% \right)
% }
% ds
% \right)
% \right]
\mathbb{P} \bigl[ a \leq e^{ W(t) }
\leq b \bigr] \cdot\mathbb{E} \biggl[ \exp\biggl( tc \cdot\min(
a, 1 ) \cdot
\frac{1}{t} \int_0^t e^{
( W(s) - ({ s }/{ t }) W(t)
)
}
\,ds \biggr) \biggr]
\\
&&\qquad \geq\mathbb{P} \bigl[ a \leq e^{ W(t) } \leq b \bigr] \nonumber\\
&&\qquad\quad{}\times
\mathbb{E}
\biggl[ \exp\biggl( tc \cdot\min( a, 1 ) \cdot\exp\biggl( \frac{1}{t}
\int_0^t W(s) - \frac{ s }{ t } W(t) \,ds
\biggr) \biggr) \biggr]\nonumber
\end{eqnarray}
for all
$
t, a, b, c \in(0,\infty)
$
with $ a < b $ where we
used Jensen's inequality and
convexity of the exponential
function in the last step.
The time integrated Brownian bridge
$
\int_0^t W(s)-\frac{s}{t} W(t) \,ds
$
on the right-hand side of
\eqref{eqestimatelemeWintegrate1}
is
%as the limit
% \int_0^t W(s)-\frac{s}{t} W(t) ds=\lim_{n\to\infty}\sum_{k=0}^{n-1}
%of a linear functional of a multivariate normal distribution again
normally distributed with mean $0$ and variance
%
%e2.4 #&#
\begin{eqnarray}
&& \E\biggl[ \biggl( \int_0^t W(s) -
\frac{ s }{ t } W(t) \,ds \biggr)^{ 2 } \biggr] \nonumber\\
&&\qquad= \E\biggl[ \int
_0^t \int_0^t
\biggl( W(s) - \frac{ s }{ t } W(t) \biggr) \biggl( W(r) - \frac{
r }{ t }
W(t) \biggr) \,dr \,ds \biggr]
\nonumber\\
&&\qquad = \int_0^t \int_0^t
\E\biggl[ W(s) W(r) - \frac{ r }{ t } W(s) W(t) - \frac{ s }{ t } W(r)
W(t)
\nonumber
\\[-8pt]
\\[-8pt]
\nonumber
&&\hspace*{205pt}{} + \frac{ s r }{ t^2 } \bigl( W(t) \bigr)^2 \biggr] \,dr \,ds
\\
&&\qquad = \int_0^t \int_0^t
\biggl( \min( r, s) - \frac{ r s }{ t } - \frac{ s r }{ t } +
\frac{ s r }{ t } \biggr) \,dr \,ds % \\ &
\nonumber\\
&&\qquad= 2 \int_0^t
\int_0^s \biggl( r - \frac{ r s }{ t }
\biggr) \,dr \,ds = \int_0^t \biggl(
s^2 - \frac{ s^3 }{ t } \biggr) \,ds = % \frac{ t^3 }{ 3 } -
% \frac{ t^3 }{ 4 }
% =
\frac{ t^3 }{ 12 } \in(0,\infty)\nonumber
\end{eqnarray}
for every
$ t \in(0,\infty) $.
As the double exponential normal
distribution has infinite mean,
we conclude that the
right-hand side of~\eqref{eqestimatelemeWintegrate1} is infinite
for all $t,a,b,c\in(0,\infty)$. This finishes the proof Lemma~\ref
{leWintegrate1}.
\end{pf}

The proof of the following Lemma~\ref{lemcount1}
makes use of
Lemma~\ref{leWintegrate1}.
Using Lemma~\ref{lemcount1},
the proof of
Theorem~\ref{thmcount1}
is then completed at the end
of this section.

%le2.3 #&#
\begin{lemma}
\label{lemcount1}
Let
$
X^x
\colonn
[0,\infty)
\times\Omega
\rightarrow\mathbb{R}^2
$,
$ x \in\mathbb{R}^2 $,
be solution processes
of the SDE~\eqref{eqSDEbsp1}
with continuous sample paths and
with
$ X^x(0) = x $
for all $ x \in\R^2 $.\vadjust{\goodbreak}
Then
$
\sup_{
x\in\{ y \in\R^2 \colonn\| y \| \leq p \}
}
\E[
\sup_{
t \in[0,p]
}
\|
X^{ x }(t)
\|^p
]
< \infty
$
for all $ p \in[0,\infty) $
and
%
%e2.5 #&#
\begin{eqnarray}
\label{eqpstrongirregularity}&& \lim_{0\neq x_1\to0 } \biggl( \frac
{ 1 }{ x_1 } \cdot
\mathbb{E} \bigl[ X^{(x_1,x_2)}_1(t) % \right]
- % \mathbb{E} \left[
X^{(0,x_2)}_1(t) \bigr] \biggr)
\nonumber
\\[-8pt]
\\[-8pt]
\nonumber
&&\qquad= \infty= \lim
_{0\neq x_1\to0 } \biggl( \frac{
1
}{ | x_1 | } \cdot\bigl\| X^{(x_1, x_2)}_1(t)
- X^{(0, x_2)}_1(t) \bigr\|_{
L^p( \Omega; \mathbb{R} )
} \biggr)
\end{eqnarray}
for all
$ t, x_2, p \in(0,\infty) $
and
there exists
an infinitely often differentiable
function
$
\varphi
\in\C^{ \infty}_{ \mathrm{cpt} }( \R^2, \R)
$
with compact support\vadjust{\goodbreak}
such that
$
\lim_{0\neq x_1\to0 }
(
\frac{ 1 }{ x_1 } \times\break 
\mathbb{E} [
\varphi (
X^{ (x_1, x_2) }(t)
)
% \right]
-
% \mathbb{E}\left[
\varphi (
X^{(0,x_2)}(t)
)
]
)
=
\infty
$
for all
$ t, x_2 \in(0,\infty) $.
\end{lemma}

\begin{pf}%{Proof
%of Lemma~\ref{lemcount1}}
%
% Corollary~2.6 in Gy\"ongy
% \citationand\ Krylov~\cite{GyoengyKrylov1996}
% guarantees the existence of
% a family of
% up to indistinguishability
% unique adapted stochastic
% processes
% $
% X^x \colonn[0,\infty)
% \times\Omega
% \rightarrow\R^2
% $,
% $ x \in\R^2 $,
% with continuous sample
% paths
% satisfying~\eqref{eqSDEbsp1}.
% Moreover,
%@@ OK this way? @@.
The global Lipschitz continuity of
$ \sigma$, the local Lipschitz
continuity of $ \mu$ and
$\anglel x, \mu(x) \angler = 0$
for all $ x \in\mathbb{R}^2 $
imply
that
\[
\sup_{
x\in\{ y \in\R^2 \colonn\| y \| \leq p \}
}
\E\Bigl[
\sup_{
t \in[0,p]
}
\bigl\|
X^{ x }(t)
\bigr\|^p
\Bigr]
< \infty
\]
for all $ p \in[0,\infty) $.
%(see, e.g., @@).
Next,
we disprove local Lipschitz continuity of
the mapping
$
\mathbb{R}^2
\ni
x \mapsto
X^x_1(t)
\in L^p( \Omega; \mathbb{R} )
$
for every
$ t, p \in(0,\infty) $.
More precisely,
aiming at a contradiction,
we assume that
the second equality in~\eqref{eqpstrongirregularity}
is false.
Then there exist positive real numbers $t,x_2,p\in(0,\infty)$ and a
sequence of real numbers
$
h_n
\in
\R
\setminus
\{ 0 \}
$,
$ n \in\N$,
such that
$
\lim_{ n \to\infty}
h_n = 0
$
and such that
$
\limsup_{n\to\infty}
\frac{
1
}{
\vert h_n \vert
}
\|
X_1^{(h_n,x_2)}(t)
-
X_1^{(0,x_2)}(t)
\|_{
L^p( \Omega; \mathbb{R} )
}
<
\infty
$.
Theorem~1.7 in Krylov\break \cite{Krylov1999}
(see also
Proposition 3.2.1 in Pr{\'e}v{\^o}t and R\"ockner~\cite{PrevotRoeckner2007})
yields that\break $\sup_{s\in[0,t]}\|X^{(h_n,x_2)}(s)-X^{(0,x_2)}(s)\|\to0$
in probability as $n\to\infty$.
Hence,\vspace*{1pt} there exists a
strictly increasing
sequence
$ n_k \in\N$,
$ k \in\N$,
of natural numbers
such that
$
\lim_{ k \to\infty}
\sup_{ s \in[0,t] }
\|
X^{ (h_{n_k},x_2) }(s) -
X^{(0,x_2)}(s)
\|
= 0$, $\P$-a.s.;
see, for example,~Corollary 6.13 in Klenke~\cite{Klenke2008}.
Applying this,
Fatou's lemma
and Lemma~\ref{leWintegrate1}
implies
%
%e2.6 #&#
\begin{eqnarray}
\label{eqstrongirregularity} \infty& >& \limsup_{ k\to\infty}
\biggl(
\frac{
1
}{
\vert h_{n_k} \vert
} \bigl\| X_1^{(h_{n_k},x_2)}(t) - X_1^{(0,x_2)}(t)
\bigr\|_{
L^p( \Omega; \mathbb{R} )
} \biggr) \nonumber\\
&= &\limsup_{ k \to\infty} \biggl(
\frac{
1
}{
\vert h_{n_k}
\vert
} \bigl\| X_1^{(h_{n_k},x_2)}(t) \bigr\|_{
L^p( \Omega; \mathbb{R} )
} \biggr)
\nonumber\\
&=& \limsup_{ k \to\infty} \biggl\llVert\exp\biggl( \int
_0^t X^{(h_{n_k},x_2)}_2(s) \,ds
\biggr) \biggr\rrVert_{
L^p( \Omega; \mathbb{R} )
} \nonumber\\
&\geq&\biggl\llVert\liminf
_{ k\to\infty} \biggl\{ \exp\biggl( \int_0^t
X^{(h_{n_k},x_2)}_2(s) \,ds \biggr) \biggr\} \biggr\rrVert
_{
L^p( \Omega; \mathbb{R} )
}
\nonumber
\\
\nonumber
& =& \biggl\llVert\exp\biggl( \int_0^t
X^{(0,x_2)}_2(s) \,ds \biggr) \biggr\rrVert_{
L^p( \Omega; \mathbb{R} )
}\\
& =&
\biggl\llVert\exp\biggl( \int_0^t
e^{
(
W(s)
-
s/2
)
} \,ds \cdot x_2 \biggr) \biggr\rrVert_{
L^p( \Omega; \mathbb{R} )
}
\\
& \geq&\biggl( \mathbb{E} \biggl[ \exp\biggl( \int_0^t
e^{
W(s)
} \,ds \cdot\frac{
p x_2
}{
e^{ t/2 }
} \biggr) \cdot\1_{
\{
1 \leq e^{ W(t) } \leq2
\}
}
\biggr] \biggr)^{
1 / p
} \nonumber\\
&=& \infty.\nonumber
\end{eqnarray}
This contradiction implies that
the second equality in~\eqref{eqpstrongirregularity}\vspace*{1pt}
is true.
The first equality
in \eqref{eqpstrongirregularity}
follows from
the second equality
in \eqref{eqpstrongirregularity}
as
$
\frac{ 1 }{ x_1 }
(
X_1^{ (x_1, x_2) }(t)
-
X_1^{ (0, x_2) }(t)
)
\in[0, \infty)
$
for all
$
x_1 \in
\R\setminus
\{ 0 \}
$
and all
$
x_2\in(0,\infty)
$.
%@@ OK this way?
In the next step, let
$ c \in(0,\infty) $
be an arbitrary
fixed real number and
let
$
\psi_1
\colonn
\mathbb{R}
\rightarrow\mathbb{R}
$
and
$
\psi_2
\colonn
\mathbb{R}
\rightarrow[0,\infty)
$
be two infinitely often differentiable
functions with
$
x \cdot\psi_1(x) \geq0
$
for all $ x \in\R$,
with
$
\psi_1(x) = \psi_2(x) = 0
$
for all
$
x \in
\mathbb{R} \setminus
[-c-1,c+1]
$
and
with
$
\psi_1( x ) = x
$
and
\mbox{$
\psi_2( x ) = 1
$}
for all $ x \in[ - c, c ] $.
Due to partition of
unity, such functions indeed
exist.
Next, let
$
\varphi\colonn\mathbb{R}^2
\rightarrow\mathbb{R}
$
be given by
$
\varphi(x_1,x_2)
=
\psi_1( x_1 ) \cdot
\psi_2( x_2 )
$
for all
$
x = (x_1, x_2) \in\mathbb{R}^2
$.
Note that
$
\varphi
% \colonn\mathbb{R}^2
% \rightarrow\mathbb{R}
\in
\C^{ \infty}_{ \mathrm{cpt} }(
\mathbb{R}^2,
\mathbb{R}
)
$
is an infinitely often
differentiable function
with compact support.
We now show that
$
\lim_{ 0 \neq x_1\to0 }
(
\frac{ 1 }{ x_1 } \cdot
\mathbb{E} [
\varphi (
X^{ (x_1, x_2) }(t)
)
-
\varphi (
X^{(0,x_2)}(t)
)
]
)
=
\infty
$
for all
$ t, x_2 \in(0,\infty) $.
Aiming at a contradiction,
assume that
there exist positive
real numbers
$ t, x_2 \in(0,\infty) $
and a sequence
$
h_n \in
\R
\setminus
\{ 0 \}
$,
$ n \in\N$,
such that
$ \lim_{ n \to\infty} h_n = 0 $
and such that
%
%e2.7 #&#
\begin{equation}
\limsup_{n\to\infty} \biggl( \frac{
1
}{ h_n } \cdot\E\bigl[
\varphi\bigl( X_1^{ (h_n, x_2) }(t) \bigr) - \varphi\bigl(
X_1^{ (0,x_2) }(t) \bigr) \bigr] \biggr) < \infty.
\end{equation}
Theorem~1.7 in
Krylov~\cite{Krylov1999}
yields that
$
\sup_{ s \in[0,t] }
\|
X^{(h_n,x_2)}(s)-
X^{(0,x_2)}(s)
\|\to0
$
in probability as $n\to\infty$.
Hence, there exists a
strictly increasing\vspace*{1pt}
sequence
$ n_k \in\N$, $ k \in\N$,
of natural
numbers
such that
$
\lim_{ k \to\infty}
$
$
\sup_{ s \in[0,t] }
$
$
\|
X^{ ( h_{ n_k }, x_2)
}(s) -
X^{ (0, x_2) }(s)
\|
= 0
$, $\mathbb{P} $-a.s.;
see, for example,
Corollary~6.13 in
Klenke~\cite{Klenke2008}.
Applying this,
% the identity
% \begin{equation}
% \left\langle\big(\nabla\varphi\big)(0,x_2),y\right\rangle
% =
% \psi_1'( 0 ) \cdot
% \psi_2( x_2 ) \cdot
% y_1
% +
% \psi_1( 0 ) \cdot
% \psi_2'( x_2 ) \cdot
% y_2
% =
% \psi_2( x_2 ) \cdot
% y_1
% \end{equation}
% for all
% $
% x_2 \in\R
% $
% and all
% $
% y = (y_1, y_2)
% \in\mathbb{R}^2
% $,
the fact
$
% (
% \varphi(x_1,x_2) -
% \varphi(0,x_2)
% )/
% x_1
\frac{
1
}{
x_1
}
(
\varphi(x_1,x_2) -
\varphi(0,x_2)
)
% =
% \frac{
% \varphi(x_1,x_2)
% }{
% x_1
% }
\in[0,\infty)
$
for all $ x_1\in\R\setminus\{ 0 \} $
and all $ x_2\in(0,\infty) $,
%@@ Ok this way?
Fatou's lemma
and Lemma~\ref{leWintegrate1}
then results in
%
%e2.8 #&#
\begin{eqnarray}
\infty& >& \limsup_{
k \to\infty
} \biggl( \frac{
1
}{
h_{n_k}
}
\mathbb{E} \bigl[ \varphi\bigl( X^{(h_{n_k},x_2)}(t) \bigr) -
\varphi\bigl(
X^{(0,x_2)}(t) \bigr) \bigr] \biggr)\nonumber\\
& = &\limsup_{
k\to\infty
}
\mathbb{E} \biggl[ \biggl\vert\frac{
\varphi( X^{(h_{n_k},x_2)}(t) )
-
\varphi( X^{(0,x_2)}(t) )
}{ h_{n_k} } \biggr\vert\biggr]
\nonumber\\
& \geq&\mathbb{E} \biggl[ \liminf_{
k\to\infty
} \biggl\vert
\frac{
\varphi( X^{(h_{n_k},x_2)}(t) )
-
\varphi( X^{(0,x_2)}(t) )
}{ h_{n_k} } \biggr\vert\biggr]\nonumber\\
& =& \mathbb{E} \biggl[ \liminf
_{
k\to\infty
} \biggl( \frac{
\varphi( X^{(h_{n_k},x_2)}(t) )
-
\varphi( X^{(0,x_2)}(t) )
}{ h_{n_k} } \biggr) \biggr]
\nonumber
\\
\nonumber
& =& \mathbb{E} \biggl[ \psi_2 \bigl( X^{ (0, x_2) }_2(t)
\bigr) \biggl( \liminf_{k\to\infty} \frac{
X^{(h_{n_k},x_2)}_1(t)
% -
% X^{(0,x_2)}_1(t)
}{h_{n_k}} \biggr)
\biggr]\\
& =& \mathbb{E} \biggl[ \psi_2 \bigl( X^{ (0,x_2) }_2(t)
\bigr) \cdot\exp\biggl( \int_0^t
e^{
(
W(s) - s/2
)
} \,ds \cdot x_2 \biggr) \biggr]
\\
& \geq&\mathbb{E} \biggl[ \mathbh{1}_{
\{
{c}/{2}
\leq
x_2 \cdot
\exp(
W(t) - t / 2
)
\leq c
\}
} \cdot\exp\biggl( \int
_0^t e^{
(
W(s) - s/2
)
} \,ds \cdot x_2
\biggr) \biggr] \nonumber\\
&=& \infty.\nonumber
\end{eqnarray}
This contradiction implies
that
$
\lim_{ 0 \neq x_1\to0 }
(
\frac{ 1 }{ x_1 } \cdot
\mathbb{E} [
\varphi (
X^{ (x_1, x_2) }(t)
)
-\break 
\varphi (
X^{(0,x_2)}(t)
)
]
)
=
\infty
$
for all
$ t, x_2 \in(0,\infty) $.
The proof
of Lemma~\ref{lemcount1}
is thus completed.
\end{pf}

\begin{pf*}{Proof
of Theorem~\ref{thmcount1}}
Theorem~1.7 in
Krylov~\cite{Krylov1999}
(see also
Proposition~3.2.1 in Pr{\'e}v{\^o}t and R\"ockner~\cite{PrevotRoeckner2007}),
in particular, shows
for every
$ t \in[0,\infty) $
that the mapping
%
%e2.9 #&#
\begin{equation}
\label{eqcontinuityprobability} \R^2 \ni x \mapsto X^x(t) \in
L^0\bigl( \Omega; \R^2 \bigr)
\end{equation}
is continuous.
This implies
for every
$
\varphi\in
\C^{ \infty}_{ \mathrm{cpt} }( \R^2, \R)
$
and every
$ t \in[0,\infty) $
that the mapping
$
\R^2 \ni x
\mapsto
\E[
\varphi( X^x(t) )
]
\in\R
$
is continuous.
Moreover, Lemma~\ref{lemcount1}
proves that
$
\sup_{
x \in
\{ y \in\R^2 \colonn\| y \| \leq p \}
}
$
$
\E[
\sup_{
t \in[0,p]
}
\|
X^{ x }(t)
\|^p
]
< \infty
$
for all $ p \in[0,\infty) $.
Combining this,
\eqref{eqcontinuityprobability},
Corollary~6.21
in Klenke~\cite{Klenke2008}
and
Theorem~6.25
in Klenke~\cite{Klenke2008}
shows for every
$ t, p \in[0,\infty) $
that the mappings
$
\R^2 \ni x
\mapsto
X^x(t) \in L^p( \Omega; \R^2 )
$
and
$
\R^2 \ni x
\mapsto
\E[ X^x(t)
]
\in\R^2
$
are continuous.
Furthermore,
Lemma~\ref{lemcount1}
implies that
there exists
an infinitely often differentiable
function
$
\varphi
\in\C^{ \infty}_{ \mathrm{cpt} }( \R^2, \R)
$
with compact support
such that
for every $ t, p \in(0,\infty) $
the mappings
$
\mathbb{R}^2
\ni x
\mapsto
\mathbb{E} [
X^x(t)
]
\in\mathbb{R}^2
$,
$
\mathbb{R}^2
\ni x
\mapsto
\mathbb{E} [
\varphi( X^x(t) )
]
\in\mathbb{R}
$
and
$
\mathbb{R}^2
\ni x
\mapsto
X^x(t)
\in
L^p( \Omega; \mathbb{R}^2 )
$
are not locally Lipschitz
continuous and not differentiable.
The proof of
Theorem~\ref{thmcount1}
is thus completed.
\end{pf*}

%We remark that
%the solution of simpler second-order linear
%PDEs exhibit the roughening effect as well.
%For example,
%consider
%the SDE
In the remainder of this section,
we briefly consider slightly modified
versions of the SDE~\eqref{eqSDEbsp1}.
The generator of the SDE~\eqref{eqSDEbsp1}
is nowhere elliptic.
We remark that the phenomenon
of loss of regularity
may also appear for an SDE
whose generator is
in many points of the
state space
elliptic.
For example, let
$
( \Omega, \mathcal{F}, \P
)
$
be a probability space
with a normal filtration
$
( \mathcal{F}_t )_{ t \in[0,\infty) }
$,
let
$
W = (W_1, W_2)
\colonn[0,\infty)
\times\Omega\to
\R^2
$
be a two-dimensional
standard
$
( \mathcal{F}_t )_{ t \in[0,\infty) }
$-Brownian motion
and let
$
X^x = (X^x_1, X^x_2)
\colonn
[0,\infty)
\times\Omega
\rightarrow\mathbb{R}^2
$,
$ x \in\mathbb{R}^2 $,
be the up to indistinguishability
unique solution processes
with continuous sample paths
of the SDE
%
%e2.10 #&#
\begin{eqnarray}
\label{eqSDEbsp1b} d X_1^x(t) &=& X^x_1(t)
\cdot X_2^x(t) \,dt + X_1^x(t)
\,dW_1(t),
\nonumber
\\[-8pt]
\\[-8pt]
\nonumber
d X_2^x(t) &=& - \bigl( X^x_1(t)
\bigr)^2 \,dt + X_2^x(t) \,dW_2(t)
\end{eqnarray}
for $ t \in(0,\infty) $
and $ x \in\mathbb{R}^2 $
satisfying
$ X^x(0) = x $
for
all $ x \in\R^2 $.
The generator of the
SDE~\eqref{eqSDEbsp1b}
is in every point
$
x = (x_1,x_2) \in\R^2
$
with
$
x_1 \cdot x_2 \neq0
$
elliptic
but
%for the SDE~\eqref{eqSDEbsp1b}
%that
there exists a
function
$
\varphi\in
\C^{ \infty}_{ \mathrm{cpt} }( \R^d, \R)
$
such that
for every
$ t \in(0,\infty) $
the functions
$
\R^2 \ni x \mapsto
\E[
X^x(t)
]
\in\R^2
$
and
$
\R^2 \ni x \mapsto
\E[
\varphi( X^x(t) )
]
\in\R
$
are not locally Lipschitz
continuous.
The proof of this statement
is completely analogous
as in the case
of the SDE~\eqref{eqSDEbsp1}.
Furthermore, the same statement
holds if the two independent
standard Brownian motion
in \eqref{eqSDEbsp1b} are
replaced by one and the same
standard Brownian motion.
More precisely, if
$
( \Omega, \mathcal{F},
\mathbb{P} )
$
is a probability space
with a normal filtration
$
( \mathcal{F}_t )_{
t \in[0,\infty)
}
$
and if
$
W \colonn[0,\infty)
\times\Omega\rightarrow
\mathbb{R}
$
is a one-dimensional standard
$
( \mathcal{F}_t )_{
t \in[0,\infty)
}
$-Brownian motion,
then the up to indistinguishability
unique solution processes
$
X^x
= (X^x_1, X^x_2 )
\colonn[0,\infty) \times\Omega
\to\R^2
$,
$ x \in\R^2 $,
of the SDE
%
%e2.11 #&#
\begin{equation}
\label{eqSDEbsp1c} d X^x(t) =\pmatrix{
X^x_1(t) \cdot X_2^x(t)
\vspace*{2pt}\cr
- \bigl( X^x_1(t) \bigr)^2
} \,dt + X^x(t) \,dW(t)
\end{equation}
for $ t \in(0,\infty) $
and $ x \in\mathbb{R}^2 $
with continuous sample paths
and with $ X^x( 0 ) = x $
for all $ x \in\R^2 $
fulfill that
there exists a function
$
\varphi\in
\C^{ \infty}_{ \mathrm{cpt} }( \R^2, \R)
$
such that
for every
$ t \in(0,\infty) $
the functions
$
\R^2 \ni x \mapsto
\E[
X^x(t)
]
\in\R^2
$
and
$
\R^2 \ni x \mapsto
\E[
\varphi( X^x(t) )
]
\in\R
$
are not locally Lipschitz
continuous.

%s3 #&#
\section{Counterexamples to regularity preservation with
degenerate additive
noise}
\label{secex2}

In this section, we show
the roughening effect for an example SDE
with globally bounded and infinitely often differentiable coefficients.
For this, it suffices to consider the following counterexample to
regularity preservation.
Let
$
( \Omega, \mathcal{F}, \P
)
$
be a probability space,
let
$
W \colonn[0,\infty) \times\Omega
\to\R
$
be a one-dimensional
standard Brownian
motion
and let
$
X^x =
( X^x_1, X^x_2, X^x_3 )
\colonn[0,\infty) \times
\Omega\to\R^3
$,
$ x \in\R^3 $,
be the up to indistinguishability
unique solution processes
with continuous sample paths
of the SDE
%
%e3.1 #&#
\begin{eqnarray}
\label{eqex2bSDE}
\nonumber
d X^x_1(t) & =& \cos\bigl(
X^x_3(t) \cdot\exp\bigl( X^x_2(t)^3
\bigr) \bigr) \,dt,
\\
d X^x_2(t) & = &\sqrt{2} \,dW(t),
\\
\nonumber
d X^x_3(t) & = &0 \,dt
\end{eqnarray}
for $ t \in[0,\infty) $
and $ x \in\R^3 $
satisfying
$ X^x(0) = x $
for all $ x \in\R^3 $.
Observe that
%
%e3.2 #&#
\begin{equation}
\label{eqexactrepsec3} X^{ x }_1(t) = x_1 + \int
_0^t \cos\bigl( x_3 \cdot\exp
\bigl( \bigl[ x_2 + \sqrt{2} W(s) \bigr]^3 \bigr) \bigr)
\,ds,
\end{equation}
$\P$-a.s. for
all $ t \in[0,\infty) $
and all $ x = (x_1, x_2, x_3) \in\R^3 $.
%Note also that the Kolmogorov
%PDE associated with~\eqref{eqex2bSDE}
%reads as
% \frac{ \partial}{ \partial t }
% u(t,x)
%=
% \frac{ \partial^2 }{ \partial x_2^2 }
% u(t,x)
% +
% \cos\left(
% x_3
%% \cdot
% \exp\left( [ x_2 ]^3 \right)
% \right)
% \cdot
% \frac{ \partial}{ \partial x_1 }
% u(t,x)
%for
%$ t \in(0,\infty) $
%and
%$ x = (x_1, x_2, x_3) \in\R^3 $.
%Corollary~\ref{corKolmogorovviscosity}
%below establishes
%existence and uniqueness
%of viscosity solutions
%of the PDE~\eqref{eqex2bKolmogorov}.
%

%th3.1 #&#
\begin{theorem}[(A counterexample
to regularity preservation with
degenerate additive noise)]
\label{corex2bfinalcorollary}
Let
$ T \in(0,\infty) $
and let
$
X^x \colonn[0,\infty) \times\Omega
\to\R^3
$,
$ x \in\R^3 $,
be solution processes
of the SDE~\eqref{eqex2bSDE}
satisfying
$
X^x(0) = x
$
for all $ x \in\R^3 $.
Then
there exists
an infinitely often differentiable
function
$
\varphi\in
\C^{ \infty}_{ \mathrm{cpt} }( \R^3, \R)
$
with compact support
such that
for every $ t \in(0,T] $
the functions
$
\R^3 \ni x \mapsto
\mathbb{E} [
X^x(t)
]
\in\R^3
$
and
$
\R^3 \ni x \mapsto
\mathbb{E} [
\varphi( X^x(t) )
]
\in\R
$
are
continuous but
not locally
H\"{o}lder continuous.
\end{theorem}

%The proof of
%Theorem~\ref{corex2bfinalcorollary}
%is deferred
%to the end of this section.
%For every $ t \in(0,T] $
%the function
%$
% \R^3
% \ni
% x \mapsto
% \mathbb{E}\big[
% \varphi( X^x(t) )
% \big]
% \in\mathbb{R}
%$
%in Theorem~\ref{thmirr2}
%is continuous and globally
%bounded but
%not locally H\"{o}lder
%continuous.
%For illustration we add a simple
%example of a function with
%such properties.
%More precisely,
%the function
%$
% \psi\colonn
% [ - \frac{ 1 }{ 9 }, \frac{ 1 }{ 9 } ]
% \rightarrow
% [0,\infty)
%$
%given by
%$ \psi(0) = 0 $
%and
%$
% \psi(x) =
% |x|^{
% \left\{
% 1 / \ln( \ln( 1 / |x| ) )
% \right\}
% }
% =
% \exp\big(
% \frac{ \ln(|x|) }{
% \ln( \ln( \frac{1}{|x|} ) )
% }
% \big)
%$
%for all
%$
% x \in
% [ - \frac{ 1 }{ 9 }, \frac{ 1 }{ 9 } ]
% \setminus\{ 0 \}
%$
%is continuous and globally bounded
%but not H\"{o}lder continuous.

In the following,
regularity properties
of the solution processes
$
X^x = (X^x_1, X^x_2,  X^x_3)
\colonn[0,\infty) \times
\Omega\to\R^3
$,
$ x \in\R^3 $,
of the SDE~\eqref{eqex2bSDE}
are investigated
in order to prove
Theorem~\ref{corex2bfinalcorollary}.
To \,do so, we first establish
a few auxiliary results.
We begin with a simple lemma on
trigonometric integrals.

%le3.2 #&#
\begin{lemma}
\label{lem1rotateNEW}
Let
$ a, b \in\mathbb{R} $
be real numbers
with $ a < b $,
let
$ \psi\colonn[a,b]
\rightarrow[0,\infty) $
be a continuously differentiable
function
and let
$
\varphi\colonn[a,b] \rightarrow\mathbb{R}
$
be a twice continuously differentiable
function
with
$ e^{ i \cdot\varphi(a) } = i $
and with
$ \varphi'(x) \geq0 $,
$ \varphi''(x) \geq0 $
and
$
\psi'(x) \leq0
$
for all $ x \in[a,b] $.
Then
$
\int_a^b
\cos( \varphi(x) ) \psi(x)
\,dx
\leq0
$.
\end{lemma}

\begin{pf}
%{Proof
%of Lemma~\ref{lem1rotateNEW}}
First, assume
w.l.o.g. that
$
\varphi(b) \geq\varphi(a) + \pi
$
(otherwise we have
$ \cos( \varphi(x) ) \leq0 $
for all $ x \in[a,b]$, and
hence $
\int_a^b
\cos( \varphi(x) ) \psi(x)
\,dx
\leq0
$).
Moreover,
assume w.l.o.g. that
$ \varphi'(x) > 0 $
for all $ x \in(a,b] $
(otherwise consider
$
\varphi|_{ [\tilde{a}, b] }
\colonn
[ \tilde{a}, b ]
\to
\R
$
where
$
\tilde{a}:=
\inf(
\{
x \in[a,b] \colonn
\varphi'(x) > 0
\}
\cup
\{ b \}
)
$
and observe that
$
\int_a^b \cos( \varphi(x) )
\psi(x) \,dx
=
\int_{ \tilde{a} }^b \cos( \varphi(x) )
\psi(x) \,dx
$).
In particular,
$
\varphi\colonn[a,b]
\rightarrow\mathbb{R}
$
is strictly
increasing and
there exists a
unique
strictly increasing
continuous function
$
\varphi^{ - 1 }
\colonn
[ \varphi(a), \varphi(b) ]
\rightarrow
[a,b]
$
with
$
\varphi^{-1}( \varphi(x) ) = x
$
for all $ x \in[a,b] $
and with
$
\varphi( \varphi^{-1}(x) ) = x
$
and
$
( \varphi^{-1} )'( x )
=
\frac{ 1 }{
\varphi'( \varphi^{-1}(x) )
}
> 0
$
for all $ x \in( \varphi(a),\varphi(b) ) $.
Integration by substitution
and integration by parts
therefore imply
%
%e3.3 #&#
\begin{eqnarray}
\label{eqlemcos1NEW} &&\int_a^b \cos\bigl(
\varphi(x) \bigr) \psi(x) \,dx \nonumber\\
&&\qquad= \int_{
\varphi(a)
}^{
\varphi(b)
} \cos(
x ) \cdot\psi\bigl( \varphi^{-1}(x ) \bigr) \cdot\bigl(
\varphi^{-1}\bigr)'(x ) \,dx \nonumber\\
&&\qquad= \int_{
\varphi(a)
}^{
\varphi(b)
}
\frac{
\cos( x ) \cdot
\psi( \varphi^{-1}(x) )
}{
\varphi'( \varphi^{-1}(x) )
} \,dx
\nonumber
\\[-8pt]
\\[-8pt]
\nonumber
&&\qquad = \frac{
[
\sin( \varphi(b) ) - 1
]
\psi(
\varphi^{-1}( \varphi(b) )
)
}{
\varphi' (
\varphi^{-1}( \varphi(b) )
)
}\nonumber\\
&&\qquad\quad{} - \int_{
\varphi(a)
}^{
\varphi(b)
} \bigl[
\sin( x ) - 1 \bigr] \biggl[ \frac{
\psi'( \varphi^{-1}(x) )
}{
[ \varphi'( \varphi^{-1}(x) ) ]^2
} - \frac{
\psi( \varphi^{-1}(x) )
\varphi''( \varphi^{-1}(x) )
}{
[
\varphi'( \varphi^{-1}(x) )
]^3
} \biggr] \,dx\nonumber\\
&&\qquad\leq0.\nonumber
\end{eqnarray}
This completes the proof
of Lemma~\ref{lem1rotateNEW}.
\end{pf}

The next lemma analyzes
suitable regularity properties
of the solution processes
$
X^x = (X^x_1, X^x_2, X^x_3)
\colonn[0,\infty) \times
\Omega\to\R^3
$,
$ x \in\R^3 $,
of the SDE~\eqref{eqex2bSDE}.
Its proof is based on
%an application
%of
Lemma~\ref{lem1rotateNEW}.

%le3.3 #&#
\begin{lemma}[(A lower bound)]
\label{lemex2lowerbound}
Let
$
( \Omega, \mathcal{F}, \P
)
$
be a probability space
and let
$
W \colonn[0,\infty)
\times\Omega\to\R
$
be a one-dimensional
standard Brownian
motion. Then
%
%e3.4 #&#
\begin{equation}
\label{eqlowerboundcos} 1 - \mathbb{E} \bigl[ \cos\bigl( h \cdot
\exp\bigl( \bigl[ x +
W(t) \bigr]^3 \bigr) \bigr) \bigr] \geq\exp\biggl(
\frac{
- 8
}{ t } \biggl[ \biggl| \ln\biggl( \frac{ \pi}{ 2 h } \biggr)
\biggr|^{ 2 / 3 } + x^2 \biggr] \biggr)
\end{equation}
for all
$
h
\in
( 0,
\frac{ \pi}{ 2 }
\exp (
-
| [ \sqrt{t} + x ] \vee0 |^3
)
]
$,
$ t \in(0,\infty) $
and all
$ x \in\R$
and
%
%e3.5 #&#
\begin{eqnarray}
&&\int_0^t \mathbb{E} \bigl[
\mathbh{1}_{
\{
W(t) \in A
\}
} \bigl( 1 - \cos\bigl( h \cdot e^{
[ x + W(s) ]^3
} \bigr)
\bigr) \bigr] \,ds
\nonumber
\\[-8pt]
\\[-8pt]
\nonumber
&&\qquad\geq\frac{ t }{ 3 } \cdot\E\bigl[ \mathbh{1}_{
\{
W(t) \in A
\}
}
e^{
{ - 64
\vert W(t) \vert^2 }/{ t }
} \bigr] \cdot\exp\biggl( \frac{
- 64
}{ t } \biggl[ \biggl| \ln
\biggl( \frac{ \pi}{ 2 h } \biggr) \biggr|^{ 2 / 3 } + x^2 \biggr]
\biggr)
\end{eqnarray}
for all
$
h
\in
( 0,
\frac{ \pi}{ 2 }
\exp (
-
[
\sqrt{ t } + | x | +
\sup_{ a \in A } | a |
]^3
)
]
$,
$ x \in\R$,
$ t \in(0,\infty) $
and all bounded
and Borel measurable sets
$ A \subset\R$.
\end{lemma}

\begin{pf}%{Proof
%of
%Lemma~\ref{lemex2lowerbound}}
First of all, define a family
$
\varphi_{ t, x, h } \colonn
[
\frac{
[
\ln( \pi/ ( 2 h ) )
]^{ 1 / 3 }
- x
}{
\sqrt{ t }
}
,
\infty
)
\to\R
$,\break 
$
(t,x,h) \in
\{
(0,\infty) \times\R\times(0,\infty)
\colonn
h \leq
\frac{ \pi}{ 2 }
\exp( - | x \vee0 |^3 )
\}
$,
of functions by
%
%e3.6 #&#
\begin{equation}
\varphi_{ t, x, h }( y ):= h \cdot\exp\bigl( [ x + \sqrt{ t } y
]^3 \bigr)
\end{equation}
for all
$
y \in
[
\frac{
[
\ln( \pi/ ( 2 h ) )
]^{ 1 / 3 }
- x
}{
\sqrt{ t }
}
,
\infty
)
$,
$ t \in(0,\infty) $,
$
h \in
( 0,
\frac{ \pi}{ 2 }
\exp( - | x \vee0 |^3 )
]
$
and all
\mbox{$ x \in\R$}.
Observe that
%
%e3.7 #&#
\begin{equation}
\varphi_{ t, x, h }'( y ) = 3 \sqrt{ t } [ x + \sqrt{ t } y
]^2 \varphi_{ t, x, h }( y ) \geq0
\end{equation}
and
%
%e3.8 #&#
\begin{equation}
\varphi_{ t, x, h }''( y ) = 6 t [ x + \sqrt{ t
} y ] \varphi_{ t, x, h }( y ) + 9 t [ x + \sqrt{ t } y ]^4
\varphi_{ t, x, h }( y ) \geq0
\end{equation}
for all
$
y \in
[
\frac{
[
\ln( \pi/ ( 2 h ) )
]^{ 1 / 3 }
- x
}{
\sqrt{ t }
}
,
\infty
)
$,
$ t \in(0,\infty) $,
$
h \in
( 0,
\frac{ \pi}{ 2 }
\exp( - | x \vee0 |^3 )
]
$
and all
$ x \in\R$.
In addition, note that
$
\varphi_{ t, x, h } (
\frac{
[
\ln( \pi/ ( 2 h ) )
]^{ 1 / 3 }
- x
}{
\sqrt{ t }
}
)
= \frac{ \pi}{ 2 }
$
for all
$ t \in(0,\infty) $,
$
h \in
( 0,
\frac{ \pi}{ 2 }
\exp( - | x \vee0 |^3 )
]
$
and all
$ x \in\R$.
We can thus apply
Lemma~\ref{lem1rotateNEW}
to obtain that
%
%e3.9 #&#
\begin{equation}
\frac{ 1 }{ \sqrt{ 2 \pi} } \int_{
{
([
\ln( \pi/ ( 2 h ) )
]^{ 1 / 3 }
- x)
}/{
\sqrt{ t }
}
}^{ \infty} \cos\bigl( h
\cdot\exp\bigl( [ x + \sqrt{ t } y ]^3 \bigr) \bigr)
e^{ { - y^2 }/{ 2 } } \,dy \leq0
\end{equation}
for all
$
t \in(0,\infty)
$,
$
h
\in
( 0,
\frac{ \pi}{ 2 }
\exp( - | x \vee0 |^3 )
]
$
and all
$ x \in\R$.
This implies
%
%e3.10 #&#
\begin{eqnarray}
\label{eqestcrucial0Hoelder} && \mathbb{E} \bigl[ \cos\bigl( h
\cdot\exp\bigl( \bigl[ x +
W(t) \bigr]^3 \bigr) \bigr) \bigr] \nonumber\\
&&\qquad= \frac{ 1 }{ \sqrt{ 2 \pi} }
\int
_{ - \infty}^{ \infty} \cos\bigl( h \cdot\exp\bigl( [ x +
\sqrt{ t } y ]^3 \bigr) \bigr) e^{ { - y^2 }/{ 2 } } \,dy
\nonumber\\
&&\qquad \leq\frac{ 1 }{ \sqrt{ 2 \pi} } \int_{ - \infty}^{
{
([
\ln( \pi/ ( 2 h ) )
]^{ 1 / 3 }
- x)
}/{
\sqrt{ t }
}
} \cos
\bigl( h \cdot\exp\bigl( [ x + \sqrt{ t } y ]^3 \bigr) \bigr)
e^{ { - y^2 }/{ 2 } } \,dy
\\
&&\qquad \leq\mathbb{P} \biggl[ W_1 \leq\frac{
[
\ln( \pi/ ( 2 h ) )
]^{ 1 / 3 }
- x
}{
\sqrt{ t }
} \biggr] \nonumber\\
&&\qquad= 1 -
\mathbb{P} \biggl[ W_1 > \frac{
[
\ln( \pi/ ( 2 h ) )
]^{ 1 / 3 }
- x
}{
\sqrt{ t }
} \biggr]\nonumber
\end{eqnarray}
for all
$
t \in(0,\infty)
$,
$
h
\in
( 0,
\frac{ \pi}{ 2 }
\exp( - | x \vee0 |^3 )
]
$
and all
$ x \in\R$.
Moreover,
Lemma~22.2 in
Klenke~\cite{Klenke2008}
yields
%
%e3.11 #&#
\begin{equation}
\P[ W_1 > y ] \geq\frac{
e^{ - { y^2 }/{ 2 } }
}{
y \sqrt{ 2 \pi}
( 1 + y^{ - 2 } )
} \geq\frac{
e^{ - { y^2 }/{ 2 } }
}{
y \sqrt{ 8 \pi}
} \geq
e^{ - 4 y^2 }
\end{equation}
for all $ y \in[1,\infty) $.
Combining this and
inequality~\eqref{eqestcrucial0Hoelder}
then shows
%
%e3.12 #&#
\begin{eqnarray}
 1 - \mathbb{E} \bigl[ \cos\bigl( h \cdot\exp\bigl( \bigl[ x + W(t)
\bigr]^3 \bigr) \bigr) \bigr] &\geq&\mathbb{P} \biggl[ W_1
> \frac{
[
\ln( \pi/ ( 2 h ) )
]^{ 1 / 3 }
- x
}{
\sqrt{ t }
} \biggr]
\nonumber
\\[-8pt]
\\[-8pt]
\nonumber
 % \\ &
&\geq&\exp\biggl( \frac{
- 4
\vert
[
\ln( \pi/ ( 2 h ) )
]^{ 1 / 3 }
- x
\vert^2
}{
t
}
\biggr)
\end{eqnarray}
for all
$
h
\in
( 0,
\frac{ \pi}{ 2 }
\exp (
-
| [ \sqrt{t} + x ] \vee0 |^3
)
]
$,
$ t \in(0,\infty) $
and all
$ x \in\R$
and the estimate
$
- \vert a + b \vert^2
\geq
- 2 a^2 - 2 b^2
$
for all $ a, b \in\R$
therefore results in
the first inequality
in~\eqref{eqlowerboundcos}.
Next, the first inequality
in \eqref{eqlowerboundcos}
implies
%
%e3.13 #&#
\begin{eqnarray}
&& \mathbb{E} \bigl[ \mathbh{1}_{
\{
W(t) \in A
\}
} \bigl\vert1 - \cos\bigl( h
\cdot\exp\bigl( \bigl[ x + W(s) \bigr]^3 \bigr) \bigr) \bigr
\vert
\bigr] %\\ & =
% \mathbb{E}\Big[
% \mathbh{1}_{
% \left\{
% W(t) \in A
% \right\}
% }
% \big|
% 1
% -
% \cos\big(
% h \cdot
% \exp\big(
% \left[
% x
% + \frac{ s }{ t } W(t)
% + W(s) - \frac{ s }{ t } W(t)
% \right]^3
% \big)
% \big)
% \big|
% \Big]
\nonumber\\
&&\qquad = \mathbb{E} \biggl[ \mathbh{1}_{
\{
W(t) \in A
\}
} \E\biggl[ 1 - \cos\biggl( h
\cdot\exp\biggl( \biggl[ x + \frac{ s }{ t } W(t) + W(s)
\nonumber
\\[-8pt]
\\[-8pt]
\nonumber
&&\hspace*{226pt}{}- \frac{
s }{ t }
W(t) \biggr]^3 \biggr) \biggr)\Big | W(t) \biggr] \biggr]
\\
&&\qquad \geq\mathbb{E} \biggl[ \mathbh{1}_{
\{
W(t) \in A
\}
} \exp\biggl(
\frac{
- 8 t
}{ s ( t - s ) } \biggl[ \biggl| \ln\biggl( \frac{ \pi}{ 2 h } \biggr)
\biggr|^{ 2 / 3 } + \biggl[ x + \frac{ s }{ t } W(t) \biggr]^2
\biggr] \biggr) \biggr]\nonumber
\end{eqnarray}
for all
$
h
\in
( 0,
\frac{ \pi}{ 2 }
\exp (
-
[
\sqrt{ t } + | x | +
\sup_{ a \in A } | a |
]^3
)
]
$,
$ x \in\R$,
$ s, t \in(0,\infty) $
with $ s < t $
and all bounded and Borel
measurable sets
$ A \subset\R$.
Hence, we get
\begin{eqnarray}
\label{eqlowerbound002} & &\int_{ 0 }^{ t } \mathbb{E} \bigl[
\mathbh{1}_{
\{
W(t) \in A
\}
} \bigl\vert1 - \cos\bigl( h \cdot\exp\bigl( \bigl[
x + W(s) \bigr]^3 \bigr) \bigr) \bigr\vert\bigr] \,ds
\nonumber\\
& &\qquad\geq\int_{
{ t }/{ 3 }
}^{
{ 2 t }/{ 3 }
} \mathbb{E} \bigl[
\mathbh{1}_{
\{
W(t) \in A
\}
} \bigl\vert1 - \cos\bigl( h \cdot\exp\bigl( \bigl[
x + W(s) \bigr]^3 \bigr) \bigr) \bigr\vert\bigr] \,ds
\nonumber
\\
&&\qquad \geq\int_{
{ t }/{ 3 }
}^{
{ 2 t }/{ 3 }
} \mathbb{E} \biggl[
\mathbh{1}_{
\{
W(t) \in A
\}
} \exp\biggl( \frac{
- 8 t
}{ s ( t - s ) } \biggl[\biggl | \ln\biggl(
\frac{ \pi}{ 2 h } \biggr) \biggr|^{ 2 / 3 } \\
&&\hspace*{195pt}{}+ \biggl[ x + \frac{ s }{ t
} W(t)
\biggr]^2 \biggr] \biggr) \biggr] \,ds
\nonumber\\
&&\qquad \geq\frac{ t }{ 3 } \cdot\E\biggl[ \mathbh{1}_{
\{
W(t) \in A
\}
} \exp
\biggl( \frac{
- 64
}{ t } \biggl[ \biggl| \ln\biggl( \frac{ \pi}{ 2 h } \biggr)
\biggr|^{ 2 / 3 } + x^2 + \bigl\vert W(t) \bigr\vert
^2 \biggr] \biggr) \biggr]\nonumber
\end{eqnarray}
for all
$
h
\in
( 0,
\frac{ \pi}{ 2 }
\exp (
-
[
\sqrt{ t } + | x | +
\sup_{ a \in A } | a |
]^3
)
]
$,
$ x \in\R$,
$ t \in(0,\infty) $
and all bounded
and Borel measurable sets
$ A \subset\R$.
%It thus remains to
%show \eqref{eqlimithepsilon}
%to complete the proof
%of Lemma~\ref{lemex2lowerbound}.
%For this note that
%ensures that
%&
% \liminf_{ h \searrow0 }
% \left(
% \frac{
% \int_0^t
% \mathbb{E}\big[
% \mathbh{1}_{
% \left\{
% W(t) \in A
% \right\}
% }
% \left|
% 1
% -
% \cos\left(
% h \cdot
% \exp\left(
% [ x + W(s) ]^3
% \right)
% \right)
% \right|
% \big]
% \,ds
% }{ h^{ \varepsilon} }
% \right)
% \frac{ t }{ 3 }
% \cdot
% \E\left[
% \mathbh{1}_{
% \left\{
% W(t) \in A
% \right\}
% }
% e^{
% \frac{
% - 72
% }{ t }
% \left[
% x^2
% +
% \left| W(t) \right|^2
% \right]
% }
% \right]
% \cdot
% \liminf_{ h \searrow0 }
% \exp\left(
% \frac{
% -72
% }{ t }
% \left|
% \ln( \frac{ \pi}{ 2 h } )
% \right|^{ 2 / 3 }
% - \varepsilon\ln( h )
% \right)
% \frac{ t }{ 3 }
% \cdot
% \E\left[
% \mathbh{1}_{
% \left\{
% W(t) \in A
% \right\}
% }
% e^{
% \frac{
% - 72
% }{ t }
% \left[
% x^2
% +
% \left| W(t) \right|^2
% \right]
% }
% \right]
% \cdot
% \exp\left(
% \liminf_{ h \searrow0 }
% \left[
% \varepsilon\ln( \frac{ 1 }{ h } )
% -
% \frac{
% 72
% }{ t }
% |
% \ln( \frac{ \pi}{ 2 h } )
% |^{ 2 / 3 }
% \right]
% \right)
% = \infty
%for all $ \varepsilon, t \in(0,\infty) $,
%$ x \in\R$
%and all $ A \in\mathcal{B}( \R) $
%with strictly positive Lebesgue measure.
This completes the proof of
Lem\-ma~\ref{lemex2lowerbound}.
\end{pf}

We are now ready to prove
Theorem~\ref{corex2bfinalcorollary}
stated at the beginning of this section.
Its proof uses the lower
bound established
in Lemma~\ref{lemex2lowerbound}
above.

\begin{pf*}{Proof
of
Theorem~\ref{corex2bfinalcorollary}}
First of all,
note that
\eqref{eqexactrepsec3}
and
Lemma~\ref{lemex2lowerbound}
imply
that
%
%e3.14 #&#
\begin{eqnarray}
 &&\lim_{ h \searrow0 } \biggl( \frac{
\E[
X^{ (0,0,0) }_1(t)
-
X^{ (0,0,h) }_1(t)
]
}{
h^{ \varepsilon}
} \biggr) \nonumber\\
&&\qquad= \lim
_{ h \searrow0 } \biggl( \frac{
\E[
\int_0^t
1
-
\cos (
h \cdot
\exp (
[ \sqrt{2} W(s) ]^3
)
)
\,ds
]
}{
h^{ \varepsilon}
} \biggr)
\nonumber\\
&&\qquad = \lim_{ h \searrow0 } \biggl( \frac{
\int_0^t
1
-
\E[
\cos (
h \cdot
\exp (
[ W(2 s) ]^3
)
)
]
\,ds
}{
h^{ \varepsilon}
} \biggr)\nonumber\\
&&\qquad = \lim
_{ h \searrow0 } \biggl( \frac{
\int_0^{ 2 t }
1
-
\E[
\cos (
h \cdot
\exp (
[ W(s) ]^3
)
)
]
\,ds
}{
2 h^{ \varepsilon}
} \biggr)
\nonumber
\\[-8pt]
\\[-8pt]
\nonumber
&&\qquad \geq\lim_{ h \searrow0 } \biggl( \frac{
\int_t^{ 2 t }
1
-
\E[
\cos (
h \cdot
\exp (
[ W(s) ]^3
)
)
]
\,ds
}{
2 h^{ \varepsilon}
} \biggr) \\
&&\qquad\geq
\lim_{ h \searrow0 } \biggl( \frac{
\int_t^{ 2 t }
\exp (
({- 8
}/{ t })
\vert
\ln( { \pi}/{( 2 h) } )
\vert^{ 2 / 3 }
)
\,ds
}{
2 h^{ \varepsilon}
} \biggr)
\nonumber\\
&&\qquad = \lim_{ h \searrow0 } \biggl( \frac{ t }{ 2 } \cdot\exp\biggl(
\frac{
- 8
}{ t } \biggl\vert\ln\biggl( \frac{ \pi}{ 2 h } \biggr) \biggr
\vert^{ 2 / 3 } + \ln\bigl( h^{ - \varepsilon} \bigr) \biggr)
\biggr)\nonumber\\
&&\qquad =
\frac{ t }{ 2 } \cdot\lim_{ h \searrow0 } \biggl( \exp\biggl(
\frac{
- 8
}{ t } \biggl\vert\ln\biggl( \frac{ \pi}{ 2 h } \biggr) \biggr
\vert^{ 2 / 3 } - \varepsilon\cdot\ln( h ) \biggr) \biggr) =
\infty\nonumber
\end{eqnarray}
for all
$ \varepsilon, t \in(0,\infty) $.
We hence get for every
$ t \in(0,\infty) $ that
the function
$
\R^3 \ni x \mapsto
\mathbb{E} [
X^x(t)
]
\in\R^3
$
is not locally
H\"{o}lder continuous.
Moreover, let
$ \psi\colonn\R\to[0,1] $
be an infinitely often differentiable
function with compact support
and with
$ \psi(x) = 1 $
for all
$ x \in[ -T, T] $
and let
$ \varphi\colonn\R^3 \to\R$
be a function given by
$
\varphi( x_1, x_2, x_3 )
=
x_1 \psi( x_1 ) \psi( x_2 )
\psi( x_3 )
$
for all $ x_1, x_2, x_3 \in\R$.
Again
\eqref{eqexactrepsec3}
and Lemma~\ref{lemex2lowerbound}
then show
%
%e3.15 #&#
\begin{eqnarray}
&& \lim_{ h \searrow0 } \bigl( h^{ - \varepsilon} \cdot\E\bigl[
\varphi
\bigl( X^{ (0,0,0) }(t) \bigr) - \varphi\bigl( X^{ (0,0,h) }(t)
\bigr)
\bigr] \bigr)\nonumber\\
&&\qquad = \lim_{ h \searrow0 } \bigl( h^{ - \varepsilon}
\cdot\E
\bigl[ \bigl( X^{ (0,0,0) }_1(t) - X^{ (0,0,h) }_1(t)
\bigr) \psi\bigl( \sqrt{2} W(t) \bigr) \bigr] \bigr)
\nonumber\\
&&\qquad \geq\lim_{ h \searrow0 } \bigl( h^{ - \varepsilon} \cdot\E
\bigl[
\mathbh{1}_{
\{
\vert \sqrt{2} W(t) \vert \leq T
\}
} \bigl( X^{ (0,0,0) }_1(t) -
X^{ (0,0,h) }_1(t) \bigr) \bigr] \bigr)
\\
&&\qquad = \lim_{ h \searrow0 } \biggl( h^{ - \varepsilon} \cdot\E\biggl
[ \int
_0^t \mathbh{1}_{
\{
\vert \sqrt{2} W(t) \vert \leq T
\}
} \bigl( 1 -
\cos\bigl( h \cdot\exp\bigl( \bigl[ \sqrt{ 2 } W(s) \bigr]^3
\bigr)
\bigr) \bigr) \,ds \biggr] \biggr)
\nonumber\\
&&\qquad = \lim_{ h \searrow0 } \biggl( \frac{ 1 }{ 2 h^{ \varepsilon} }
\cdot\E\biggl[
\int_0^{ 2 t } \mathbh{1}_{
\{
\vert W( 2 t ) \vert \leq T
\}
} \bigl( 1 -
\cos\bigl( h \cdot\exp\bigl( \bigl[ W(s) \bigr]^3 \bigr) \bigr)
\bigr) \,ds \biggr] \biggr) = \infty\nonumber
\end{eqnarray}
for all $ t \in(0,T] $.
The proof
of Theorem~\ref{corex2bfinalcorollary}
is thus completed.
\end{pf*}

In the remainder of this section,
we briefly consider a slightly
modified version of
the SDE~\eqref{eqex2bSDE}.
More formally, let
$
(
\Z_n
)_{
n \in\N_0
}
$
be a family
of sets
defined by
$
\mathbb{Z}_0:= \mathbb{Z}
:= \{ \ldots, -2, -1, 0, 1, 2, \ldots\}
$
and by
$
\mathbb{Z}_n
:=
\{
z \in\mathbb{Z} \colonn
\frac{z}{2} \notin\mathbb{Z}
\}
$
$
=
\{ \ldots, -3, -1, 1, 3, \ldots
\}
$
for all $ n \in\mathbb{N} $.
Then let
$
\mu= (\mu_1, \mu_2, \mu_3)
\colonn\mathbb{R}^3
\rightarrow\mathbb{R}^3
$
and
$
B \in\mathbb{R}^{ 3 }
$
be given by
%
%e3.16 #&#
\begin{eqnarray}
\label{eqex2defmu} \mu(x) &=& \pmatrix{
 \displaystyle\sum
_{ n = 0 }^{ \infty} \sum_{ m \in\mathbb{Z}_n }
\frac{ 1 }{ 4^{ (n + |m|) } } \cos\biggl( \biggl( x_3 - \frac{ m
}{ 2^n }
\biggr) \exp\bigl( [ x_2 ]^3 \bigr) \biggr)
\vspace*{2pt}\cr
0
\vspace*{2pt}\cr
0
}
\quad \mbox{and}
\nonumber
\\[-8pt]
\\[-8pt]
\nonumber
 %\end{equation}
%and by
B &=& \pmatrix{ 0
\vspace*{2pt}\cr
1
\vspace*{2pt}\cr
0
}
\end{eqnarray}
for all
$
x = (x_1,x_2,x_3) \in
\mathbb{R}^3
$.
Note that
$ \mu\colonn\R^3 \to\R^3 $
is infinitely often differentiable
and
globally bounded by~$ 2 $.
%Indeed, we have
%%&
% \sup_{ x \in\R^3 }
% \left\|
% \mu( x )
% \right\|
% & \leq
% \sum_{ n = 0 }^{ \infty}
% \sum_{ m \in\mathbb{Z}_n }
% \frac{ 1 }{ 4^{ (n + |m|) } }
% \left[
% \sup_{ x_2, x_3 \in\R}
% \left|
% \cos\Big(
% ( x_3 - \frac{ m }{ 2^n } )
% \exp\left( [ x_2 ]^3 \right)
% \Big)
% \right|
% \right]
% \sum_{ n = 0 }^{ \infty}
% \sum_{ m \in\mathbb{Z}_n }
% \frac{ 1 }{ 4^{ (n + |m|) } }
% =
% \sum_{ m \in\mathbb{Z} }
% \frac{ 1 }{ 4^{ |m| } }
% +
% \left(
% \sum_{ n = 1 }^{ \infty}
% \frac{ 1 }{ 4^n }
% \right)
% \left(
% \sum_{ m \in\mathbb{Z}_1 }
% \frac{ 1 }{ 4^{ |m| } }
% \right)
% \leq2.
Moreover, let
$
( \Omega, \mathcal{F}, \P
)
$
be a probability space,
let
$
W \colonn[0,\infty) \times\Omega
\to\R
$
be a one-dimensional
standard Brownian
motion and let
$
X^x \colonn[0,\infty) \times
\Omega\to\R^3
$,
$ x \in\R^3 $,
be the up to indistinguishability
unique
solution processes
with continuous sample paths
of the SDE
%
%e3.17 #&#
\begin{equation}
\label{eqSDESec3B} d X^x(t) = \mu\bigl( X^x(t) \bigr) \,dt + B
\,d W(t)
\end{equation}
for $ t \in[0,\infty) $
and $ x \in\R^3 $
satisfying
$ X^x(0) = x $
for all $ x \in\R^3 $.
%The stochastic
%processes
%$
% X^x = (X^x_1, X^x_2, X^x_3)
% \colonn[0,\infty) \times
% \Omega\to\R^3
%$,
%$ x \in\R^3 $,
%are thus solution
%processes of the SDE
% d X^x_1(t)
%& =
% \sum_{ n = 0 }^{ \infty}
% \sum_{ m \in\mathbb{Z}_n }
% \frac{
% 1
% }{ 4^{ (n + |m|) } }
% \cos\left(
% ( X^x_3(t) - \frac{ m }{ 2^n } )
% \exp\big( [ X^x_2(t) ]^3 \big)
% \right)
% dt
% +
% 0 dW_1(t)
% d X^x_2(t)
%& =
% 0 dt + dW_2(t)
% d X^x_3(t) & = 0 dt
% + 0 dW_3(t)
%for $ t \in[0,\infty) $
%and $ x \in\R^3 $
%where
%$ X^x(0) = x $
%for all $ x \in\R^3 $.
The following Theorem~\ref{thmirr2}
establishes
that the function
$
[0,\infty) \times\R^3
\to
\E[
X^x(t)
]
\in\R^3
$
is nowhere locally
H\"{o}lder continuous.
Its proof is a
straightforward
consequence of
Lemma~\ref{lemex2lowerbound}
and, therefore, omitted.

%th3.4 #&#
\begin{theorem}[(A further
counterexample to regularity
preservation with degenerate
additive noise)]
\label{thmirr2}
Let
$ c, T \in(0,\infty) $
and
let
$
X^x \colonn[0,\infty) \times\Omega
\to\R^3
$,
$ x \in\R^3 $,
be solution processes
of the SDE~\eqref{eqSDESec3B}
with continuous sample paths
and with
$
X^x(0) = x
$
for all $ x \in\R^3 $.
Then
for every $ t \in(0,\infty) $
and
every nonempty open
set $ O \subset\R^3 $,
the function
$
O \ni x \mapsto
\mathbb{E} [
X^x(t)
]
\in\R^3
$
is
continuous but
not locally
H\"{o}lder continuous.
Moreover,
there exists
an infinitely often differentiable
function
$
\varphi\in
\C^{ \infty}_{ \mathrm{cpt} }( \R^3, \R)
$
with compact support
such that
for every $ t \in(0,T] $
and
every nonempty open set
$ O \subset( -c, c )^3 $
the function
$
O \ni x \mapsto
\mathbb{E} [
\varphi( X^x(t) )
]
\in\R
$
is
continuous but
not locally
H\"{o}lder continuous.
\end{theorem}

%s4 #&#
\section{Solutions of Kolmogorov
equations}
\label{secKolmogorovequations}

If the transition
semigroup associated with an SDE is smooth,
then it satisfies the
Kolmogorov
equation
(which is a second-order linear PDE)
corresponding to the SDE
in the classical sense.
The transition semigroups in
our counterexamples are,
however, not locally Lipschitz
continuous and are therefore
no classical solutions of the
Kolmogorov equations
of the corresponding SDEs.
The purpose of this section
is to verify that the
nonsmooth transition semigroup
associated with
such an SDE still satisfies the
Kolmogorov equation
but in a certain weak sense.
More precisely,
in Section~\ref{ssecViscositysolutionsofKolmogorovequations},
we show that the
transition semigroups in
our counterexamples
are \textit{viscosity solutions}
%(see Definition~\ref{dviscositysolution})
of the associated Kolmogorov equations.
Moreover,
in Section~\ref{ssecSolutionsofKolmogorovequationsinthedistributionalsense},
we show that the
transition semigroups in our
counterexamples are solutions
of the associated Kolmogorov equations
\textit{in the distributional sense}.
Throughout this section, the notation
$
\sup( \varnothing):= - \infty
$
and
$
\inf(\varnothing):= \infty
$
is used.

%%%%%%%%%%%%%%%%%%%%%%%%%%%%
%s4.1 #&#
\subsection{Definition and basic properties of viscosity solutions}
\label{ssecDefinitionofviscositysolutions}

Viscosity solutions were first introduced
in Crandall and Lions~\cite{CrandallLions1983}
(see also~\cite{Evans1978,Evans1980,CrandallLions1981}).
The name \emph{viscosity solution} is due to the method
of vanishing viscosity; see the
discussion in Section~10.1 in
Evans~\cite{Evans2010}.
For a review of the theory and
for more references,
we refer the reader to
the well-known users's
guide
Crandall, Ishii and Lions~\cite{CrandallIshiiLions1992}.

For $ d \in\N$,
we denote by
$
\mathbb{S}_d
=
\{
A \in\R^{ d \times d }
\colonn A = A^*
\}
$
the set of all
symmetric
$ d \times d $-matrices.
Moreover,
for $ d \in\N$
and $ A, B \in\mathbb{S}_d $
we write
$ A \leq B $ in the following
if
$
\langle x, A x \rangle
\leq
\langle x, B x \rangle
$
for all
$ x \in\R^d $.
Furthermore,
for $ d \in\N$
and an open set
$ O \subset\R^d $
we call
a function
$
F
\colonn
O
\times
\R
\times
\R^d
\times
\mathbb{S}_d
\to\R
$
%is called \emph{proper} if
% F(x,r,p,A)
% \leq
% F(x,s,p,B)
%for all
%$ x \in O $,
%$ p \in\R^d $,
%$ r, s \in\R$
%with $ r \leq s $
%and all
%$ A,B \in\mathbb{S}_d $ with $A\geq B$
%and $ F $
%is called
\emph{degenerate elliptic}
(see, e.g., (0.3) in Crandall,
Ishii
and Lions~\cite{CrandallIshiiLions1992})
if
$
F(x, r, p, A) \leq F(x, r, p, B)
$
for all
$ x \in O $,
$ r \in\R$,
$ p \in\R^d $
and all
$
A, B \in\mathbb{S}_d
$
with
$ A \geq B $.
For convenience
of the reader, we recall
the definition of a viscosity
solution
(see, e.g., Section~2 in Crandall,
Ishii
and Lions~\cite{CrandallIshiiLions1992}
and also Definition~1.2
in Appendix~C
in Peng~\cite{Peng2010}).

%

% TODO: User Guide always assumes $F$ to be continuous,
% see line 2 of Section 2.
% Also Peng 2010 assumes continuity. Ishii 1989 does not.

%de4.1 #&#
\begin{definition}[(Viscosity solution)]
\label{dviscositysolution}
Let $ d \in\N$, let $ O \subset\R^d $
be an open set
and
let
$
F \colonn
O
\times\R\times\R^d
\times\mathbb{S}_d
\to\R
$
be a degenerate elliptic function.
A
function $ u \colonn O \to\R$
is
said to be
a \emph{viscosity subsolution} of $F=0$
(or, equivalently, a viscosity solution
of $ F \leq0 $) if
$ u $ is upper semicontinuous and if
it holds
for all $ x \in O $
and all
$
\phi\in\C^2(O,\R)
$
with
$ \phi\geq u $
and
$ \phi(x) = u(x) $
that
%
%e4.1 #&#
\begin{equation}
F \bigl( x, \phi(x), (\nabla\phi) (x), (\operatorname{Hess} \phi)
(x) \bigr)
\leq0.
\end{equation}
Similarly,
a
function $ u \colonn O \to\R$
is said to be
a \emph{viscosity supersolution}
of $ F = 0 $ (or, equivalently,
a viscosity solution
of $F \geq0$)
if
$ u $ is lower semicontinuous
and if
it holds
for all $ x \in O $
and all
$
\phi\in\C^2(O,\R)
$
with
$ \phi\leq u $
and
$ \phi(x) = u(x) $
that
%
%e4.2 #&#
\begin{equation}
\label{eqFleq0} F \bigl( x, \phi(x), (\nabla\phi) (x),
(\operatorname{Hess}
\phi) (x) \bigr) \geq0.\vadjust{\goodbreak}
\end{equation}
Finally, a
function $ u \colonn O\to\R$
is said to be
a \emph{viscosity solution} of $F=0$
if $u$ is both a viscosity subsolution
and a viscosity supersolution of
$F=0$.
\end{definition}

In the proof of
Corollary~\ref{ctechnicallemmauniqueness}
below,
% Below in
% Subsection~\ref{ssecUniquenessof%viscositysolutionsof%Kolmogorovequations}
the following elementary lemma (Lem\-ma~\ref{lemsignchanges}) is used.
The proof of Lemma~\ref{lemsignchanges} is clear and, therefore, omitted.

%le4.2 #&#
\begin{lemma}[(Sign changes of viscosity solutions)]
\label{lemsignchanges}
Let $ d \in\N$, let $ O \subset\R^d $ be an open set,
let
$
F \colonn O \times\R\times\R^d \times\mathbb{S}_d \to\R
$
be a degenerate elliptic function and let
$ u \colonn O \to\R$ be a viscosity solution of
$ F \geq0 $.
Then the function
$
\tilde{F} \colonn O \times\R\times\R^d \times\mathbb{S}_d \to\R
$
defined by
$
\tilde{F}( x, r, p, A )
:= - F( x, - r, - p, - A )
$
for all
$ ( x, r, p, A ) \in O \times\R\times\R^d \times\mathbb{S}_d $
is degenerate elliptic and
the function
$
O \ni x \mapsto- u(x) \in\R
$
is a viscosity solution of $ \tilde{F} \leq0 $.
The corresponding statement holds for viscosity solutions
of $ F \leq0 $ and $ F = 0$,
respectively.
\end{lemma}

Above in Definition~\ref{dviscositysolution},
the concept of viscosity solutions is presented
via suitable test functions. An alternative instrument to characterize
viscosity solutions are so-called \emph{semijets}
(see, e.g., Definition~2.2 in
Crandall, Ishii and Lions~\cite{CrandallIshiiLions1992}).
They are recalled in the next definition.

%de4.3 #&#
\begin{definition}[(Second-order semijets)]
\label{dsemijets}
Let $ d \in\N$, let $ O \subset\R^d $ be an open set
and let $ u \colonn O \to\R$ be a function. Then
we define functions
$
J^2_{ + } u \colonn O \to
\mathcal{P} ( \R^d \times\mathbb{S}_d )
$,
$
J^2_{ - } u \colonn O \to
\mathcal{P} ( \R^d \times\mathbb{S}_d )
$,
$
\hat{J}^2_{ + } u \colonn O \to
\mathcal{P} ( \R^d \times\mathbb{S}_d )
$
and
$
\hat{J}^2_{ - } u \colonn O \to
\mathcal{P} ( \R^d \times\mathbb{S}_d )
$
by
\begin{eqnarray*}
&&\hspace*{-4pt}\bigl( J^2_{ + } u \bigr) ( x ) \\
&&\hspace*{-4pt}\qquad:= \biggl\{ ( p, A ) \in
\R^d \times\mathbb{S}_d \colonn\\
&&\hspace*{38pt}\limsup
_{
O \setminus\{ x \}
\ni y \to x
} \biggl( \frac{
u( y ) - u( x )
-
\anglel p, x - y \angler
-
({ 1 }/{ 2 })
\anglel x - y, A ( x - y ) \angler
}{
\llVert x - y \rrVert^2
} \biggr) \leq0 \biggr\}\hspace*{-0.5pt},
\\
&&\hspace*{-4pt}\bigl( \hat{J}^2_{ + } u \bigr) ( x ) \\
&&\hspace*{-8pt}\qquad:= \Biggl\{
(p, A) \in\R^d \times\mathbb{S}_d \colonn\\
&&\hspace*{34pt}\pmatrix{
\exists( x_n, p_n, A_n
)_{ n \in\N} \subset O \times\R^d \times\mathbb{S}_d
\colonn\bigl( \forall n \in\N\colonn(p_n, A_n) \in
\bigl( J^2_{ + } u\bigr) ( x_n ) \bigr)
\vspace*{2pt}\cr
\mbox{and }\displaystyle \lim_{ n \to\infty} \bigl( x_n,
u(x_n), p_n, A_n \bigr) = \bigl( x, u(x), p,
A \bigr)}\hspace*{-0.5pt}
 \Biggr\}\hspace*{-0.5pt},
\\
&&\hspace*{-4pt}\bigl( J^2_{ - } u \bigr) ( x ) \\
&&\hspace*{-4pt}\qquad:= \biggl\{ ( p, A ) \in
\R^d \times\mathbb{S}_d \colonn\\
&&\hspace*{38pt}\liminf
_{
O \setminus\{ x \}
\ni y \to x
} \biggl( \frac{
u( y ) - u( x )
-
\anglel p, x - y \angler
-
({ 1 }/{ 2 })
\anglel x - y, A ( x - y ) \angler
}{
\llVert x - y \rrVert^2
} \biggr) \geq0 \biggr\}
\end{eqnarray*}
and
\begin{eqnarray*}
&&\hspace*{-4pt}\bigl( \hat{J}^2_{ - } u \bigr) ( x ) \\
&&\hspace*{-4pt}\qquad:= \Biggl\{
(p, A) \in\R^d \times\mathbb{S}_d \colonn\\
&&\hspace*{38pt}\pmatrix{
\exists( x_n, p_n, A_n
)_{ n \in\N} \subset O \times\R^d \times\mathbb{S}_d
\colonn\bigl( \forall n \in\N\colonn(p_n, A_n) \in
\bigl( J^2_{ - } u\bigr) ( x_n ) \bigr)
\vspace*{2pt}\cr
\mbox{and } \displaystyle\lim_{ n \to\infty} \bigl( x_n,
u(x_n), p_n, A_n \bigr) = \bigl( x, u(x), p,
A \bigr)} \Biggr\}
\end{eqnarray*}
for all $ x \in O $.
\end{definition}

The next lemma (Lemma~\ref{lemsemijets}),
which is essentially one of the statements in
Remark~2.3 in
Crandall, Ishii and Lions~\cite{CrandallIshiiLions1992},
illustrates the relationship between semijets in the sense
of Definition~\ref{dsemijets} and
suitable test functions in the sense of Definition~\ref{dviscositysolution}.

%le4.4 #&#
\begin{lemma}[(Properties of semijets)]
\label{lemsemijets}
Let $ d \in\N$, let $ O \subset\R^d $
be an open set and let
$
u \colonn O \to\R
$
be a function. Then
%
%e4.3 #&#
%e4.4 #&#
\begin{eqnarray}
\bigl( J^2_{ + } u\bigr) ( x ) & =& \bigl\{ \bigl( (\nabla
\phi) (x), ( \operatorname{Hess} \phi) (x) \bigr) \in\R^d \times
\mathbb{S}_d \colonn\nonumber\\
&&\hspace*{6pt}\bigl( \phi\in C^2( O, \R) \mbox{
with } u(x) = \phi(x) \mbox{ and } u \leq\phi\bigr) \bigr\}
\nonumber
\\[-8pt]
\\[-8pt]
\nonumber
& = &\bigl\{ \bigl( (\nabla\phi) (x), ( \operatorname{Hess} \phi)
(x) \bigr) \in
\R^d \times\mathbb{S}_d \colonn\\
&&\hspace*{6pt}\bigl( \phi\in
C^2( O, \R) \mbox{ and } u - \phi\mbox{ has a local maximum at } x
\bigr) \bigr\}\nonumber
\end{eqnarray}
and
%
%e4.5 #&#
%e4.6 #&#
\begin{eqnarray}
\bigl( J^2_{ - } u\bigr) ( x ) & = &\bigl\{ \bigl( (\nabla
\phi) (x), ( \operatorname{Hess} \phi) (x) \bigr) \in\R^d \times
\mathbb{S}_d \colonn\nonumber\\
&&\hspace*{6pt}\bigl( \phi\in C^2( O, \R) \mbox{
with } u(x) = \phi(x) \mbox{ and } u \geq\phi\bigr) \bigr\}
\nonumber
\\[-8pt]
\\[-8pt]
\nonumber
& =& \bigl\{ \bigl( (\nabla\phi) (x), ( \operatorname{Hess} \phi)
(x) \bigr) \in
\R^d \times\mathbb{S}_d \colonn\\
&&\hspace*{6pt}\bigl( \phi\in
C^2( O, \R) \mbox{ and } u - \phi\mbox{ has a local minimum at } x
\bigr) \bigr\}\nonumber
\end{eqnarray}
for all $ x \in O $.
\end{lemma}

The next corollary, which is essentially one of the statements
in Remark~2.3 in
Crandall, Ishii and~Lions~\cite{CrandallIshiiLions1992},
is an immediate consequence of Lemma~\ref{lemsemijets} above.

%co4.5 #&#
\begin{corollary}[(Characterizations of viscosity solutions)]
\label{corsemijetsequivalence}
Let
$ d \in\N$,
let
$ O \subset\R^d $
be an open set,
let
$
F \colonn O \times\R\times\R^d \times\mathbb{S}_d \to\R
$
be a degenerate elliptic function
and let $ u \colonn O \to\R$ be an upper semicontinuous function.
Then the following three assertions are equivalent:
\begin{itemize}
\item
$ u $ is a viscosity subsolution of $ F = 0 $
($ u $ is a viscosity solution of $ F \leq0 $),
\item
for every $ x \in O $
and every
$
\phi
\in
\{
\psi\in\C^2( O, \R) \colonn
x \mbox{ is a local maximum of }
( u - \psi) \colonn O \to\R
\}
$
it holds that
$
F (
x, u(x), (\nabla\phi)(x),
(\operatorname{Hess} \phi)(x)
)
\leq0
$,
\item
for every $ x \in O $
and every $ (p,A) \in( J^2_{ + } u)( x ) $
it holds that
$
F( x, u(x), p, A ) \leq0
$.
\end{itemize}
The corresponding statement holds for
viscosity supersolutions and
viscosity solutions.
\end{corollary}

The next corollary, which is
Remark~2.4
in Crandall, Ishii and~Lions~\cite{CrandallIshiiLions1992},
illustrates a further characterization of viscosity solutions
under the assumption that $ F $ is continuous.
It follows immediately from Corollary~\ref{corsemijetsequivalence}
and from the semicontinuity of $ F $.

% If the function $ F $ in Corollary~\ref{corsemijetsequivalence}
% is semicontinuous, then the next corollary
% (Corollary~\ref{corsemijetsequivalence2}) illustrates a further
%characterization
% of viscosity solutions.
% Corollary~\ref{corsemijetsequivalence2} is
% Remark~2.4
% in Crandall, Ishii \citationand~Lions~\cite{CrandallIshiiLions1992}.

%co4.6 #&#
\begin{corollary}[(Characterizations of viscosity solutions
for semicontinuous left-hand sides)]
\label{corsemijetsequivalence2}
Let
$ d \in\N$,
let
$ O \subset\R^d $
be an open set,
let
$
F \colonn O \times\R\times\R^d \times\mathbb{S}_d \to\R
$
be a degenerate elliptic and lower semicontinuous function
and let $ u \colonn O \to\R$ be an upper semicontinuous function.
Then
$ u $ is a viscosity solution of $ F \leq0 $
if and only if it holds
for every $ x \in O $
and every $ (p,A) \in( \hat{J}^2_{ + } u)( x ) $
that
$
F( x, u(x), p, A ) \leq0
$.
The corresponding statement holds for
viscosity solutions of $ F \geq0 $ and $ F = 0$,
respectively.
\end{corollary}

The next well-known remark
(see, e.g., Section~2 in Crandall, Ishii and Lions~\cite
{CrandallIshiiLions1992})
illustrates that
classical solutions are viscosity solutions.
We will use it in the proof of
Lemma~\ref{lemKolmogorovviscosity2} below.

%re4.1 #&#
\begin{remark}[(Classical
solutions are viscosity
solutions)]
\label{remclassicalsolution}
Let $ d \in\N$, let $ O \subset\R^d $
be an open set,
let
$
F \colonn
O
\times\R\times\R^d
\times\mathbb{S}_d
\to\R
$
be a degenerate elliptic function
and let
$ u \in\C^2( O, \R) $
be a classical subsolution
of $ F = 0 $, that is, suppose that
%
%e4.7 #&#
\begin{equation}
\label{eqclassical} F \bigl( x, u(x), (\nabla u) (x), (\operatorname
{Hess} u) (x)
\bigr) \leq0
\end{equation}
for all $ x \in O $.
Then $ u $
is also
a viscosity subsolution
of $ F = 0 $.
Indeed,
for every $ x \in O $
and every
$
\phi
\in
\{
\psi\in\C^2( O, \R) \colonn
x \mbox{ is a local maximum of }
( u - \psi) \colonn O \to\R
\}
$
it holds that
$ ( \nabla( u - \phi) )(x) = 0 $
and
$ (\operatorname{Hess} ( u - \phi) )(x) \leq0$
and, therefore,
%
%e4.8 #&#
\begin{eqnarray}
F \bigl( x, u(x), (\nabla\phi) (x), (\operatorname{Hess} \phi) (x)
\bigr) &= &F
\bigl( x, u(x), (\nabla u) (x), (\operatorname{Hess} \phi) (x) \bigr)
\nonumber
\\
&\leq& F \bigl( x, u(x), (\nabla u) (x), (\operatorname{Hess} u) (x)
\bigr)\\
& \leq&
0\nonumber
\end{eqnarray}
due to \eqref{eqclassical}
and due to the degenerate ellipticity
assumption on $ F $.
The corresponding statement holds
for classical supersolutions
and classical solutions of $ F = 0 $.
\end{remark}

For the convenience of the reader, we also state a special
case of Theorem 3.2 in
Crandall, Ishii and~Lions~\cite{CrandallIshiiLions1992}
in the next lemma.
It
%Lemma~\ref{lTheorem32UserGuidy}
will be used in
the proof of
Lemma~\ref{ltechnicallemmauniqueness} below.

%le4.7 #&#
\begin{lemma}[(Construction of suitable semijets)]
\label{lTheorem32UserGuidy}
Let $ d, k \in\N$,
$ \eps\in(0, \infty) $,
let
$\CO\subset\R^{ d } $
be an open set,
let
$
\Phi
\in C^2(\CO^k, \R)
$,
let
$ u_i \colonn\CO\to\R$,
$ i \in\{ 1, \ldots, k \} $,
be upper semicontinuous functions
and let
$
\hat{x} =
( \hat{x}_1, \ldots, \hat{x}_k )
\in
\CO^k
$
be a local maximum point of the
function
$
\CO^k
\ni( x_1, \ldots, x_k)
\mapsto
(
\sum_{ i = 1 }^k
u_i(x_i)
)
- \Phi(x_1, \ldots, x_k )
\in\R
$.
Then there
exist
matrices
$
A_1 \in\mathbb{S}_{d}
,
 \ldots,
A_k \in\mathbb{S}_{d}
$
% $
% A_i
% \in
% \mathbb{S}_{ d_i }
% $,
% $ i \in\{ 1, \dots, k \} $,
such that
for all $ i \in\{1, \ldots, k \} $
it holds that
$
(
( \nabla_{ x_i } \Phi)( \hat{x} ), A_i
)
\in
(
\hat{J}^2_{ + } u_i
)( \hat{x}_i )
$
and such that
%
%e4.9 #&#
\begin{eqnarray}
- \biggl( \frac{ 1 }{ \eps} + \bigl\| ( \operatorname{Hess} \Phi) (
\hat{x} )
\bigr\|_{
L (
% \times_{ i = 1 }^k
\R^{ kd }
)
} \biggr) I % _{ \R^{ d \times d } }
% \left[
% \sum_{ i = 1 }^k
% \| z_i \|^2
% \right]
&\leq&\pmatrix{ A_1 & \cdots& 0
\vspace*{2pt}\cr
\vdots& \ddots& \vdots
\vspace*{2pt}\cr
0 & \cdots& A_k
}
 % \sum_{ i = 1 }^k
% \left\langle z_i, A_i^{ (0) } z_i \right\rangle
\nonumber
\\[-8pt]
\\[-8pt]
\nonumber
&\leq&% \left\langle
% z,
% \left[
(\operatorname{Hess}
\Phi) (\hat{x}) + \eps\bigl[ ( \operatorname{Hess} \Phi) ( \hat
{x} )
\bigr]^2. % \right] z
% \right\rangle
\end{eqnarray}
%
% \begin{equation}
% -\left(\frac{1}{\eps}
% +\|(\operatorname{Hess} \eta)(\hat{x})\|_{L\left(\times_{i=1}^k
% \right)\sum_{i=1}^k\|z_i\|^2
% \leq\sum_{i=1}^k\left\langle z_i,A_i^{(0)} z_i\right\rangle
% \leq
% (\operatorname{Hess} \eta)(\hat{x})\left(z,z\right)
% +\eps\left\|(\operatorname{Hess} \eta)(\hat{x})(z)\right\|_{L(\R^d,
% \end{equation}
% for all $z=(z_1,\ldots,z_k)\in\times_{i=1}^k\R^{d_i}$.
\end{lemma}

%%%%%%%%%%%%%%%%%%%%%%%%%%%%
%s4.2 #&#
\subsection{An approximation result for viscosity solutions}
\label{ssecAnapproximationresultforviscositysolutions}

The following approximation result
for viscosity solutions is essentially well known
(see Proposition~1.2 in
Ishii~\cite{Ishii1989}
which refers to the first-order case
in Theorem A.2 in
Barles and~Perthame~\cite{BarlesPerthame1987};
see also Lemma 6.1 in Crandall, Ishii and~Lions~\cite{CrandallIshiiLions1992}
and the remarks thereafter).
For completeness, we give the proof here following the line of arguments
for the first-order case in
Theorem~A.2
in Barles
and~Perthame~\cite{BarlesPerthame1987}.
In the remainder of this article, we
use the notation
$
\dist( x, A ):=
\inf (
\{
\| x - y \| \in[0,\infty)
\colonn
y \in A
\}
\cup
\{ \infty\}
)
\in[0,\infty]
$
for all
$ x \in\R^d $,
all
$ A \subset\R^d $
and all $ d \in\N$.

%le4.8 #&#
\begin{lemma}
\label{thmlimitsofviscositysolutions}
Let $ d \in\N$,
let $ O \subset\R^d $
be an open set,
let
$ u_n \colonn O \to\R$,
$ n \in\N_0 $,
be functions
and let
$
F_n \colonn
O \times\R\times\R^d \times
\mathbb{S}_d \to\R
$,
$ n \in\N_0 $, be degenerate elliptic
functions such that
$ F_0 $ is continuous.
Moreover, assume that
%
%e4.10 #&#
\begin{eqnarray}
\label{eqconvergencey}&& \limsup_{ n \to\infty} \sup_{ (x,r,p,A)
\in K }
\bigl\vert F_n(x,r,p,A) - F_0(x,r,p,A) \bigr\vert
\nonumber
\\[-8pt]
\\[-8pt]
\nonumber
&&\qquad= 0
= \limsup_{ n \to\infty} \sup_{ x \in\bar{K} } \bigl\vert
u_n(x) - u_0(x) \bigr\vert
\end{eqnarray}
for all nonempty compact sets
$ K \subset O \times\R\times
\R^d \times\mathbb{S}_d $
and all nonempty compact sets
$ \bar{K} \subset O $
and assume
for every $ n \in\N$ that
$ u_n $ is a viscosity solution
of $ F_n = 0 $.
Then $ u_0 $ is a viscosity solution
of $ F_0 = 0 $.
\end{lemma}

\begin{pf}%{Proof
%of
%Lemma~\ref{thmlimitsofviscositysolutions}}
The proof is divided into
two steps.

\textit{Step}~1:
% In this first step
Let
% assume that
% there exists
% an
$ x_0 \in O $
and
% a function
let
$
\phi\in\C^2(O,\R)
$
be a function
such that
$ x_0 $
is a strict maximum of
$ u_0 - \phi$, that is,
%
%e4.11 #&#
\begin{equation}
\label{eqstrictmaximum} u_0(x) - \phi(x) < u_0(
x_0 ) - \phi( x_0 )
\end{equation}
for all
$
x \in
O \setminus\{ x_0 \}
$.
Then we define
$
r:=
\min ( 1,
\frac{ 1 }{ 2 }
\operatorname{dist}( x_0,
\R^d \setminus O)
) \in[0,1]
$.
Since $ O \subset\R^d $
is an open set,
we obtain that $ r \in(0,1] $.
Furthermore,
continuity of the function $ \phi$
and of the functions
$ u_n $, $ n \in\N$,
together with
compactness
of the set
$
\{
y \in\R^d \colonn
\| y - x_0 \| \leq r
\}
\subset O
$
proves that
there exists a
sequence
$
x_n \in
\{
y \in\R^d \colonn
\| y - x_0 \| \leq r
\}\subset O $,
$ n \in\N$,
of vectors such that
%
%e4.12 #&#
\begin{equation}
\label{eqlocalmaxun} u_n(x) - \phi(x) \leq u_n(x_n)-
\phi(x_n)
\end{equation}
for all
$
x \in\{ y \in\R^d
\colonn\| y - x_0 \| \leq r \}
$
and all $ n \in\N$.
We now prove that
the sequence $ (x_n)_{ n \in\N} $
converges to $ x_0 $.
Aiming at a contraction,
we assume that
the sequence $ (x_n)_{ n \in\N} $
does not converge to
$ x_0 $.
Due to compactness of
$
\{ y \in\R^d \colonn
\| y - x_0 \| \leq r \}
$,
there exists a vector
$
\bar{x}_0 \in
\{
y \in\R^d \colonn0 < \| y - x_0 \| \leq r
\}
\subset O
$
and an increasing sequence
$ n_k \in\N$, $ k \in\N$,
such that
$
\lim_{ k \to\infty} x_{ n_k }
= \bar{x}_0
$.
In particular, we obtain that the set
$
\{ \bar{x}_0 \}
\cup
(
\bigcup_{ k \in\N}
\{ x_{ n_k } \}
)
$
is compact.
Assumption~\eqref{eqconvergencey},
inequality~\eqref{eqlocalmaxun}
and
inequality~\eqref{eqstrictmaximum}
hence imply that
\begin{eqnarray*}
u_0( x_0 ) -
\phi(x_0)
&=&
\lim_{ k \to\infty} \bigl(
u_{ n_k }(x_0) - \phi(x_0)
\bigr)
\leq
\limsup_{ k \to\infty}
\bigl(
u_{n_k}( x_{n_k} ) -
\phi(x_{n_k})
\bigr)
% \\&
\\
&=&
u_0( \bar{x}_0 ) -
\phi( \bar{x}_0 )
<
u_0(x_0) - \phi(x_0).
\end{eqnarray*}
From this contradiction,
we infer that
$ \lim_{ n \to\infty} x_n = x_0 $.
Assumption~\eqref{eqconvergencey}
and
continuity of
$
\nabla\phi\colonn O
\to\R^d
$
and of
$
\operatorname{Hess} \phi\colonn
O \to\mathbb{S}_d
$
hence imply that
%
%e4.13 #&#
\begin{eqnarray}
\label{eqlimit}&& \lim_{ n \to\infty} \bigl( x_n,
u_n(x_n), (\nabla\phi) (x_n), (
\operatorname{Hess} \phi) (x_n) \bigr)
\nonumber
\\[-8pt]
\\[-8pt]
\nonumber
&&\qquad = \bigl( x_0,
u_0(x_0), (\nabla\phi) (x_0), (
\operatorname{Hess} \phi) (x_0) \bigr).
\end{eqnarray}
In addition,
$ \lim_{ n \to\infty} x_n = x_0 $
and \eqref{eqlocalmaxun}
show that
there exists a
natural number $ n_0 \in\N$
such that we have
for all
$ n \in\{n_0, n_0 + 1, \ldots\} $
that
$ \| x_n - x_0 \| < r $
and that
$ x_n \in O $ is a local maximum
of the function
$ ( u_n - \phi) \colonn O \to\R$.
Hence,
Corollary~\ref{corsemijetsequivalence}
and
the assumption that
$ u_n $ is a
viscosity solution of $ F_n = 0 $
% hence
show that
%
%e4.14 #&#
\begin{equation}
\label{eqviscosity0} F_n \bigl( x_n, u_n(x_n),
(\nabla\phi) (x_n), ( \operatorname{Hess} \phi) (x_n)
\bigr) \leq0
\end{equation}
for all $ n \in\{n_0,n_0+1,\ldots\} $.
Continuity of
$ F_0 $,
equation~\eqref{eqlimit},
assumption~\eqref{eqconvergencey},
inequality~\eqref{eqviscosity0}
and
compactness of the set
$
\bigcup_{ n \in\N_0 }
\{
(
x_n, u_n(x_n),
(\nabla\phi)(x_n),\break 
(\operatorname{Hess} \phi)(x_n)
)
\}
$
therefore yield
that
%
%e4.15 #&#
\begin{eqnarray}
&& F_0 \bigl( x_0, u_0(x_0),
(\nabla\phi) (x_0), (\operatorname{Hess} \phi) (x_0)
\bigr)\nonumber\\
&&\qquad = \lim_{ n \to\infty} F_0 \bigl( x_n,
u_n(x_n), (\nabla\phi) (x_n), (
\operatorname{Hess} \phi) (x_n) \bigr)
\\
&&\qquad = \lim_{ n \to\infty} F_n \bigl( x_n,
u_n(x_n), (\nabla\phi) (x_n), (
\operatorname{Hess} \phi) (x_n) \bigr) \leq0.\nonumber
\end{eqnarray}
We thus have proved
that
$
F_0 (
x, u_0(x), (\nabla\phi)(x),
(\operatorname{Hess} \phi)(x)
)
\leq0
$
for all
$
\phi\in
\{
\psi\in\C^2( O, \R)
\colonn
x
$
is a strict maximum of
$
( u_0 - \psi) \colonn O \to\R
\}
$
and all
$ x \in O $.

\textit{Step}~2:
% In this second step
Let
% assume that
% there exists
% an
$ x_0 \in O $
and let
% a function
$
\phi\in\C^2(O,\R)
$
be a function
such that
$
\phi( x_0 ) = u_0( x_0 )
$
and
$ \phi\geq u_0 $.
Next define functions
$
\phi_\eps\colonn O \to\R
$,
$ \eps\in(0,1) $,
by
$
\phi_\eps(x) = \phi(x)
+ \eps\| x - x_0 \|^2
$
for all $x\in O$
and all $\eps\in(0,1)$.
Note for every
$ \eps\in(0,1) $ that
$ x_0 $ is
a strict
maximum of the function
$ ( u_0 - \phi_\eps) \colonn O \to\R$.
Step~1 can thus be applied to obtain
%
%e4.16 #&#
\begin{equation}
\label{eqphiepsinequality} F_0 \bigl( x_0, u_0(x_0),
( \nabla\phi_{ \eps}) (x_0), ( \operatorname{Hess}
\phi_{ \eps}) (x_0) \bigr) \leq0
\end{equation}
for all $ \eps\in(0,1) $.
Moreover, observe that
$
(\nabla\phi_\eps)(x_0)
= (\nabla\phi)(x_0)
$
and that
$
(\operatorname{Hess}
\phi_\eps)(x_0)
=
(\operatorname{Hess} \phi)(x_0)
+ 2 \eps I_{ d }
$
for all $\eps\in(0,1)$
where
$
I_{ d }
\in\mathbb{S}^d
$
is the
$ d \times d $-unit matrix.
Consequently,
we see that
$
\lim_{ \varepsilon\searrow0 }
( \nabla\phi_\eps
)(x_0)
=
( \nabla\phi)(x_0)
$
and that
$
\lim_{ \varepsilon\searrow0 }
( \operatorname{Hess} \phi_\eps)(x_0)
=
( \operatorname{Hess} \phi)(x_0)
$.
Continuity of
$ F_0 $
and
inequality~\eqref{eqphiepsinequality}
hence yield
%
%e4.17 #&#
\begin{eqnarray}
&& F_0 \bigl( x_0, u_0(x_0), (
\nabla\phi) (x_0), (\operatorname{Hess} \phi) (x_0)
\bigr)
\nonumber
\\[-8pt]
\\[-8pt]
\nonumber
&&\qquad = \lim_{ \eps\searrow0 } F_0 \bigl( x_0,
u_0(x_0), ( \nabla\phi_\eps)
(x_0), ( \operatorname{Hess} \phi_\eps)
(x_0) \bigr) \leq0.
\end{eqnarray}
We thus have proved that
$
F_0 (
x, u_0(x), (\nabla\phi)(x),
(\operatorname{Hess} \phi)(x)
)
\leq0
$
for all
$
\phi\in \C^2( O, \R)
$
with $ \phi(x) = u_0(x) $ and
$ \phi\geq u_0 $
and all
$ x \in O $.
This shows that
$ u_0 $ is a viscosity subsolution
of $ F_0 = 0 $.
In the same way, it can be shown
that $ u_0 $ is a viscosity supersolution
of $ F_0 = 0 $ and we thereby obtain
that $ u_0 $ is a viscosity solution
of $ F_0 = 0 $.
The proof
of
Lemma~\ref{thmlimitsofviscositysolutions}
is thus completed.
\end{pf}

%%%%%%%%%%%%%%%%%%%%%%%%%%%%
%s4.3 #&#
\subsection{Uniqueness of viscosity
solutions of Kolmogorov
equations}
\label{ssecUniquenessofviscositysolutionsofKolmogorovequations}

A key result of this
subsection
(Corollary~\ref{coruniqueness2})
%(Lemma~\ref{lcomparisonviscositysolution}
%below)
establishes uniqueness
of viscosity solutions
of a second-order linear PDE within
a certain class of functions
and is apparently new.
This uniqueness result is based
on the well-known concept
of superharmonic functions
or---in the PDE language---on
the idea of dominating
supersolutions.
More precisely,
let $d\in\N$ and let
$
(
\Omega, \mathcal{F},
\mathbb{P}
)
$
be a probability space
with a normal filtration
$
(
\mathcal{F}_t
)_{
t \in[0,\infty)
}
$.
For solution processes
$X^x\colonn[0,\infty)\times\Omega\to\R^d$, $x\in\R^d$, of many SDEs,
there exists a function
$
V \in\C^2 ( \R^d, (0,\infty) )
$
[often
$
\R^d \ni x \mapsto1 + \|x\|^2
\in(0,\infty)
$]
and
a real number
$ \rho\in\R$
such that
the stochastic
processes
$
[0, \infty) \times\Omega\ni
(t, \omega) \to
e^{ - \rho t }
\cdot
V( X^x(t)(\omega) )
\in(0,\infty)
$,
$ x \in\R^d $,
are nonnegative
supermartingales
(so that
$
\E [ V(X^x(t)) ]
\leq
e^{ \rho t } \cdot V(x)
$
for all
$ (t,x) \in[0,\infty) \times\R^d
$);
see, for example, the examples
in Section~4 in
\cite{HutzenthalerJentzen2014Memoires}.
For these stochastic processes
to be supermartingales,
it suffices that the
Lyapunov function $V$ satisfies
%
%e4.18 #&#
\begin{equation}
\CL V (x) \le\rho V(x)
\end{equation}
for all
$
x \in\R^d
$,
where $\CL$ is the generator of the SDE under consideration.
In other words,
it suffices that the map
$
(0,\infty) \times\R^d
\ni(t,x) \to
e^{ \rho t } \cdot V(x)
\in(0,\infty)
$
is a classical supersolution
of the Kolmogorov
equation.
For $ T \in(0,\infty) $,
$ d \in\N$ and an open
set $ O \subset\R^d $,
a function
$
G \colonn
(0,T)
\times
O
\times\R\times\R^d
\times\mathbb{S}_d
\to\R
$
is here called
\textit{degenerate elliptic} if
$
G(t, x, r, p, A) \leq G(t, x, r, p, B)
$
for all
$ t \in(0,T)$,
$ x \in O $,
$ r \in\R$,
$ p \in\R^d $
and all
$ A, B \in\mathbb{S}_d $
with
$ A \leq B $
(see, e.g.,
inequality~(1.2)
in Appendix~C
in Peng~\cite{Peng2010}
and compare also
with Section~\ref{ssecDefinitionofviscositysolutions}
above).
% Lemma~\ref{lcomparisonviscositysolution}
% is our main uniqueness result for viscosity solutions
% of parabolic PDEs.
%This result will then be applied in T O D O
%where the function
% $
% G \colonn
% (0,T)
% \times
% O
% \times\R\times\R^d
% \times\mathbb{S}_d
% \to\R
% $
To establish Corollary~\ref{coruniqueness2},
we first state a few auxiliary results.
For the convenience of the reader,
we first state Proposition 3.7 from
Crandall, Ishii and~Lions~\cite{CrandallIshiiLions1992}
in the next lemma.

%le4.9 #&#
\begin{lemma}
\label{lProp37UserGuide}
Let
$ d \in\N$,
let
$ O \subset\R^d $
be a set,
let
$
\eta\colonn O \to\R
$
be an upper semicontinuous function,
let
$
\phi
\colonn O \to[0,\infty)
$
be a lower semicontinuous function
satisfying
$
\lim_{ \alpha\to\infty}
\sup_{ y \in O }
(
\eta(y) - \alpha\cdot\phi(y)
)
\in\R
$
and
let
$
x \colonn(0,\infty) \to O
$
be a function satisfying
%
%e4.19 #&#
\begin{equation}
\label{eqassumptionProp37} \lim_{\alpha\to\infty} \Bigl( \sup_{
y \in O}
\bigl( \eta(y) - \alpha\cdot\phi(y) \bigr) - \bigl( \eta\bigl(
x(\alpha)
\bigr) - \alpha\cdot\phi\bigl( x(\alpha) \bigr) \bigr) \Bigr)
= 0.
\end{equation}
Then
$
\lim_{ \alpha\to\infty}
\alpha\cdot
\phi( x( \alpha) ) = 0
$
and
for all
$ \alpha_n \in( 0, \infty) $,
$ n \in\N$,
with\break 
$
\lim_{ n \to\infty}
\alpha_n = \infty
$
and
$
\lim_{ n \to\infty}
x( \alpha_n )
=:
x_0
\in O
$
it holds that
% $
% x_0 \in O \cap
% \big[
% \cap_{ \alpha> 0 }
% \overline{
% x( ( \alpha, \infty) )
% }
%
% \big]
% $
% every limit point $x_0\in O$ of the family
% $x(\alpha)\in O$, $\alpha\in(0,\infty)$, as $\alpha\to\infty$
% satisfies
% $\psi(x_0)=0$ and
% it holds that
% \begin{equation} \label{eq320}
% \forall
% x_0 \in O \cap
% \left(
% \bigcap_{ \alpha> 0 }
% \overline{
% x^{ - 1 }( ( \alpha, \infty) )
% }
% \right)
% \colonn
%
$
\phi( x_0 ) = 0
$
%
% \mbox{and}
%
and\break
$
\eta(x_0)
=
\lim_{
\alpha\to\infty
}
\sup_{ y \in O}
( \eta(y) - \alpha\cdot\phi( y ) )
= 
\sup_{ y \in\phi^{ - 1 }( 0 )
}
\eta( y )
$.
% \end{equation}
\end{lemma}

The next lemma essentially generalizes Theorem 2.2 in
Appendix C in Peng~\cite{Peng2010} (which assumes the
functions $ G_1, \ldots, G_k $ to be uniformly continuous
in the second argument uniformly in the last argument)
and is a generalized analog of
Theorem 8.2 in
Crandall, Ishii and~Lions~\cite{CrandallIshiiLions1992}
for unbounded domains.
Given an open set $ O \subset\R^d $, we define a sequence $ O_n
\subset O $, $ n \in\N$,
of compact sets by
%
%e4.20 #&#
\begin{equation}
O_n:= \biggl\{ x \in O \colonn\operatorname{dist}\bigl(x,
\R^d \setminus O\bigr) \ge\frac{ 1 }{ n } \mbox{ and } \llVert x
\rrVert\le n \biggr\}
\end{equation}
for all $ n \in\N$. We also write $ O_n^c:= O \setminus O_n $ for
the complement of $ O_n $ in $ O $.
% We then have

% The proof of
% Lemma~\ref{lcomparisonviscositysolution}
% below is an example application
% of the technical
% Lemma~\ref{ltechnicallemmauniqueness}.

%le4.10 #&#
\begin{lemma}[(A domination result
for viscosity subsolutions)]
\label{ltechnicallemmauniqueness}
Let $ T \in(0,\infty) $,
$ d, k \in\N$,
let
$ O \subset\R^d $
be an open set,
let
$
G_1, \ldots, G_k
\colonn
(0,T)
\times
O
\times
\R
\times
\R^d
\times\mathbb{S}_d
\to\R
$
be degenerate elliptic and upper semicontinuous functions
and let
$
u_1, \ldots, u_k \colonn[0,T]\times O \to\R
$
be upper semicontinuous functions such
that for every
$ i \in\{ 1, \ldots, k \} $
it holds that
$ u_i|_{ (0,T) \times O } $
is
a viscosity subsolution of
%
%e4.21 #&#
\begin{equation}
\label{eqparabolicequationi} \frac{ \partial}{ \partial t }
u_i(t,x) - G_i
\bigl( t, x, u_i(t,x), (\nabla_x u_i)
(t,x), (\operatorname{Hess}_xu_i) (t,x) \bigr) = 0
\end{equation}
for
$ (t,x) \in(0,T) \times O $.
Moreover, assume that
%
%e4.22 #&#
\begin{eqnarray}
\label{eq00assumption} &&\limsup_{
n \to\infty
} \Biggl[ \sum
_{ i = 1 }^k G_i \bigl(
t_i^{(n)}, x_i^{(n)},
r_i^{ (n) },\nonumber \\
&& \hspace*{ 70pt}n \bigl( \1_{ [ 2, k] }( i ) \cdot\bigl[
x_i^{ (n) } - x_{ i - 1 }^{ (n) } \bigr] +
\1_{
[1, k-1]
}( i ) \cdot\bigl[ x_{ i }^{ (n) } -
x_{ i + 1 }^{ (n) } \bigr] \bigr), %\nabla_{x_i}\big(\langle
%x^{(n)},B_{k,d}x^{(n)}
n
A_i^{ (n) } \bigr) \Biggr] \\
&&\qquad\leq0\nonumber
\end{eqnarray}
for all
$
(
t_i^{ (n) }, x_i^{ (n) }, r_i^{ (n) }, A_i^{ (n) }
)
\in
(0,T) \times O \times\R\times\mathbb{S}_d
$,
$ n \in\N$,
$ i \in\{ 1, \ldots, k \} $,
satisfying that
$
\lim_{ n \to\infty}
( t_1^{ (n) }, x_1^{ (n) } )
\in( 0, T ) \times O
$,
that
$
\lim_{ n \to\infty}
(
\sqrt{ n }
\sum_{ i = 2 }^k
\|
( t_i^{ (n) }, x_i^{ (n) } )
-
( t_{ i - 1 }^{ (n) }, x_{ i - 1 }^{ (n) } )
\|
)
= 0
$,
that
$
\lim_{ n \to\infty}
\sum_{ i = 1 }^k
r_i^{ (n) }
> 0
$,
that
$
\sup_{ n \in\N}
$
$
\sum_{ i = 1 }^k
| r_i^{ (n) } |
< \infty
$
and that
$
\forall
n \in\N\colonn
\forall
z_1, \ldots, z_k \in\R^d
\colonn
- 5
\sum_{ i = 1 }^k
\| z_i \|^2
\leq
\sum_{ i = 1 }^k
\langle z_i, A_i^{(n)} z_i
\rangle
\leq\break 
5
\sum_{ i = 2 }^k
\| z_i - z_{ i - 1 }
\|^2
$.
Furthermore, assume that
$
\sum_{ i = 1 }^k
u_i(0,x)
\leq0
$
for all
$ x \in O $
and that
%
%e4.23 #&#
\begin{equation}
\label{eqattainsmaximum} \lim_{ n \to\infty} \sup_{
(t,x) \in
(0,T) \times O_n^c
} \sum
_{ i = 1 }^k u_i(t,x) \leq0.
\end{equation}
Then
$
\sum_{ i = 1 }^k u_i(t, x) \leq0
$
for all
$
(t,x) \in[0,T) \times O
$.
\end{lemma}

\begin{pf}%{Proof
%of
%Lemma~\ref{ltechnicallemmauniqueness}}
%
If $ O = \varnothing$,
then the assertion is trivial.
So for the rest of
the proof, we assume that
$ O \neq\varnothing$.
We will show
that
$
\sum_{ i = 1 }^k
u_i(t,x)
\leq
\frac{ k \delta}{ (T - t) }
$
for all
$
(t,x) \in[0,T) \times O
$
and all
$
\delta\in(0,1]
$.
Letting $\delta\to0 $
will then yield that
$
\sum_{ i = 1 }^k
u_i( t, x )
\leq0
$
for all
$
(t,x) \in[0,T) \times O
$.
In the following, we thus fix
$ \delta\in(0,1] $.
In a first step of this proof,
we modify the problem.
% through an application of
% Lemma~\ref{lemaffinelinear}.
More precisely,
define functions
$
\tilde{u}_1,
\ldots,
\tilde{u}_k
\colonn
[0,T) \times O
\to[ - \infty, \infty)
$
by
$
\tilde{u}_i(t,x)
:=
u_i( t, x ) -
\frac{ \delta}{ (T - t) }
$
and
functions
$
\tilde{G}_1, \ldots,
\tilde{G}_k \colonn
(0,T)
\times
O
\times\R\times\R^d
\times\mathbb{S}_d
\to\R
$
by
%
%e4.24 #&#
\begin{eqnarray}
\label{eqdeftildeG} \tilde{G}_i ( t, x, r, p, A )&:=& G_i
\biggl( t, x, r + \frac{ \delta}{ (T - t) }, p, A \biggr) - \frac
{ \partial}{ \partial t } \biggl(
\frac{ \delta}{ (T - t) } \biggr)
\nonumber
\\[-8pt]
\\[-8pt]
\nonumber
& = &G_i \biggl( t, x, r +
\frac{ \delta}{ (T - t) }, p, A \biggr) - \frac{ \delta}{ (T -
t)^2 }.
\end{eqnarray}
Then it holds
for every $ i \in\{ 1, \ldots, k \} $
that
$
\tilde{u}_i|_{ (0,T) \times O }
$
is a viscosity subsolution of
%
%e4.25 #&#
\begin{equation}
\label{eqparabolicequationtildeui} \frac{ \partial}{ \partial t
} \tilde{u}_i(t,x) -
\tilde{G}_i \bigl( t, x, \tilde{u}_i(t,x), (
\nabla_x \tilde{u}_i) (t,x), (\operatorname{Hess}_x
\tilde{u}_i) (t,x) \bigr) = 0
\end{equation}
for
$ (t,x) \in(0,T) \times O $.
It remains to prove that
$
\sum_{ i = 1 }^k
\tilde{u}_i( z ) \leq0
$
for all
$
z \in[0,T) \times O
$.
Aiming at a contradiction, we assume that the
extended real number
$
S_0:=
\sup_{ z \in[0,T) \times O }
\sum_{ i = 1 }^k
\tilde{u}_i( z )
\in
( - \infty, \infty]
$
satisfies that
$
S_0
\in
(0,\infty]
$.
Assumption~\eqref{eqattainsmaximum}
then implies that
there exists a natural
number
$ n_0 \in\N$
such that $K:= O_{n_0}$ is nonempty
and such that
$
\sum_{ i = 1 }^k
\tilde{u}_i( z )
\leq
\sum_{ i = 1 }^k
u_i( z )
\leq
\min( 1, \frac{ S_0 }{ 2 } )
$
for all
$
z \in
(0,T) \times K^c
$.
This, together
with
$
\sum_{ i = 1 }^k
\tilde{u}_i( 0, x)
\leq\break 
\sum_{ i = 1 }^k
u_i( 0, x)
\leq
0
$
and
$
\sum_{ i = 1 }^k
\tilde{u}_i( T, x )
= - \infty
$
for all
$
x \in O
$
implies that
%
%e4.26 #&#
\begin{equation}
\label{eqsumtildeuoutside} \sup_{
z \in
[0,T] \times K^c
} \sum_{ i = 1 }^k
\tilde{u}_i( z ) \leq\min\biggl( 1, \frac{ S_0 }{ 2 } \biggr)
\leq\frac{ S_0 }{ 2 }.
\end{equation}
Moreover, the
function
$
\sum_{ i = 1 }^k
\tilde{u}_i
\colonn
[0,T] \times O
\to[ - \infty, \infty)
$
is
upper
semicontinuous
and
is hence bounded from above
on the compact set
$ [0,T] \times K$.
Combining this with \eqref{eqsumtildeuoutside}
proves that
$
S_0 < \infty
$
and we thus get
$
S_0 \in( 0, \infty)
$.
In the next step, we
define a function
$
\phi\colonn
(
[0,T] \times O
)^k \to[0,\infty)
$
by
% \begin{equation}
% \label{eqdefpsi}
$
\phi(
z_1, \ldots, z_k
)
=
\frac{ 1 }{ 2 }
\sum_{ i = 2 }^k
\llVert
z_i - z_{ i - 1 }
\rrVert^2
$
% \end{equation}
for all
$
z_1, \ldots, z_k \in[0,T] \times O
$.
For several
$ n \in\N$,
we will apply
Lemma~\ref{lTheorem32UserGuidy}
with
$
\CO= (0,T) \times O
$,
$
\varepsilon= \frac{ 1 }{ n }
$
and with
$
\Phi=
n \cdot
\phi|_{
(
(0,T) \times O
)^k
}
$
below.
For this, we
now check the
assumptions
of Lemma~\ref{lTheorem32UserGuidy}.
Define a function
$
\eta\colonn
( [0,T] \times K )^k
\to[-\infty,\infty)
$
by
$
\eta( z_1, \ldots, z_k )
=
\sum_{ i = 1 }^k
\tilde{u}_i( z_i )
$
for all
$
z_1, \ldots, z_k \in[0,T] \times K
$.
Note for every
$ \alpha\in(0,\infty) $
that the function
$
( [0,T] \times K )^k
\ni
\underline{z}
\mapsto
\eta(
\underline{z}
)
-
\alpha\cdot
\phi( \underline{z} )
\in
[ - \infty, \infty)
$
is upper semicontinuous
with a compact domain of definition
and therefore,
attains its maximum
$
S_{ \alpha}
:=
\sup_{
\underline{z}
\in
( [0,T] \times K )^k
}
(
\eta(
\underline{z}
)
-
\alpha\cdot
\phi(
\underline{z}
)
)
< \infty
$
in a point
$
\underline{z}^{ ( \alpha) }
=
(
(
t^{ ( \alpha) }_1,
x^{ ( \alpha) }_1
),
\ldots,
(
t^{ ( \alpha) }_k,
x^{ ( \alpha) }_k
)
)
\in( [0,T] \times K )^k
$.
Next observe that
% the fact that $ S_0 > 0 $
% ensures that
%
%e4.27 #&#
\begin{equation}
\label{eqPropertySalpha} \infty> S_\alpha\geq\sup_{ z \in[0,T)
\times K }
\eta( z, z, \ldots, z ) = \sup_{ z \in[0,T) \times K } \sum
_{ i = 1 }^k \tilde{u}_i( z ) =
S_0 > 0
\end{equation}
for all
$
\alpha\in(0,\infty)
$.
This together with monotonicity
of the function
$
( 0, \infty) \ni\alpha
\mapsto S_{ \alpha}
\in( 0, \infty)
$
implies that
the limit
$
\lim_{ \alpha\to\infty}
S_{ \alpha}
$
exists in
$
(0, \infty)
$, that is, it holds that
$
\lim_{ \alpha\to\infty}
S_{ \alpha}
\in( 0, \infty)
$.
The set
$
\{
\underline{z}^{ (n) }
\colonn
n \in\N
\}
\subset
( [0,T] \times K )^k
$
is
relatively compact and, therefore,
there exists a limit point
$
\underline{ \hat{z} }
=
% \big(
% \hat{z}_1, \ldots, \hat{z}_k
% \big)
% =
(
( \hat{t}_1, \hat{x}_1 ),
\ldots,\break 
( \hat{t}_k, \hat{x}_k )
)
\in
( [0,T] \times K )^k
$
of this set.
Let
$
n_j \in\N
$,
$ j \in\N$,
be a strictly increasing sequence such that
$
\lim_{ j \to\infty}
\underline{z}^{ ( n_j ) }
=
\underline{ \hat{z} }
$.
Clearly,
$
\tilde{u}_i(T,x) = - \infty
$
for all
$ x \in K $ and all
$ i \in\{ 1, \ldots, k \} $
implies
that
$
t^{ ( \alpha) }_1, \ldots, t^{ ( \alpha) }_k
\in[0, T)
$
for all
$
\alpha\in(0,\infty)
$.
In addition, observe that if
$
( \hat{t}_1, \ldots, \hat{t}_k )
\in
[0,T]^k \setminus[0,T)^k
$,
then
\eqref{eqPropertySalpha}
implies that
%
%e4.28 #&#
\begin{eqnarray}
 0 &< &% S_0
% \leq
% \lim_{ \alpha\to\infty}
% S_{ \alpha}
% =
\lim
_{ j \to\infty} S_{ n_j } = \lim_{ j \to\infty} \bigl(
\eta\bigl( \underline{z}^{ ( n_j ) } \bigr) - n_j \cdot\phi
\bigl( \underline{z}^{ ( n_j ) } \bigr) \bigr) % \\ &
\leq\lim
_{ j \to\infty} \eta\bigl( \underline{z}^{ ( n_j ) } \bigr)
\nonumber
\\[-8pt]
\\[-8pt]
\nonumber
&\leq&
\Biggl( \sum_{ i = 1 }^k \Bigl[ \sup
_{
z \in[0,T] \times K
} u_i( z ) \Bigr] \Biggr) - \infty= -
\infty
\end{eqnarray}
and this contradiction shows that
$
( \hat{t}_1, \ldots, \hat{t}_k )
\in[ 0, T )^k
$.
Next observe that
%
%e4.29 #&#
\begin{eqnarray}
&&\lim_{ \alpha\to\infty} \Bigl[ \sup_{ z \in( [0,T) \times K )^k
} \bigl( \eta(
z ) - \alpha\cdot\phi( z ) \bigr) - \bigl( \eta\bigl( \underline
{z}^{ (\alpha) }
\bigr) - \alpha\cdot\phi\bigl( \underline{z}^{ (\alpha) } \bigr
) \bigr)
\Bigr]
\nonumber
\\[-8pt]
\\[-8pt]
\nonumber
&&\qquad= \lim_{ \alpha\to\infty} [ S_{ \alpha} - S_{ \alpha} ] = 0
.
\end{eqnarray}
Hence, Lemma~\ref{lProp37UserGuide}
applied to
$
\eta|_{
( [0,T) \times K )^k
}
$
and to
$
\phi|_{
( [0,T) \times K )^k
}
$
yields that
%
%e4.30 #&#
\begin{equation}
\label{eqetaconsequencelimes} 0 = \lim_{
\alpha
\to\infty
} \bigl[ \alpha\cdot\phi\bigl(
\underline{z}^{ (\alpha) } \bigr) \bigr] = \lim_{
\alpha\to\infty
} \Biggl[
\frac{ \alpha}{ 2 } \sum_{ i = 2 }^k \bigl\|
\bigl( t^{ (\alpha) }_i, x^{ (\alpha) }_i \bigr) -
\bigl( t^{ (\alpha) }_{ i - 1 }, x^{ (\alpha) }_{ i - 1 } \bigr)
% z^{ ( \alpha) }_i
% -
% z^{ ( \alpha) }_{ i - 1 }
\bigr\|^2 \Biggr]
\end{equation}
and
that
$
\phi(
\underline{ \hat{z} }
) = 0
$.
The definition of $ \phi$
% (see \eqref{eqdefpsi})
therefore ensures that
$
( \hat{t}_i, \hat{x}_i )
=
( \hat{t}_j, \hat{x}_j )
$
for all
$
i, j \in
\{ 1, \ldots, k \}
$.
% and that $\lim_{(0,\infty)\ni\alpha\to\infty}S_\alpha=S_0>0$.
% This shows that all assumptions of Lemma~\ref{lTheorem32UserGuidy}
% are satisfied.
Furthermore,
observe that if
$
\hat{t}_1 = 0
$,
then
\eqref{eqPropertySalpha} and
the upper semicontinuity
of
$
\eta
$
show that
%
%e4.31 #&#
\begin{eqnarray}
0 &<& S_0 \leq\lim_{ j \to\infty} S_{ n_j } \leq
\limsup_{ j \to\infty} \eta\bigl( \underline{z}^{ ( n_j ) } \bigr)
\leq\eta( \underline{ \hat{z} } ) = \sum_{ i = 1 }^k
\tilde{u}_i ( \hat{t}_1, \hat{x}_1 )
\nonumber
\\[-8pt]
\\[-8pt]
\nonumber
&=& \sum
_{i=1}^k u_i ( 0,
\hat{x}_1 ) - \frac{k\delta}{T} \leq0
\end{eqnarray}
and this contradiction implies that
$
\hat{t}_1 = \hat{t}_2 = \cdots= \hat{t}_k
\in(0,T)
$.
Consequently, there exists a natural number
$
j_0 \in\N
$
such that for every
$
j \in\{ j_0, j_0 + 1, \ldots\}
$
it holds that
$
t^{ (n_j) }_1, \ldots, t^{ (n_j) }_k
\in(0,T)
$.
Next,
for every
$
n \in
\mathcal{N}
:=
\{
m \in\N\colonn
t^{ (m) }_1, \ldots, \break t^{ (m) }_k \in(0,T)
\}
$,
we apply Lemma~\ref{lTheorem32UserGuidy}
with
$
\CO= (0,T) \times O
$,
with
$
\eps= \frac{ 1 }{ n }
$,
with
the functions
$
\tilde{u}_1|_{ (0,T) \times O },
\ldots,
\tilde{u}_k|_{ (0,T) \times O }
$
and
$
\Phi
=
n \cdot\phi|_{ ( (0,T) \times O )^k
}
$
and
with the local maximum point
$
\underline{z}^{ (n) }
\in( (0,T) \times O )^k
$
to obtain the existence
of matrices
$
( A^{ (n) }_1, \ldots, A^{ (n) }_k )
=
(
(
a^{ n, 1 }_{ i, j }
)_{ i, j \in\{ 1, \ldots, d + 1 \} },
\ldots,
(
a^{ n, k }_{ i, j }
)_{ i, j \in\{ 1, \ldots, d + 1 \} }
)
\in
(
\mathbb{S}_{ d + 1 }
)^k
$,
$ n \in\mathcal{N} $,
such that
for every $ n \in\mathcal{N} $
and every $ i \in\{ 1, \ldots, k \} $
it holds that
%
%e4.32 #&#
\begin{equation}
\label{eqliesinJplus}\quad \bigl( n ( \nabla_{ (t_i, x_i) } \phi)
\bigl( \bigl(
t^{ (n) }_1, x^{ (n) }_1 \bigr), \ldots,
\bigl( t^{ (n) }_k, x^{ (n) }_k \bigr)
\bigr), n A_i^{ (n) } \bigr) \in\bigl( \hat{J}^2_{ + }
\tilde{u}_i \bigr) \bigl( t^{ (n) }_i,
x^{ (n) }_i \bigr)
\end{equation}
and
\begin{eqnarray*}
- \bigl[ n + n\bigl \| ( \operatorname{Hess} \phi) \bigl( \underline
{z}^{ (n) }
\bigr)\bigr \|_{
L (
\R^{ ( d + 1 ) k }
)
} \bigr] I &\leq&\pmatrix{ n
A_1^{ (n) } & \cdots& 0
\vspace*{2pt}\cr
\vdots& \ddots& \vdots
\vspace*{2pt}\cr
0 & \cdots& n A_k^{ (n) }
}
% \\ &
\\
&\leq& n ( \operatorname{Hess} \phi) \bigl( \underline{z}^{ (n) }
\bigr) + \frac{ 1 }{ n } \bigl[ n ( \operatorname{Hess} \phi)
\bigl(
\underline{z}^{ (n) } \bigr) \bigr]^2.
\end{eqnarray*}
Combining this with
the identity
$
(\operatorname{Hess} \phi)( z ) =
(\operatorname{Hess} \phi)( 0 )
$
for all $ z \in( (0,T) \times O )^k $
then implies that
%
%e4.33 #&#
\begin{eqnarray}
\label{eqmatrixAest}  - \bigl[ 1 + \bigl\| ( \operatorname{Hess} \phi
) ( 0 )
\bigr\|_{
L (
\R^{ ( d + 1 ) k }
)
} \bigr] I &\leq&\pmatrix{
A_1^{ (n) } & \cdots& 0
\vspace*{2pt}\cr
 \vdots& \ddots& \vdots
\vspace*{2pt}\cr
0 & \cdots& A_k^{ (n) }
}
\nonumber
\\[-8pt]
\\[-8pt]
\nonumber
&\leq&( \operatorname{Hess} \phi) (0) + \bigl[ ( \operatorname{Hess}
\phi) ( 0 )
\bigr]^2
\end{eqnarray}
for all $ n \in\mathcal{N} $.
To simplify the notation we define matrices
$
B^{ (n) }_l \in\mathbb{S}_d
$,
$ l \in\{ 1, \ldots, k \} $,
$ n \in\mathcal{N} $,
by\vspace*{1pt}
$
B^{ (n) }_l
:=
(
a^{ n, l }_{ i + 1, j + 1 }
)_{
i, j \in\{ 1, \ldots, d \}
}
$
for all
$ l \in\{ 1, \ldots, k \} $
and all
$ n \in\mathcal{N} $.
Corollary~\ref{corsemijetsequivalence2} together
with \eqref{eqliesinJplus}
and the fact that it holds for every
$ i \in\{ 1, \ldots, k \} $
that
$
\tilde{u}_i|_{
(0, T) \times O
}
$
is a viscosity subsolution \eqref{eqparabolicequationtildeui}
then proves that
%
%e4.34 #&#
\begin{eqnarray}
&&n \biggl( \frac{ \partial}{ \partial t_i } % \nabla_{ t_i }
\phi\biggr) \bigl(
\underline{z}^{ (n) } \bigr) - \tilde{G}_i \bigl(
t_i^{ (n) }, x_i^{ (n) },
\tilde{u}_i \bigl( t_i^{ (n) },
x_i^{ (n) } \bigr), n ( \nabla_{ x_i } \phi)
\bigl( \underline{z}^{ (n) } \bigr), n B_i^{(n)}
\bigr) 
\nonumber
\\[-8pt]
\\[-8pt]
\nonumber
&&\qquad\leq0
\end{eqnarray}
for all
$ i \in\{ 1, \ldots, k \} $
and all
$ n \in\mathcal{N} $.
Summing over $ i \in\{ 1, \ldots, k \} $
hence results in
%
%e4.35 #&#
\begin{eqnarray}
\label{eqsumG}&& n \sum_{ i = 1 }^k \biggl(
\frac{ \partial}{ \partial t_i } \phi\biggr) \bigl( \underline
{z}^{ (n) } \bigr)
\nonumber
\\[-8pt]
\\[-8pt]
\nonumber
&&\qquad\leq
\sum_{ i = 1 }^k \tilde{G}_i
\bigl( t_i^{ (n) }, x_i^{ (n) },
\tilde{u}_i \bigl( t_i^{ (n) },
x_i^{ (n) } \bigr), n ( \nabla_{ x_i } \phi)
\bigl( \underline{z}^{ (n) } \bigr), n B_i^{(n)}
\bigr)
\end{eqnarray}
for all
$ n \in\mathcal{N} $.
Next note that the definition of
$ \phi$
% (see \eqref{eqdefpsi})
ensures in the case $ k \geq2 $ that
%
%e4.36 #&#
\begin{eqnarray}
\label{eqnablat} \biggl( \frac{ \partial}{ \partial t_i } \phi
\biggr) \bigl( ( t_1,
x_1 ), \ldots, ( t_k, x_k ) \bigr)& =&
\frac{ 1 }{ 2 } % \sum_{ i = 1 }^k
\sum_{ j = 2 }^k
\frac{ \partial}{ \partial t_i } ( t_j - t_{ j - 1 } )^2
\nonumber
\\[-8pt]
\\[-8pt]
\nonumber
&=&
\cases{ t_1 - t_2, & \quad$ i = 1,$ \vspace*{2pt}
\cr
2
t_i - t_{ i - 1 } - t_{ i + 1 }, & \quad $1 < i < k,$
\vspace*{2pt}
\cr
t_k - t_{ k - 1 }, &\quad  $ i = k,$ }
\end{eqnarray}
for all
$
i \in\{ 1, \ldots, k \}
$
and all
$
(t_1,x_1), \ldots,
(t_k,x_k) \in(0,T) \times O
$
and, therefore,
we obtain that in the case $ k \geq2 $
it holds that
%
%e4.37 #&#
%e4.38 #&#
%e4.39 #&#
\begin{eqnarray}
&& \sum_{ i = 1 }^k \biggl(
\frac{ \partial}{ \partial t_i } \phi\biggr) \bigl( ( t_1, x_1 ),
\ldots, ( t_k, x_k ) \bigr) \nonumber\\
&&\qquad= t_1 -
t_2 + t_k - t_{ k - 1 } + \sum
_{ i = 2 }^{ k - 1 } ( 2 t_i - t_{ i - 1 }
- t_{ i + 1 } )
\nonumber
\\
&&\qquad= t_1 - t_2 + t_k -
t_{ k - 1 }
 + \Biggl( \sum_{ i = 2 }^{ k - 1 }
t_i - t_{ i - 1 } \Biggr) + \Biggl( \sum
_{ i = 2 }^{ k - 1 } t_i - t_{ i + 1 }
\Biggr) \\
&&\qquad= \Biggl( t_1 - t_{ k - 1 } + \sum
_{ i = 2 }^{ k - 1 } t_i - t_{ i - 1 }
\Biggr) + \Biggl( t_k - t_2 + \sum
_{ i = 2 }^{ k - 1 } t_i - t_{ i + 1 }
\Biggr) \nonumber\\
&&\qquad= 0\nonumber
\end{eqnarray}
for all
$
(t_1,x_1), \ldots,
(t_k,x_k) \in(0,T) \times O
$.
Combining this with
\eqref{eqsumG}
results in
%
%e4.40 #&#
\begin{equation}
 0 \leq\sum_{ i = 1 }^k
\tilde{G}_i \bigl( t_i^{ (n) },
x_i^{ (n) }, \tilde{u}_i \bigl(
t_i^{ (n) }, x_i^{ (n) } \bigr), n (
\nabla_{ x_i } \phi) \bigl( \underline{z}^{ (n) } \bigr), n
B_i^{ (n)} \bigr)
\end{equation}
for all
$ n \in\mathcal{N} $.
Therefore,
we obtain from
\eqref{eqdeftildeG}
and from
$
\hat{t}_1 = \cdots= \hat{t}_k
\in( 0, T )
$
and
$
t^{ (n_j) }_1,
\ldots,
t^{ (n_j) }_k
\in(0,T)
$
for all $ j \in\{ j_0, j_0 + 1, \ldots\} $
that
%
%e4.41 #&#
\begin{eqnarray}
 \label{eqGdeltaestimate} &&\sum_{ i = 1 }^k
\frac{
\delta
}{
(
T - t^{ (n_j) }_i
)^{ 2 }
} \nonumber\\
&&\qquad\leq\sum_{ i = 1 }^k
G_i \biggl( t^{ (n_j) }_i, x^{ (n_j) }_i,
\tilde{u}_i \bigl( t_i^{ (n_j) },
x_i^{ (n_j) } \bigr) \\
&&\hspace*{35pt}\qquad\quad{}+ \frac{ \delta
}{
( T - t^{ (n_j) }_i
)
}, n_j
( \nabla_{ x_i } \phi) \bigl( \underline{z}^{ (n_j) } \bigr),
n_j B_i^{ (n_j)} \biggr)\nonumber\
\end{eqnarray}
for all
$ j \in\{ j_0, j_0, \ldots\} $.
In the next step, we define
$
(
\mathbf{t}^{ (n) }_i,
\mathbf{x}^{ (n) }_i,
\mathbf{r}^{ (n) }_i,
\mathbf{A}^{ (n) }_i
)
\in
( 0, T ) \times
O \times
\R\times
\mathbb{S}_d
$,
$
i \in\{ 1, \ldots, k \}
$,
$ n \in\N$,
by
%
%e4.42 #&#
\begin{eqnarray}
&&\bigl( \mathbf{t}^{ (n) }_i, \mathbf{x}^{ (n) }_i,
\mathbf{r}^{ (n) }_i, \mathbf{A}^{ (n) }_i
\bigr)
\nonumber
\\[-8pt]
\\[-8pt]
\nonumber
&&\qquad:= \cases{\displaystyle\biggl( t^{ (n) }_i, x^{ (n) }_i,
\tilde{u}_i \bigl( t_i^{ (n) },
x_i^{ (n) } \bigr) + \frac{ \delta}{ ( T - t_i^{ (n) } ) },
B_i^{ (n)} \biggr), \vspace*{2pt}\cr
\qquad  n \in\bigl\{
n_j \in\N\colonn j \in\{ j_0, j_0 + 1,
\ldots\} \bigr\}, \vspace*{2pt}
\cr
\displaystyle\biggl( \hat{t}_1,
\hat{x}_1, \frac{
\lim_{ \alpha\to\infty}
S_{ \alpha}
}{ k } + \frac{ \delta}{
( T - \hat{t}_1 )
}, 0 \biggr), \vspace*{2pt}\cr
\qquad\mbox{else}, }
\end{eqnarray}
for all
$ i \in\{ 1, \ldots, k \} $
and all
$ n \in\N$.
Moreover, observe that in the case $ k \geq2 $ it holds that
%
%e4.43 #&#
\begin{eqnarray}
\label{eqnablax} ( \nabla_{ x_i } \phi) \bigl( ( t_1,
x_1 ), \ldots, ( t_k, x_k ) \bigr) &=&
\frac{ 1 }{ 2 } \sum_{ j = 2 }^k
\nabla_{ x_i } \bigl( \llVert x_j - x_{ j - 1 }
\rrVert^2 \bigr)
\nonumber
\\[-8pt]
\\[-8pt]
\nonumber
&=& \cases{ x_1 - x_2, & \quad$i = 1,$ \vspace*{2pt}
\cr
2 x_i - x_{ i - 1 } -
x_{ i + 1 }, &\quad $1 < i < k,$ \vspace*{2pt}
\cr
x_k -
x_{ k - 1 }, & \quad $i = k$}
\end{eqnarray}
for all
$
i \in\{ 1, \ldots, k \}
$
and all
$
( t_1, x_1 ), \ldots, ( t_k, x_k ) \in(0, T) \times O
$.
% Combining this with \eqref{eqGdeltaestimate}
% results in
% \begin{eqnarray}
% \label{eqGliminfbigger0}
% \frac{ k \delta
% }{
% \left( T - \hat{t}_1 \right)^2
% }
% \leq
% \liminf_{ j \to\infty}
% \Bigg[
% \sum_{ i = 1 }^k
% G_i\Big(
% t^{ (n_j) }_i, x^{ (n_j) }_i,
% \tilde{u}_i\big(
% t_i^{ (n_j) }, x_i^{ (n_j) }
% \big)
% +
% \frac{ \delta
% }{
% \left( T - t^{ (n_j) }_i \right)
% }
%,
% n_j
%
% \big(
% \mathbh{1}_{
% [2, \infty)
% }(i)
% \cdot
% [
% x_i^{ (n_j) }
% -
% x_{ i - 1 }^{ (n_j) }
% ]
% \\
% +
% \mathbh{1}_{
% [0, k-1]
% }(i)
% \cdot
% [
% x_i^{ (n_j) }
% -
% x_{ i + 1 }^{ (n_j) }
% ]
% \big)
%,
% n_j
% B_i^{ (n_j)}
% \Big)
% \Bigg].
% \end{eqnarray}
% Then \eqref{eqGliminfbigger0}
Then~\eqref{eqGdeltaestimate}
ensures that
%
%e4.44 #&#
\begin{eqnarray}
\label{eqGliminfbigger1}&& \frac{ k \delta
}{
( T - \hat{t}_1 )^2
} \leq\limsup_{ n \to\infty} \Biggl[
\sum_{ i = 1 }^k G_i \bigl( {
\mathbf{ t}}^{ (n) }_i, \mathbf{x}^{ (n) }_i, {
\mathbf{ r}}^{ (n) }_i,\nonumber\\
 && \hspace*{128pt}n \bigl( \mathbh{1}_{
[2, k]
}(i)
\cdot\bigl[ \mathbf{x}_i^{ (n) } - \mathbf{x}_{ i - 1 }^{ (n) }
\bigr]\\
&&\hspace*{138pt}{} + \mathbh{1}_{
[1, k-1]
}(i) \cdot\bigl[ \mathbf{x}_i^{ (n) }
- \mathbf{x}_{ i + 1 }^{ (n) } \bigr] \bigr), n \mathbf{
A}_i^{ (n)} \bigr) \Biggr].\nonumber
\end{eqnarray}
Next, we observe that the Taylor expansion
$
\phi(z)
=
\phi( 0 )
+
\anglel
( \nabla\phi)( 0 ), z
\angler
+
\frac{ 1 }{ 2 }
\langle
z, ( \operatorname{Hess} \phi)(0) z
\rangle
=
\frac{ 1 }{ 2 }
\langle
z, ( \operatorname{Hess} \phi)(0) z
\rangle
$
for all
$
z \in\R^{ (d + 1) k }
$
implies that\break 
$
(\nabla\phi)(z)
=
( \operatorname{Hess} \phi)(0)
z
$
for all
$
z \in\R^{ (d + 1) k }
$.
This
together with
\eqref{eqnablat},
\eqref{eqnablax}
and the estimate
$
(a + b)^2
\leq
2 a^2 + 2 b^2
$
for all
$ a, b \in\R$
results in
%
%e4.45 #&#
%e4.46 #&#
%e4.47 #&#
%e4.48 #&#
\begin{eqnarray}
\label{eqB2leqB} \bigl\langle z, \bigl( ( \operatorname{Hess}
\phi) (0)
\bigr)^2 z \bigr\rangle&=& \bigl\langle( \operatorname{Hess} \phi)
(0) z, ( \operatorname{Hess} \phi) (0) z \bigr\rangle= \bigl\| (
\operatorname{Hess}
\phi) (0) z \bigr\|^2\nonumber\\
& =&\bigl \| ( \nabla\phi) (z) \bigr\|^2\nonumber
\\
& =& \| z_1 - z_2 \|^2 + \Biggl[ \sum
_{ i = 2 }^{ k - 1 } \| 2 z_i -
z_{ i - 1 } - z_{ i + 1 } \|^2 \Biggr] + \|
z_k - z_{ k - 1 } \|^2
\\
& \leq&2 \| z_1 - z_2 \|^2 + \Biggl[ \sum
_{ i = 2 }^{ k - 1 } 2 \bigl( \| z_i -
z_{ i - 1 } \|^2 + \| z_i - z_{ i + 1 }
\|^2 \bigr) \Biggr] \nonumber\\
&&{}+ 2 \llVert z_k - z_{ k - 1 }
\rrVert^2\nonumber
\\
& =& 4 \Biggl[ \sum_{ i = 2 }^k \|
z_i - z_{ i - 1 } \|^2 \Biggr] \leq8 \Biggl[
\sum_{ i = 2 }^k \| z_i
\|^2 \Biggr] + 8 \Biggl[ \sum_{ i = 2 }^k
\| z_{ i - 1 } \|^2 \Biggr] \nonumber\\
&\leq&16 \llVert z \rrVert
^2\nonumber
\end{eqnarray}
for all
$
z = (z_1, \ldots, z_k)
\in\R^{ ( d + 1 ) k }
$.
Inequality~\eqref{eqB2leqB}
implies that\break 
$
\|
( \operatorname{Hess} \phi)(0)
\|_{
L( \R^{ (d+1) \times k } )
}
\leq4
$.
Consequently,
\eqref{eqmatrixAest},
\eqref{eqB2leqB}
and
$
\langle
z,\break  ( \operatorname{Hess} \phi)(0) z
\rangle
=
2 \phi(z)
$
for all
$
z \in
\R^{ (d + 1) k }
$
yield that
%
%e4.49 #&#
\begin{eqnarray}
\label{eq55inequality1} - 5 \llVert z \rrVert^2& \leq&\sum
_{ i = 1 }^k \bigl\langle z_i,
A_i^{ (n) } z_i \bigr\rangle% \\ &
\leq
2 \phi(z) + \bigl\langle z, \bigl( ( \operatorname{Hess} \phi) (
0 )
\bigr)^2 z \bigr\rangle
\nonumber
\\[-8pt]
\\[-8pt]
\nonumber
&\leq& 5 \sum_{ i = 2 }^k
\| z_i - z_{ i - 1 } \|^2
\end{eqnarray}
for all
$
z = ( z_1, \ldots, z_k )
\in\R^{ ( d + 1 ) k }
$
and all
$
n \in\mathcal{N}
$.
Inequality~\eqref{eq55inequality1},
in particular, implies
$
- 5
\llVert z \rrVert^2
\leq
\sum_{ i = 1 }^k
\langle
z_i $,
$
B_i^{ (n) }
z_i
\rangle
=
\sum_{ i = 1 }^k
\langle
z_i,
\mathbf{A}_i^{ (n) }
z_i
\rangle
\leq
5
\sum_{ i = 2 }^k
\|
z_i - z_{ i - 1 }
\|^2
$
for all
$
z =
(z_1, \ldots, z_k)
\in\R^{ d k }
$
% and all
% $
% n \in\mathcal{N}
% $
% and this ensures that
% $
% - 5
% \left\| z \right\|^2
% \leq
% \sum_{ i = 1 }^k
% \langle
% z_i,
% \mathbf{A}_i^{ (n) }
% z_i
% \rangle
% \leq
% 5
% \sum_{ i = 2 }^k
% \|
% z_i - z_{ i - 1 }
% \|^2
% $
% for all
% $
% z =
% (z_1, \ldots, z_k)
% \in\R^{ d k }
% $
and all
$
n \in\N
$.
Combining this, the identities
%
%e4.50 #&#
\begin{eqnarray}
 \lim_{ j \to\infty} \Biggl[ \sum_{ i = 1 }^k
\biggl( \tilde{u}_i \bigl( t_i^{ (n_j) },
x_i^{ (n_j) } \bigr) + \frac{ \delta}{
( T - t_i^{ (n_j) } )
} \biggr) \Biggr]
% =
% \lim_{ j \to\infty}
% \left[
% \eta\big(
% (
% t_1^{ (n_j) },
% x_1^{ (n_j) }
% ),
% \ldots,
% (
% t_k^{ (n_j) },
% x_k^{ (n_j) }
% )
% \big)
% \right]
% +
% \frac{ k \delta}{
% \left(
% T - \hat{t}_1
% \right)
% }
% \\ &
&= &\Bigl( \lim
_{ j \to\infty} S_{ n_j } \Bigr) + \frac{
k \delta
}{
( T - \hat{t}_1 )
}
% =
% S_0
% +
% \frac{ k \delta}{
% \left( T - \hat{t}_1 \right)
% }
\nonumber
\\[-8pt]
\\[-8pt]
\nonumber
&=&
\lim_{ n \to\infty} \Biggl[ \sum_{ i = 1 }^k
\mathbf{r}^{ (n) }_i \Biggr] > 0,
\end{eqnarray}
%
%and
$
\lim_{
n \to\infty
}
n
\sum_{ i = 2 }^k
\|
(
\mathbf{t}^{ ( n ) }_i, \mathbf{x}^{ ( n ) }_i
)
-
(
\mathbf{t}^{ ( n ) }_{ i - 1 },
\mathbf{x}^{ ( n ) }_{ i - 1 }
)
\|^2
= 0
$
[see \eqref{eqetaconsequencelimes}]
and the estimate
$
\sup_{ n \in\N}
\max_{ i \in\{ 1, \ldots, k \} }
| \mathbf{r}^{ (n) }_i |
< \infty
$
with
assumption~\eqref{eq00assumption}
and with \eqref{eqGliminfbigger1}
shows that
\begin{eqnarray*}
0 &<& \frac{ k \delta}{ ( T - \hat{t}_1 )^2 } \\
&\leq&\limsup_{ n \to
\infty} \Biggl[ \sum
_{ i = 1 }^k G_i \bigl( \mathbf{
t}^{ (n) }_i, \mathbf{x}^{ (n) }_i, \mathbf{
r}^{ (n) }_i, n \bigl( \mathbh{1}_{
[2, \infty)
}(i) \cdot
\bigl[ \mathbf{x}_i^{ (n) } - \mathbf{x}_{ i - 1 }^{ (n) }
\bigr] \\
&&\hspace*{141pt}{}+ \mathbh{1}_{
[0, k-1]
}(i) \cdot\bigl[ \mathbf{x}_i^{ (n) }
- \mathbf{x}_{ i + 1 }^{ (n) } \bigr] \bigr), n \mathbf{
A}_i^{ (n) } \bigr) \Biggr] \\
&\leq&0.
\end{eqnarray*}
This contradiction implies that
$ S_0 \leq0 $.
As $ \delta\in(0,1] $ was arbitrary,
we conclude that
$
\sum_{ i = 1 }^k
u_i(t,x) \leq0
$
for all
$
(t,x) \in[0,T) \times O
$.
This finishes the proof of
Lemma~\ref{ltechnicallemmauniqueness}.
\end{pf}

The next result, Corollary~\ref{ctechnicallemmauniqueness},
establishes a comparison result for certain viscosity subsolutions
and certain viscosity supersolutions of a PDE.
It is a direct consequence of
Lemma~\ref{ltechnicallemmauniqueness}
above in the case $ k = 2 $.
Corollary~\ref{ctechnicallemmauniqueness} essentially generalizes
Theorem 2.4 in Appendix C in Peng~\cite{Peng2010}
(which assumes the function $ G $ to be globally Lipschitz continuous
in the third and last argument uniformly in the remaining arguments)
and essentially generalizes Theorem 8.2 in
Crandall, Ishii and Lions~\cite{CrandallIshiiLions1992} (which assumes a
bounded domain and a globally uniform
estimate on the function $ G $).
Corollary~\ref{ctechnicallemmauniqueness} is an immediate consequence
of Lemma~\ref{lemsignchanges} and
Lemma~\ref{ltechnicallemmauniqueness} with $ k = 2 $.
Its proof is therefore omitted.

%co4.11 #&#
\begin{corollary}[(A comparison result
for viscosity subsolutions and viscosity supersolutions)]
\label{ctechnicallemmauniqueness}
Let $ T \in(0,\infty) $,
$ d \in\N$,
let
$ O \subset\R^d $
be an open set,
let
$
u_1, u_2 \in
C( [0,T]\times O, \R)
$,
let
$
G
\colonn
(0,T)
\times
O
\times
\R
\times
\R^d
\times\mathbb{S}_d
\to\R
$
be a degenerate elliptic and continuous function
and assume that
$ u_1|_{ (0,T) \times O } $
is
a viscosity subsolution of
%
%e4.51 #&#
\begin{equation}
\label{eqparabolicequationiCOR} \frac{ \partial}{ \partial t }
u(t,x) - G \bigl( t, x, u(t,x), (
\nabla_x u) (t,x), (\operatorname{Hess}_x u) (t,x) \bigr)
= 0
\end{equation}
for
$
(t,x) \in(0,T) \times O
$
and that $ u_2|_{ (0,T) \times O } $
is a viscosity supersolution of
\eqref{eqparabolicequationiCOR}.
Moreover, assume that
%
%e4.52 #&#
\begin{eqnarray}
\label{eq00assumptionCOR} &&\limsup_{
n \to\infty
} \bigl[ G \bigl( t_n
, x_n, r_n, n ( x_n - \hat{x}_n
), n A_n \bigr)
\nonumber
\\[-8pt]
\\[-8pt]
\nonumber
&&\hspace*{35pt}{} - G \bigl( \hat{t}_n,
\hat{x}_n, \hat{r}_n, n ( x_n -
\hat{x}_n ), n \hat{A}_n \bigr) \bigr] \leq0
\end{eqnarray}
for all
$
(
t_n, x_n, r_n, A_n
)
,
(
\hat{t}_n,
\hat{x}_n,
\hat{r}_n,
\hat{A}_n
)
\in
(0,T) \times O \times\R\times\mathbb{S}_d
$,
$ n \in\N$,
satisfying that
$
\lim_{ n \to\infty}
( t_n, x_n )
\in( 0, T ) \times O
$,
that
$
\lim_{ n \to\infty}
(
\sqrt{ n }
\|
(
t_n,
x_n
)
-
(
\hat{t}_n,
\hat{x}_n
)
\|
)
= 0
$,
that
$
0 <
\lim_{ n \to\infty}
(
r_n - \hat{r}_n
)
\leq
\sup_{ n \in\N}
(
| r_n |
+
| \hat{r}_n |
)
< \infty
$
and that
$
\forall
n \in\N,
z, \hat{z}
\in\R^d
\colonn
\langle z, A_n z
\rangle
-
\langle\hat{z}, \hat{A}_n \hat{z}
\rangle
\leq
5
\| z - \hat{z}
\|^2
$.
Furthermore, assume that
$
u_1(0,x)
\leq
u_2(0,x)
$
for all
$ x \in O $
and that
%
%e4.53 #&#
\begin{equation}
\label{eqattainsmaximumCOR} \lim_{ n \to\infty} \Bigl[ \sup_{
(t,x) \in
(0,T) \times
O_n^c
% \left\{ y \in O \colonn
% \dist( y, \R^d \setminus O )
% < \frac{ 1 }{ n }
% \mbox{ or }
% \| y \| > n
% \right\}
}
\bigl( u_1(t,x) - u_2(t,x) \bigr) \Bigr] \leq0.
\end{equation}
Then
$
u_1 \leq u_2
$, that is,
it holds that
$
u_1(t, x)
\leq
u_2( t, x )
$
for all
$
(t,x) \in[0,T] \times O
$.
\end{corollary}

Assumption~\eqref{eqattainsmaximumCOR}
in Corollary~\ref{ctechnicallemmauniqueness}
is in several cases difficult
to verify.
Lemma~\ref{lcomparisonviscositysolution} below
gives an extension of
Corollary~\ref{ctechnicallemmauniqueness}
which postulates a less restrictive
condition than \eqref{eqattainsmaximumCOR}
by using a suitable Lyapunov type function
[cf.~\eqref{eqattainsmaximumdifference2comparison}
in Lemma~\ref{lcomparisonviscositysolution}
and
\eqref{eqattainsmaximumCOR}
in Corollary~\ref{ctechnicallemmauniqueness}].
In the proof of Lemma~\ref{lcomparisonviscositysolution},
the following elementary
lemma is used.

%le4.12 #&#
\begin{lemma}[(Scaling
of viscosity subsolutions and viscosity supersolutions)]
\label{lquotientsubsolutionsmooth}
Let $ T \in(0,\infty) $,
$ d \in\N$,
let $ O \subset\R^d $
be an open set,
let
$
V \in
\C^2(
(0,T) \times O, (0,\infty)
)
$,
let
$
G \colonn
(0,T)
\times
O
\times\R\times\R^d
\times\mathbb{S}_d
\to\R
$
be a degenerate
elliptic function,
let
$ u \colonn(0,T) \times O \to\R$
be a viscosity subsolution (supersolution)
of \eqref{eqparabolicequationiCOR}
and let
$
\tilde{G} \colonn
(0,T)
\times
O
\times\R\times\R^d
\times\mathbb{S}_d
\to\R
$
be a function defined by
%
%e4.54 #&#
\begin{eqnarray}
\label{eqGtilde} && \tilde{G}(t,x,r,p,A)\nonumber\\
&&\qquad:= %\\ &
\frac{ 1 }{ V(t,x) } G
\bigl( t, x, r V(t,x), p V(t,x) + r (\nabla_x V) (t,x),
 A V(t,x) \nonumber\\
&&\hspace*{85pt}
{}+ p \bigl[ ( \nabla_x V ) (t,x) \bigr]^{*} + (
\nabla_x V ) (t,x) p^{*} \\
&&\hspace*{216pt}{}+ r ( \operatorname{Hess}_x
V) (t,x) \bigr) \nonumber\\
&&\qquad\quad{}- r \frac{
({ \partial}/{
\partial t
})
V(t,x)
}{
V(t,x)
}
\nonumber
\end{eqnarray}
for all
$
(t, x, r, p, A)
\in
(0,T) \times
O \times\R\times\R^d
\times\mathbb{S}_d
$.
Then
$ \tilde{G} $ is
degenerate elliptic and
the function
$
\tilde{u} \colonn(0,T)
\times O \to\R
$
defined by
$
\tilde{u}(t,x)
=
\frac{ u(t,x) }{ V(t,x) }
$
for all
$ (t,x) \in(0,T) \times O $
is a viscosity subsolution (supersolution) of
%
%e4.55 #&#
\begin{equation}
\label{eqparabolicequationtilde} \frac{ \partial}{
\partial t
} \tilde{u}( t, x ) - \tilde{G} \bigl( t, x,
\tilde{u}(t,x), (\nabla_x \tilde{u}) (t,x), (\operatorname{Hess}_x
\tilde{u}) (t,x) \bigr) = 0
\end{equation}
for $ (t,x) \in(0,T) \times O $.
\end{lemma}

\begin{pf}%{Proof
%of
%Lemma~\ref{lquotientsubsolutionsmooth}}
We proof Lemma~\ref{lquotientsubsolutionsmooth}
in the case where $ u $ is a viscosity subsolution
of~\eqref{eqparabolicequationiCOR}.
The case where $ u $ is a viscosity supersolution
of~\eqref{eqparabolicequationiCOR} follows
analogously.
We thus assume in the following that $ u $ is a viscosity
subsolution of~\eqref{eqparabolicequationiCOR}.
First, observe that
$
\tilde{u}
$
is upper
semicontinuous
and that
$
\tilde{G}
$
is degenerate elliptic.
In the next step assume
that there exist
a vector
$ (t,x) \in(0,T) \times O $
and
a function
$
\phi\in
\C^2( (0,T) \times O, \R)
$
satisfying
$
\phi( t, x ) =
\tilde{u}(t,x)
$
and
$
\phi\geq\tilde{u}
$.
Then the function
$
(0,T) \times O
\ni
(s,y) \mapsto
\phi(s,y) V(s,y) \in\R
$
is in
$
\C^{ 2 }( (0,T) \times O, \R)
$
and satisfies
$
\phi(t, x) V(t,x) =
\tilde{u}(t,x) V(t,x) =
u(t,x)
$
and
$
\phi\cdot V \geq
\tilde{u} \cdot V = u
$.
As $u$ is a viscosity
subsolution
of~\eqref{eqparabolicequationiCOR},
we get
%
%e4.56 #&#
\begin{eqnarray}
&& V(t,x) \cdot\frac{ \partial}{ \partial t
} \phi(t,x) + \phi(t,x) \cdot\frac{ \partial}{ \partial t }
V(t,x) % \\ & =
% \frac{ \partial}{ \partial t }
% \big(
% \phi(t,x) V(t,x)
% \big)
\nonumber
\\[-8pt]
\\[-8pt]
\nonumber
&&\qquad\leq G
\bigl( t, x, \phi(t,x) V(t,x), \bigl( \nabla_x( \phi V) \bigr) (t,x),
\bigl( \operatorname{Hess}_x( \phi V) \bigr) (t,x) \bigr).
\end{eqnarray}
Rearranging this inequality
results in
%
%e4.57 #&#
%e4.58 #&#
%e4.59 #&#
%e4.60 #&#
\begin{eqnarray}
\label{eqGtildephi}\nonumber \frac{ \partial
}{ \partial t
} \phi(t,x) %\\ &
&\leq&
\frac{ 1 }{ V(t,x) } G \bigl( t, x, \phi(t,x) V(t,x), \bigl( \nabla_x(
\phi V ) \bigr) (t,x),\\
&&\hspace*{126pt}{} \bigl( \operatorname{Hess}_x( \phi V ) \bigr)
(t,x) \bigr)\nonumber\\
&&{} - \phi(t,x) \frac{
({ \partial}/{ \partial t }) V(t,x)
}{
V(t,x)
}
\\
& =& \frac{ 1 }{ V(t,x) } G \bigl( t, x, \phi(t,x) V(t,x), (\nabla_x \phi)
(t,x) V(t,x)\nonumber \\
&&\hspace*{48pt}{}+ \phi(t,x) (\nabla_x V) (t,x), ( \operatorname{Hess}_x
\phi) (t,x) V(t,x)
\nonumber\\
&&\hspace*{48pt}{} + ( \nabla_x \phi) (t,x) \bigl[ (\nabla_x V) (t,x)
\bigr]^{*}\nonumber\\
&&{}\hspace*{48pt} + ( \nabla_x V ) (t,x) \bigl[ (
\nabla_x \phi) (t,x) \bigr]^{*} \nonumber\\
&&\hspace*{136pt}{}+ \phi(t,x) (
\operatorname{Hess}_x V) (t,x) \bigr)\nonumber\\
&&{} - \phi(t,x) \frac{
({ \partial}/{ \partial t })
V(t,x)
}{
V(t,x)
}
\nonumber\\
& = &\tilde{G} \bigl( t, x, \phi(t,x), (\nabla_x \phi) (t,x), (
\operatorname{Hess}_x \phi) (t,x) \bigr).\nonumber
\end{eqnarray}
This proves
inequality~\eqref{eqGtildephi}
for all
$
\phi\in
\{
\psi\in
\C^{ 2 }( (0,T) \times O, \R)
\colonn
\psi(t,x) =
\tilde{u}(t,x)
\mbox{ and }
\psi\geq\tilde{u}
\}
$
and all
$ (t,x) \in(0,T) \times O $.
Therefore,
$ \tilde{u} $ is a viscosity
subsolution
of~\eqref{eqparabolicequationtilde}
and the proof of
Lemma~\ref{lquotientsubsolutionsmooth}
is completed.
\end{pf}

%le4.13 #&#
\begin{lemma}[(A further comparison result for viscosity
subsolutions and viscosity supersolutions)]
\label{lcomparisonviscositysolution}
Let
$ T \in(0,\infty) $,
$ d \in\N$,
let
$ O \subset\R^d $
be an open set,
let
$
u_1, u_2 \in
C( [0,T] \times O, \R)
$,
$
V \in
C( [0,T] \times O, ( 0, \infty) )
$,
% with
% $
% V|_{ (0,T) \times O }
% \in C^2( (0,T) \times\R, ( 0, \infty) )
% $,
let
$
G
\colonn
(0,T)
\times
O
\times
\R
\times
\R^d
\times
\mathbb{S}_d
\to
\R
$
be a degenerate elliptic and continuous function
and assume that
$
u_1|_{ (0,T) \times O }
$
is a viscosity subsolution of
%
%e4.61 #&#
\begin{equation}
\label{eqparabolicequationunboundedcase} \frac{ \partial}{
\partial t } u(t,x) - G \bigl( t, x, u(t,x), (
\nabla_x u) (t,x), (\operatorname{Hess}_xu) (t,x) \bigr) =
0
\end{equation}
for
$ (t,x) \in(0,T) \times O $,
that
$
u_2|_{ (0,T) \times O }
$
is a viscosity supersolution
of~\eqref{eqparabolicequationunboundedcase}
and that
for every $ r \in( 0, \infty) $
it holds that
$ r V|_{(0,T)\times O}\in C^2 ((0,T)\times O,(0,\infty) ) $
is a classical supersolution
of~\eqref{eqparabolicequationunboundedcase}.
Moreover, assume that
%
%e4.62 #&#
\begin{eqnarray}
\label{eqcomparisonviscositysolutionassumption}&& \limsup_{ n \to
\infty} \biggl( \frac{
G( t_n, x_n, r_n, p_n, A_n + n B_n V( t_n, x_n ) )
}{
V( t_n, x_n )
} \nonumber\\
&&\hspace*{16pt}\qquad{}-
\frac{
G( \hat{t}_n, \hat{x}_n, \hat{r}_n, \hat{p}_n, \hat{A}_n + n \hat
{B}_n V( \hat{t}_n, \hat{x}_n ) )
}{
V( \hat{t}_n, \hat{x}_n )
} \biggr) \\
&&\qquad\leq\frac{
G( t_0, x_0, r_0, p_0, A_0 )
}{
V( t_0, x_0 )
}\nonumber
\end{eqnarray}
for all
$
( t_n, x_n, r_n, p_n, A_n, B_n ),
( \hat{t}_n, \hat{x}_n, \hat{r}_n, \hat{p}_n, \hat{A}_n, \hat
{B}_n )
\in
(0,T)
\times O
\times\R
\times\R^d
\times\mathbb{S}_d
\times\mathbb{S}_d
$,
$
n \in\N_0
$,
satisfying that
$
\lim_{ n \to\infty}(t_n,x_n)
=
(t_0, x_0)
$,
that\break 
$
\lim_{ n \to\infty}
(
\sqrt{ n }
\|
( t_n, x_n )
-
( \hat{t}_n, \hat{x}_n )
\|
)
= 0
$,
that
$
0 < r_0
=
\lim_{ n \to\infty}
( r_n - \hat{r}_n )
\leq\break 
\sup_{ n \in\N}
(
|r_n|
+
|\hat{r}_n|
)
< \infty
$,
that
$
\lim_{ n \to\infty} ( p_n - \hat{p}_n ) = p_0
$,
that
$
\lim_{ n \to\infty} ( A_n - \hat{A}_n ) = A_0
$,
that
$
\lim_{ n \to\infty}
(
n^{ - 1 / 2 }
[
\|
\hat{p}_n
\|
+
\|
\hat{A}_n
\|_{ L( \R^d ) }
]
)
= 0
$
and that
$
\forall n \in\N, z, \hat{z} \in\R^d
\colonn
\langle z, B_n z \rangle
-
\langle\hat{z}, \hat{B}_n \hat{z} \rangle
\leq
5 \llVert z - \hat{z} \rrVert^2
$.
Furthermore, assume that
$
u_1(0, x) \leq u_2(0,x)
$
for all $ x \in O $
and that
%
%e4.63 #&#
\begin{equation}
\label{eqattainsmaximumdifference2comparison} \lim_{ n \to\infty
} \biggl[ \sup_{
{
x \in O_n^c }
}
\sup_{
t \in(0,T)
} \frac{ ( u_1(t, x) - u_2(t,x) ) }{ V(t,x) } \biggr] \leq0.
\end{equation}
Then
$ u_1 \leq u_2 $, that is,
it holds that
$ u_1(t,x) \leq u_2(t,x) $
for all
$ (t,x) \in[0,T] \times O $.
\end{lemma}

\begin{pf}%{Proof
%of Lemma~\ref{lcomparisonviscositysolution}}
Define functions
$
\tilde{u}_1,
\tilde{u}_2
\colonn[0,T] \times O \to\R
$
and
$
\tilde{G}
\colonn
(0,T)
\times
O
\times\R\times\R^d
\times\mathbb{S}_d
\to\R
$
by
$
\tilde{u}_1(t,x) =
\frac{ u_1(t,x) }{ V(t,x) }
$
and
$
\tilde{u}_2(t,x) =
\frac{ u_2(t,x) }{ V(t,x) }
$
for all
$
(t,x) \in[0,T] \times O
$
and by
%
%e4.64 #&#
\begin{eqnarray}
\label{eqGtilde1}&& \tilde{G}( t, x, r, p, A)\nonumber\\
&&\qquad:= \frac{ 1 }{ V(t,x)
} G \bigl( t, x,
r V(t,x), p V(t,x) + r (\nabla_x V) (t,x),A V(t,x)
\nonumber
\\
&&{}\hspace*{59pt}\qquad  + p \bigl[ ( \nabla_x V ) (t,x) \bigr]^{*} + (
\nabla_x V ) (t,x) p^{*}\\
&&\hspace*{192pt}\qquad{} + r ( \operatorname{Hess}_x
V) (t,x) \bigr)\nonumber\\
&&\quad\qquad{} - r \frac{
({ \partial}/{
\partial t
})
V(t,x)
}{
V(t,x)
}
\nonumber
\end{eqnarray}
for all
$
(t, x, r, p, A)
\in
(0,T) \times
O \times\R\times\R^d
\times\mathbb{S}_d
$.
% Then
% $ \tilde{G} $ is continuous
% and
Lemma~\ref{lquotientsubsolutionsmooth}
then
ensures that
$ \tilde{G} $
is degenerate elliptic,
that
$
\tilde{u}_1|_{ (0,T) \times O }
$
is a viscosity subsolution of
%
%e4.65 #&#
\begin{equation}
\label{eqparabolicequationunboundedcaseiPROOF} \frac{ \partial
}{ \partial t } u(t,x) - \tilde{G} \bigl( t, x, u(t,x), (
\nabla_x u) ( t, x ), (\operatorname{Hess}_x u) (t,x)
\bigr) = 0
\end{equation}
for
$
(t,x) \in(0,T) \times O
$
and that
$
\tilde{u}_2|_{ (0,T) \times O }
$
is viscosity supersolution of
\eqref{eqparabolicequationunboundedcaseiPROOF}.
Below we will finish this proof
by an application
of Corollary~\ref{ctechnicallemmauniqueness}
with $ \tilde{u}_1 $, $ \tilde{u}_2 $
and $ \tilde{G} $.
For this, we now check the assumptions
of Corollary~\ref{ctechnicallemmauniqueness}.
First, observe that
assumption~\eqref{eqattainsmaximumdifference2comparison}
ensures that
%implies
% \begin{equation} \begin{eqnarray}
% &
% \lim_{ n \to\infty}
% \Bigg[
% \sup_{
% \substack{
% (t,x) \in
% (0,T)\times\left\{ y \in O \colonn
% \dist( y, \R^d \setminus O )
% < \frac{ 1 }{ n }
% \mbox{ or }
% \| y \| > n
% \right\}
% }
% }
% \Big(
% \tilde{u}_1(t,x)
% -
% \tilde{u}_2(t,x)
% \Big)
% \Bigg]
% % \\ &
% % =
% % \lim_{ n \to\infty}
% % \Bigg[
% % \sup_{
% % \substack{
% % x \in
% % \left\{ y \in O \colonn
% % \dist( y, \R^d \setminus O )
% % < \frac{ 1 }{ n }
% % \mbox{ or }
% % \| y \| > n
% % \right\}
% % }
% % }
% % \sup_{ t \in(0,T) }
% % \frac{\big(u_1(t,x)-u_2(t,x)\big)}{V(t,x)}
% % \Bigg]
% \leq0
% \end{eqnarray}
% \end{equation}
% and this proves~
\eqref{eqattainsmaximumCOR} is fulfilled.
In addition, observe that
the assumption
$
u_1(0,x) \leq u_2(0,x)
$
for all $ x \in O $
ensures that
$
\tilde{u}_1(0,x)
\leq
\tilde{u}_2(0,x)
$
for all
$ x \in O $.
In the next step,
we verify \eqref{eq00assumptionCOR}.
For this, let
$
(
t_n, x_n, r_n, A_n
)
,
(
\hat{t}_n,
\hat{x}_n,
\hat{r}_n,
\hat{A}_n
)
\in
(0,T) \times O \times\R\times\mathbb{S}_d
$,
$ n \in\N_0 $,
be sequences
satisfying that
$
\lim_{ n \to\infty}
( t_n, x_n )
=
( t_0, x_0 )
=
( \hat{t}_0, \hat{x}_0 )
\in( 0, T ) \times O
$,
that
$
\lim_{ n \to\infty}
(
\sqrt{ n }
\|
( t_n, x_n )
-
( \hat{t}_n, \hat{x}_n )
\|
)
= 0
$,
that
$
0 <
r_0
=
\hat{r}_0
=
\lim_{ n \to\infty}
(
r_n - \hat{r}_n
)
\leq
\sup_{ n \in\N}
(
| r_n |
+
| \hat{r}_n |
)
< \infty
$
and that
$
\forall
n \in\N,
z, \hat{z}
\in\R^d
\colonn
\langle z, A_n z
\rangle
-
\langle\hat{z}, \hat{A}_n \hat{z}
\rangle
\leq
5
\| z - \hat{z}
\|^2
$.
To verify \eqref{eq00assumptionCOR},
we will apply assumption~\eqref{eqcomparisonviscositysolutionassumption}.
For this, we define
$ \tilde{V} \colonn[0,T] \times O \to(0,\infty) $
and
$
(
\mathbf{t}_n,
\mathbf{x}_n,
\mathbf{r}_n,
\mathbf{p}_n,
\mathbf{A}_n,
\mathbf{B}_n
)
,
(
\hat{\mathbf{t}}_n,
\hat{\mathbf{x}}_n,
\hat{\mathbf{r}}_n,
\hat{\mathbf{p}}_n,
\hat{\mathbf{A}}_n,\break 
\hat{\mathbf{B}}_n
)
\in
(0,T) \times
O \times
\R\times
\R^d \times
\mathbb{S}_d
$,
$ n \in\N_0 $,
by
$
\tilde{V}( t, x ) =
r_0 \cdot V( t, x )
$
for all
$
( t, x ) \in[0,T] \times O
$
and by
$
(\mathbf{t}_n,
\mathbf{x}_n,
\mathbf{r}_n
)
:= (t_n,
x_n,
r_n V( t_n, x_n )
)
$,
$
(\hat{\mathbf{t}}_n,
\hat{\mathbf{x}}_n,
\hat{\mathbf{r}}_n
)
:= (\hat{t}_n,
\hat{x}_n,
\hat{r}_n V( \hat{t}_n, \hat{x}_n )
)
$,
% $
% \hat{\mathbf{t}}_n
%:=
% \hat{t}_n
% $,
% $
% \mathbf{x}_n:= x_n
% $,
% $
% \hat{\mathbf{x}}_n
%:=
% \hat{x}_n
% $,
% $
% \mathbf{r}_n
%:=
% r_n V( t_n, x_n )
% $,
% $
% \hat{\mathbf{r}}_n
%:=
% \hat{r}_n
% V( \hat{t}_n, \hat{x}_n )
% $,
$
\mathbf{B}_n:=
A_n
$,
$
\hat{\mathbf{B}}_n:=
\hat{A}_n
$,
%
%e4.66 #&#
%e4.67 #&#
%e4.68 #&#
%e4.69 #&#
\begin{eqnarray}
\mathbf{p}_n &:=& n ( x_n - \hat{x}_n ) V(
t_n, x_n ) + r_n ( \nabla_x V )
(t_n, x_n ),
\\
\hat{\mathbf{p}}_n &:=& n ( x_n - \hat{x}_n ) V(
\hat{t}_n, \hat{x}_n ) + \hat{r}_n (
\nabla_x V ) ( \hat{t}_n, \hat{x}_n ),
\\
\mathbf{A}_n &:=& n ( x_n - \hat{x}_n ) \bigl[
( \nabla_x V ) ( t_n, x_n )
\bigr]^{*} + ( \nabla_x V ) ( t_n,
x_n ) n ( x_n - \hat{x}_n )^{*}
\nonumber
\\[-8pt]
\\[-8pt]
\nonumber
&&{} +
r_n ( \operatorname{Hess}_x V) ( t_n,
x_n ),
\\
\hat{\mathbf{A}}_n &:= &n ( x_n - \hat{x}_n )
\bigl[ ( \nabla_x V ) ( \hat{t}_n, \hat{x}_n
) \bigr]^{*} + ( \nabla_x V ) ( \hat{t}_n,
\hat{x}_n ) n ( x_n - \hat{x}_n
)^{*}
\nonumber
\\[-8pt]
\\[-8pt]
\nonumber
&&{} + \hat{r}_n ( \operatorname{Hess}_x V) (
\hat{t}_n, \hat{x}_n )
\end{eqnarray}
for all $ n \in\N_0 $.
Continuity of $ V $
and
%the fact
$
0 <
r_0 =
\lim_{ n \to\infty}
(
r_n - \hat{r}_n
)
\leq
\sup_{ n \in\N}
(
| r_n |
+
| \hat{r}_n |
)
< \infty
$
then imply that
%
%e4.70 #&#
\begin{eqnarray}
\label{eqlimr} 0 < \mathbf{r}_0 & =& r_0 V(
t_0, x_0 ) % =
% \lim_{ n \to\infty}
% \left(
% (
% r_n - \hat{r}_n
% )
% V( t_n, x_n )
% \right)
= \lim
_{ n \to\infty} \bigl( r_n V( t_n,
x_n ) - \hat{r}_n V( \hat{t}_n,
\hat{x}_n ) \bigr) % \\ &
\nonumber\\
&=& \lim_{ n \to\infty} ( \mathbf{
r}_n - \hat{\mathbf{r}}_n ) \\
&\leq&\sup_{ n \in\N}
\bigl( | \mathbf{r}_n | + | \hat{\mathbf{r}}_n | \bigr) < \infty.\nonumber
\end{eqnarray}
Moreover, note that the local Lipschitz continuity of
$ V $ and $ \nabla_x V $ and the continuity of $ \operatorname
{Hess}_x V $
together with the assumptions
$
\lim_{ n \to\infty}
(
\sqrt{ n }
\|
( t_n, x_n )
-
( \hat{t}_n, \hat{x}_n )
\|
)
=
\lim_{ n \to\infty}
(
\sqrt{ n }
\|
x_n
-
\hat{x}_n
\|
)
= 0
$,
$
\lim_{ n \to\infty}
(
r_n - \hat{r}_n
)
= r_0
$
and\break 
$
\sup_{ n \in\N}
| \hat{r}_n | < \infty
$
imply that
%
%e4.71 #&#
%e4.72 #&#
%e4.73 #&#
%e4.74 #&#
%e4.75 #&#
\begin{eqnarray}
\label{eqlimp} \lim_{ n \to\infty} ( \mathbf{p}_n - \mathbf{
\hat{p}}_n ) & =& \lim_{ n \to\infty} \bigl[ n (
x_n - \hat{x}_n ) \bigl( V( t_n,
x_n ) - V( \hat{t}_n, \hat{x}_n ) \bigr)
\bigr]\nonumber\\
&&{} + \lim_{ n \to\infty} \bigl[ ( r_n -
\hat{r}_n ) ( \nabla_x V ) ( t_n,
x_n ) \bigr]
\nonumber
\\[-8pt]
\\[-8pt]
\nonumber
&&{} + \lim_{ n \to\infty} \bigl[ \hat{r}_n \bigl( (
\nabla_x V ) ( t_n, x_n ) - (
\nabla_x V ) ( \hat{t}_n, \hat{x}_n ) \bigr)
\bigr] \\
&= & r_0 ( \nabla_x V ) ( t_0,
x_0 ) = \mathbf{p}_0,\nonumber\\
\label{eqlimAhatA}  \lim_{ n \to\infty} ( \mathbf{A}_n -
\mathbf{
\hat{A}}_n ) &=& \lim_{ n \to\infty} \bigl( n (
x_n - \hat{x}_n ) \bigl( \bigl[ ( \nabla_x V
) ( t_n, x_n ) \bigr]^{*} - \bigl[ (
\nabla_x V ) ( \hat{t}_n, \hat{x}_n )
\bigr]^{*} \bigr) \bigr)
\nonumber\\
&&{} + \lim_{ n \to\infty} \bigl( \bigl[ ( \nabla_x V ) (
t_n, x_n ) - ( \nabla_x V ) (
\hat{t}_n, \hat{x}_n ) \bigr] n ( x_n -
\hat{x}_n )^* \bigr) \nonumber\\
&&{}+ \lim_{ n \to\infty} \bigl( [
r_n - \hat{r}_n ] ( \operatorname{Hess}_x V)
( t_n, x_n ) \bigr)
\\
& &{}+ \lim_{ n \to\infty} \bigl( \hat{r}_n \bigl[ (
\operatorname{Hess}_x V) ( t_n, x_n ) - (
\operatorname{Hess}_x V) ( \hat{t}_n,
\hat{x}_n ) \bigr] \bigr) \nonumber\\
&=& r_0 ( \operatorname{Hess}_x
v ) ( t_0, x_0 ) = \mathbf{A}_0\nonumber
\end{eqnarray}
and
$
\lim_{ n \to\infty}
(
n^{ - 1 / 2 }
[
\| \hat{\mathbf{p}}_n \|
+
\| \hat{\mathbf{A}}_n \|_{ L( \R^d ) }
]
)
= 0
$.
Combining this and \eqref{eqlimr}
%and \eqref{eqlimAhatA}
with
assumption \eqref{eqcomparisonviscositysolutionassumption}
shows that
%
%e4.76 #&#
\begin{eqnarray}
&&\limsup_{ n \to\infty} \biggl( \frac{
G(
\mathbf{t}_n,
\mathbf{x}_n,
\mathbf{r}_n,
\mathbf{p}_n,
\mathbf{A}_n
+ n \mathbf{B}_n
V( \mathbf{t}_n, \mathbf{x}_n )
)
}{
V( \mathbf{t}_n, \mathbf{x}_n )
} \nonumber\\
&&\hspace*{36pt}{}-
\frac{
G(
\hat{\mathbf{t}}_n,
\hat{\mathbf{x}}_n,
\hat{\mathbf{r}}_n,
\hat{\mathbf{p}}_n,
\hat{\mathbf{A}}_n
+ n \hat{\mathbf{B}}_n
V( \hat{\mathbf{t}}_n, \hat{\mathbf{x}}_n )
)
}{
V( \hat{\mathbf{t}}_n, \hat{\mathbf{x}}_n )
} \biggr) \\
&&\qquad\leq\frac{
G( \mathbf{t}_0, \mathbf{x}_0, \mathbf{r}_0, \mathbf{p}_0, \mathbf
{A}_0 )
}{
V( \mathbf{t}_0, \mathbf{x}_0 )
}.\nonumber
\end{eqnarray}
The definition of $ \tilde{G} $ hence
implies that
%
%e4.77 #&#
%e4.78 #&#
%e4.79 #&#
\begin{eqnarray}
&& \limsup_{
n \to\infty
} \bigl( \tilde{G} \bigl( t_n,
x_n, r_n, n ( x_n - \hat{x}_n )
, n A_n \bigr) - \tilde{G} \bigl( \hat{t}_n,
\hat{x}_n, \hat{r}_n, n ( x_n -
\hat{x}_n ), n \hat{A}_n \bigr) \bigr)
\nonumber\\
&&\qquad = \limsup_{
n \to\infty
} \biggl( \frac{
G(
t_n, x_n, \mathbf{r}_n,
\mathbf{p}_n,
\mathbf{A}_n + n \mathbf{B}_n V( t_n, x_n )
)
-
r_n
({ \partial}/{ \partial t }) V( t_n, x_n )
}{
V( t_n, x_n )
}\nonumber\\
&&\hspace*{70pt}{} -
\frac{
G(
\hat{t}_n, \hat{x}_n, \hat{\mathbf{r}}_n,
\hat{\mathbf{p}}_n,
\hat{\mathbf{A}}_n
+
n \hat{\mathbf{B}}_n V( \hat{t}_n, \hat{x}_n )
)
-
\hat{r}_n
({ \partial}/{ \partial t })
V( \hat{t}_n, \hat{x}_n )
}{
V( \hat{t}_n, \hat{x}_n )
} \biggr)
\nonumber\\
&&\qquad \leq% \frac{
% G(
% t_0, x_0, \mathbf{r}_0,
% \mathbf{p}_0,
% \mathbf{A}_0
% )
% }{
% V( t_0, x_0 )
% }
% +
% \limsup_{
% n \to\infty
% }
% \left[
% \frac{
% \hat{r}_n
% \frac{ \partial}{ \partial t }
% V( \hat{t}_n, \hat{x}_n )
% }{
% V( \hat{t}_n, \hat{x}_n )
% }
% -
% \frac{
% r_n
% \frac{ \partial}{ \partial t } V( t_n, x_n )
% }{
% V( t_n, x_n )
% }
% \right]
% \\ & =
% \frac{
% G(
% t_0, x_0, \mathbf{r}_0,
% \mathbf{p}_0,
% \mathbf{A}_0
% )
% }{
% V( t_0, x_0 )
% }
% +
% \lim_{
% n \to\infty
% }
% \left[
% \frac{
% \left(
% \hat{r}_n
% -
% r_n
% \right)
% \frac{ \partial}{ \partial t }
% V( \hat{t}_n, \hat{x}_n )
% }{
% V( \hat{t}_n, \hat{x}_n )
% }
% \right]
% +
% \lim_{
% n \to\infty
% }
% \left(
% r_n
% \left[
% \frac{
% \frac{ \partial}{ \partial t }
% V( \hat{t}_n, \hat{x}_n )
% }{
% V( \hat{t}_n, \hat{x}_n )
% }
% -
% \frac{
% \frac{ \partial}{ \partial t } V( t_n, x_n )
% }{
% V( t_n, x_n )
% }
% \right]
% \right)
% \\ & =
\frac{
G(
t_0, x_0, \mathbf{r}_0,
\mathbf{p}_0,
\mathbf{A}_0
)
}{
V( t_0, x_0 )
} - \frac{
r_0 ({ \partial}/{ \partial t }) V( t_0, x_0 )
}{
V( t_0, x_0 )
} % \\ &
\nonumber
\\[-8pt]
\\[-8pt]
\nonumber
&&\qquad=
\bigl(
-
\bigl[
({ \partial}/{ \partial t })
\tilde{V}( t_0, x_0 )
\\
&&\hspace*{50pt}{}-
G\bigl(
t_0, x_0,
\tilde{V}( t_0, x_0 ),
( \nabla_x \tilde{V} )( t_0, x_0 ),
( \operatorname{Hess}_x \tilde{V} )( t_0, x_0 )
\bigr)
\bigr]
\bigr)\nonumber\\
&&\hspace*{36pt}{}/\bigl(
V( t_0, x_0 )
\bigr) \nonumber\\
&&\qquad\leq0\nonumber
\end{eqnarray}
as $ \tilde{V} $ is by assumption a classical supersolution
of \eqref{eqparabolicequationunboundedcase}.
We can thus apply Corollary~\ref{ctechnicallemmauniqueness}
to obtain that
$
\tilde{u}_1( t, x )
=
\frac{ u_1( t, x ) }{ V( t, x ) }
\leq
\frac{ u_2( t, x ) }{ V( t, x ) }
=
\tilde{u}_2( t, x )
$
for all $ (t,x) \in[0,T] \times O $.
This finishes the proof
of Lemma~\ref{lcomparisonviscositysolution}.
\end{pf}

The next result,
Corollary~\ref{coruniqueness2},
asserts uniqueness of the solution of
a linear second-order PDE.
We assume that the
Lyapunov-type function
$
V \colonn[0,T] \times O
\to(0,\infty)
$
in
Lemma~\ref{lcomparisonviscositysolution}
is of the form
$
V(t,x) = e^{ \rho t } \cdot\tilde{V}(x)
$
for all
$ (t,x) \in[0,T] \times O $
where
$ \rho\in\R$ is a real number
and where
$ \tilde{V} \colonn O \to(0,\infty) $
is a twice continuously differentiable
function.

%co4.14 #&#
\begin{corollary}[(Uniqueness
of viscosity solutions of Kolmogorov type equations)]
\label{coruniqueness2}
Let
$ T \in(0,\infty) $,
$ d, m \in\N$,
$ \rho\in\R$,
let
$ O \subset\R^d $
be an open set,
let
$
\varphi\in C( O, \R)
$,
$
v \in C( (0,T) \times O, \R)
$,
let
$
\mu\colonn
(0,T) \times O \to\R^d
$
and
$
\sigma\colonn
(0,T) \times O \to
\R^{ d \times m }
$
be locally Lipschitz continuous functions
and
let $ V \in C^2( O, (0,\infty) ) $
satisfy
%
%e4.80 #&#
\begin{eqnarray}
\label{eqassumptionV}&& v(t,x) V(x) + \bigl\anglel\mu(t,x), (\nabla
V) (x) \bigr\angler
+ \operatorname{tr} \bigl( \sigma(t,x) \bigl[ \sigma(t,x) \bigr
]^* (
\operatorname{Hess} V) (x) \bigr)
\nonumber
\\[-8pt]
\\[-8pt]
\nonumber
&&\qquad \leq\rho\cdot V(x)
\end{eqnarray}
for all $ (t,x) \in(0,T) \times O $.
Then there exists
at most one
continuous function
$
u
\colonn[0,T] \times O \to\R
$
which fulfills
$
u(0, x) = \varphi(x)
$
for all $ x \in O $,
which fulfills\
$
\lim_{ n \to\infty}
\sup_{
{
(t,x) \in(0,T) \times
O_n^c
}
}
\frac{
| u(t, x) |
}{
V(x)
}
$
$
= 0
$
and which fulfills that
$
u|_{ (0,T) \times O }
$
is a viscosity solution of
%
%e4.81 #&#
\begin{eqnarray}
\label{eqsecond-orderPDE} &&\frac{ \partial}{ \partial t } u(t,x) -
v(t,x) u(t,x) - \bigl\langle\mu(t,x),
(\nabla_x u) (t,x) \bigr\rangle\nonumber\\
&&\quad{}- \operatorname{tr} \bigl(
\sigma(t,x) \bigl[ \sigma(t,x) \bigr]^* (\operatorname{Hess}_xu) (t,x)
\bigr)
\\
&&\qquad= 0\nonumber
\end{eqnarray}
for $ (t,x) \in(0,T) \times O $.
\end{corollary}

\begin{pf}%{Proof
%of Corollary~\ref{coruniqueness2}}
Let
$
u_1,
u_2
\colonn[0,T] \times O \to\R
$
be two continuous
functions such that
$
u_1(0,x) = \varphi(x) = u_2(0,x)
$
for all
$ x \in O $,
such that
\[
\lim_{ n \to\infty}
\sup_{
{
(t,x) \in(0,T) \times
O_n^c
}
}
\frac{
| u_1(t, x) |
+
| u_2(t, x) |
}{
V(x)
}
= 0
\]
and such that
$
u_1|_{ (0,T) \times O }
$
and
$
u_2|_{ (0,T) \times O }
$
are viscosity solutions
of~\eqref{eqsecond-orderPDE}.
Then define a function
$
G \colonn(0,T) \times O \times
\R\times\R^d \times
\mathbb{S}_d
\to\R
$
by
$
G( t, x, r, p, A)
=
v(t,x) r
+
\langle
\mu(t,x), p
\rangle
+
\operatorname{tr} (
\sigma(t,x)
[ \sigma(t,x) ]^*
A
)
$.
We show Corollary~\ref{coruniqueness2} by applying
Lemma~\ref{lcomparisonviscositysolution}.
To this end we now verify
\eqref{eqcomparisonviscositysolutionassumption}.
%
% check the assumptions
% of Lemma~\ref{lcomparisonviscositysolution}.
% It follows from the equation
% \begin{equation}
% \begin{eqnarray}
% &
% \operatorname{tr}
% \Big(
% \sigma(t,x)
% [ \sigma(t,x) ]^*
% B
% \Big)
% =
% \operatorname{tr}
% \Big(
% [ \sigma(t,x) ]^*
% B
% \sigma(t,x)
% \Big)
% =
% % \sum_{ i = 1 }^m
% % \left\langle
% % e_i^{ (m) },
% % [ \sigma(t,x) ]^*
% % B
% % \sigma(t,x)
% % e_i^{ (m) }
% % \right\rangle
% % \\ & =
% \sum_{ i = 1 }^m
% \left\langle
% \sigma(t,x)
% e_i^{ (m) },
% B
% \sigma(t,x)
% e_i^{ (m) }
% \right\rangle
% \geq0
% \end{eqnarray}
% \end{equation}
% for all
% $ (t,x,B) \in(0,T) \times O \times\mathbb{S}_d $ with $ B \geq0 $
% that $ G $ is degenerate elliptic.
% HIER: This implies that
% \begin{equation} \begin{eqnarray}
% \limsup_{n\to\infty}\mbox{tr}(A(t_n,x_n)(A_n-\hat{A}_n))\leq
% \mbox{tr}(A(t_0,x_0)A_0)
% \end{eqnarray} \end{equation}
% for all TODO.
% In the next step we verify
For this, let
$
(t_n, x_n, r_n,\break  p_n,  A_n, B_n),
(\hat{t}_n, \hat{x}_n, \hat{r}_n,  \hat{p}_n, \hat{A}_n, \hat{B}_n )
\in
(0,T) \times O \times\R\times\R^d \times
\mathbb{S}_d \times\mathbb{S}_d
$,
$
n \in\N_0
$,
satisfy that
$
\lim_{ n \to\infty}( t_n, x_n )
=
( t_0, x_0)
$,
that
$
\lim_{ n \to\infty}
(
\sqrt{ n }
\|
( t_n, x_n )
-
( \hat{t}_n, \hat{x}_n )
\|
)
= 0
$,
that
$
0 < r_0
=
\lim_{ n \to\infty} (r_n - \hat{r}_n )
\leq
\sup_{ n \in\N}
( | r_n | + | \hat{r}_n | )
<
\infty
$,
that
$
\lim_{ n \to\infty} ( p_n - \hat{p}_n ) = p_0
$,
that
$
\lim_{ n \to\infty} ( A_n - \hat{A}_n ) = A_0
$,
that
$
\lim_{ n \to\infty}
(
n^{ - 1 / 2 }
[
\| \hat{p}_n \|
+\break 
\| \hat{A}_n \|_{
L( \R^d )
}
]
)
= 0
$
and that
$
\forall n \in\N, z, \hat{z} \in\R^d
\colonn
\anglel z, B_n z \angler
-
\langle\hat{z}, \hat{B}_n \hat{z} \rangle
\leq
5 \| z - \hat{z} \|^2
$.
Then it holds that
%
%e4.82 #&#
%e4.83 #&#
%e4.84 #&#
%e4.85 #&#
%e4.86 #&#
%e4.87 #&#
%e4.88 #&#
%e4.89 #&#
\begin{eqnarray}
& &\limsup_{ n \to\infty} \biggl( \frac{ 1 }{ V( t_n, x_n ) } G
\bigl(
t_n, x_n, r_n, p_n,
A_n + n B_n V( t_n, x_n ) \bigr)
\nonumber\\
&&\hspace*{25pt}\quad{}- \frac{ 1 }{
V( \hat{t}_n, \hat{x}_n )
} G \bigl( \hat{t}_n, \hat{x}_n,
\hat{r}_n, \hat{p}_n, \hat{A}_n + n
\hat{B}_n V( \hat{t}_n, \hat{x}_n ) \bigr)
\biggr)
\nonumber\\
&&\qquad \leq\limsup_{ n \to\infty} \biggl( \frac{
v( t_n, x_n ) r_n
}{
V( t_n, x_n )
} -
\frac{
v( \hat{t}_n, \hat{x}_n ) \hat{r}_n
}{
V( \hat{t}_n, \hat{x}_n )
} \biggr) \nonumber\\
&&\qquad\quad{}+ \limsup_{ n \to\infty} \biggl(
\frac{
\langle
\mu( t_n, x_n), p_n
\rangle
}{
V( t_n, x_n )
} - \frac{
\langle
\mu( \hat{t}_n, \hat{x}_n ),
\hat{p}_n
\rangle
}{
V( \hat{t}_n, \hat{x}_n )
} \biggr)
\nonumber\\
&&\qquad\quad{} + \limsup_{ n \to\infty} \biggl( \frac{
\operatorname{tr}(
\sigma( t_n, x_n )
[ \sigma( t_n, x_n ) ]^*
A_n
)
}{
V( t_n, x_n )
} \nonumber\\
&&\hspace*{82pt}{}-
\frac{
\operatorname{tr}(
\sigma( \hat{t}_n, \hat{x}_n )
[ \sigma( \hat{t}_n, \hat{x}_n ) ]^*
\hat{A}_n
)
}{
V( \hat{t}_n, \hat{x}_n )
} \biggr)
\nonumber\\
&&\qquad\quad{} + \limsup_{ n \to\infty} \bigl( n \bigl[ \operatorname{tr}
\bigl(
\bigl[ \sigma( t_n, x_n ) \bigr]^* B_n \sigma(
t_n, x_n ) \bigr) \nonumber\\
&&\hspace*{88pt}{}- \operatorname{tr} \bigl( \bigl[
\sigma( \hat{t}_n, \hat{x}_n ) \bigr]^*
\hat{B}_n \sigma( \hat{t}_n, \hat{x}_n )
\bigr) \bigr] \bigr)
\\
&&\qquad \leq\limsup_{ n \to\infty} \biggl( \frac{
v( t_n, x_n ) ( r_n - \hat{r}_n )
}{
V( t_n, x_n )
} \biggr)
\nonumber\\
&&\qquad\quad{} +
\limsup_{ n \to\infty} \biggl( \biggl[ \frac{
v( t_n, x_n )
}{
V( t_n, x_n )
} -
\frac{
v( \hat{t}_n, \hat{x}_n )
}{
V( \hat{t}_n, \hat{x}_n )
} \biggr] \hat{r}_n \biggr)\nonumber \\
&&\qquad\quad{}+ \limsup
_{ n \to\infty} \biggl( \frac{
\langle
\mu( t_n, x_n), p_n - \hat{p}_n
\rangle
}{
V( t_n, x_n )
} \biggr)
\nonumber\\
&&\qquad\quad{} + \limsup_{ n \to\infty} \biggl( \biggl\langle\sqrt{ n }
\biggl[
\frac{
\mu( t_n, x_n )
}{
V( t_n, x_n )
} - \frac{
\mu( \hat{t}_n, \hat{x}_n )
}{
V( \hat{t}_n, \hat{x}_n )
} \biggr], \frac{ \hat{p}_n }{ \sqrt{ n } } \biggr
\rangle\biggr)\nonumber\\
&&\quad\qquad{} + \limsup_{ n \to\infty} \biggl( \operatorname{tr}
\biggl( \frac{
\sigma( t_n, x_n )
[ \sigma( t_n, x_n ) ]^*
}{
V( t_n, x_n )
} ( A_n - \hat{A}_n ) \biggr)
\biggr)
\nonumber\\
&&\qquad\quad{} + \limsup_{ n \to\infty} \biggl( \operatorname{tr} \biggl(
\sqrt{ n
} \biggl[ \frac{
\sigma( t_n, x_n )
[ \sigma( t_n, x_n ) ]^*
}{
V( t_n, x_n )
} - \frac{
\sigma( \hat{t}_n, \hat{x}_n )
[ \sigma( \hat{t}_n, \hat{x}_n ) ]^*
}{
V( \hat{t}_n, \hat{x}_n )
} \biggr] \frac{
\hat{A}_n
}{
\sqrt{ n }
}
\biggr) \biggr)
\nonumber\\
&&\qquad\quad{} + \limsup_{ n \to\infty} \Biggl( n \sum
_{ i = 1 }^m \bigl[ \bigl\anglel\sigma(
t_n, x_n ) e^{ (m) }_i,
B_n \sigma( t_n, x_n ) e^{ (m) }_i
\bigr\angler\nonumber\\
&&\hspace*{106pt}{}- \bigl\anglel\sigma( \hat{t}_n, \hat{x}_n
) e^{ (m) }_i, \hat{B}_n \sigma(
\hat{t}_n, \hat{x}_n ) e^{ (m) }_i
\bigr\angler\bigr] \Biggr).\nonumber
\end{eqnarray}
Hence, the local Lipschitz continuity of
the functions
$
\frac{ \mu}{ V }
$ and
$
\frac{ A }{ V }
$
together with the properties of
$
(t_n, x_n, r_n, p_n, A_n, B_n)
$,
$
(\hat{t}_n, \hat{x}_n, \hat{r}_n, \hat{p}_n, \hat{A}_n, \hat{B}_n )
\in
(0,T) \times O \times\R\times\R^d \times
\mathbb{S}_d \times\mathbb{S}_d
$,
$
n \in\N_0
$,
% $
% \langle z, B_n z \rangle
% -
% \langle\hat{z}, \hat{B}_n \hat{z} \rangle
% \leq
% 5 \| z - \hat{z} \|^2
% $
% for all $ z, \hat{z} \in\R^d $
% and all $ n \in\N$
% and the boundedness of
% $ ( \hat{r}_n )_{ n \in\N} $
implies that
%
%e4.90 #&#
%e4.91 #&#
%e4.92 #&#
%e4.93 #&#
%e4.94 #&#
\begin{eqnarray}
&& \limsup_{ n \to\infty} \biggl( \frac{ 1 }{ V( t_n, x_n ) } G
\bigl(
t_n, x_n, r_n, p_n,
A_n + n B_n V( t_n, x_n ) \bigr)
\nonumber\\
&&\hspace*{36pt}{}- \frac{ 1 }{
V( \hat{t}_n, \hat{x}_n )
} G \bigl( \hat{t}_n, \hat{x}_n,
\hat{r}_n, \hat{p}_n, \hat{A}_n + n
\hat{B}_n V( \hat{t}_n, \hat{x}_n ) \bigr)
\biggr)
\nonumber\\
&&\qquad \leq\frac{
v( t_0, x_0 ) r_0
}{
V( t_0, x_0 )
} + \frac{
\langle
\mu( t_0, x_0), p_0
\rangle
}{
V( t_0, x_0 )
} + \operatorname{tr} \biggl(
\frac{
\sigma( t_0, x_0 )
[ \sigma( t_0, x_0 ) ]^*
}{
V( t_0, x_0 )
} A_0 \biggr)
\nonumber\\
&&\qquad\quad{} + \limsup_{ n \to\infty} \biggl( d \biggl[ \sqrt{ n } \biggl
\llVert
\frac{
\sigma( t_n, x_n )
[ \sigma( t_n, x_n ) ]^*
}{
V( t_n, x_n )
}
\nonumber
\\
\nonumber
&&\hspace*{116pt}{}- \frac{
\sigma( \hat{t}_n, \hat{x}_n )
[ \sigma( \hat{t}_n, \hat{x}_n ) ]^*
}{
V( \hat{t}_n, \hat{x}_n )
} \biggr\rrVert_{ L( \R^d ) } \biggr]
\frac{
\| \hat{A}_n \|_{ L( \R^d ) }
}{ \sqrt{ n } } \biggr)
\\
&&\qquad\quad{} + \limsup_{ n \to\infty} \Biggl( n \sum
_{ i = 1 }^m 5\bigl \| \sigma( t_n,
x_n ) e^{ (m) }_i - \sigma(
\hat{t}_n, \hat{x}_n ) e^{ (m) }_i
\bigr\|^2 \Biggr)
\nonumber\\
&&\qquad = \frac{
G( t_0, x_0, r_0, p_0, A_0 )
}{
V( t_0, x_0 )
} \\
&&\qquad\quad{}+ 5 \limsup_{ n \to\infty} \bigl( n \bigl\| \sigma(
t_n, x_n ) - \sigma( \hat{t}_n,
\hat{x}_n ) \bigr\|_{
HS( \R^m, \R^d )
}^2 \bigr)\nonumber\\
&&\qquad =
\frac{
G( t_0, x_0, r_0, p_0, A_0 )
}{
V( t_0, x_0 )
}.\nonumber
\end{eqnarray}
This shows assumption~\eqref{eqcomparisonviscositysolutionassumption}.
Moreover,
by assumption,
$
u_1|_{ (0,T) \times O }
$
is a viscosity subsolution
of~\eqref{eqsecond-orderPDE}
and
$
u_2|_{(0,T)\times O}
$
is a viscosity supersolution
of~\eqref{eqsecond-orderPDE}.
Furthermore,
\eqref{eqassumptionV} shows
for every $ r \in(0,\infty) $
that the function
$
(0,T)\times O \ni(t,x) \mapsto
r \cdot e^{ \rho t } \cdot V(x)
\in(0, \infty)
$
is a classical supersolution
of~\eqref{eqsecond-orderPDE}.
% At time $0$ we have $u_1(0,x)= u_2(0,x)$ for all $x\in O$.
In addition, observe that \eqref{eqattainsmaximumdifference2comparison}
follows from
$
\lim_{ n \to\infty}
\sup_{
{
(t,x) \in(0,T) \times
O_n^c
}
}\times\break 
\frac{
| u_1(t, x) |
+
| u_2(t, x) |
}{
V(x)
}
$
$
= 0
$.
% \begin{equation} \begin{eqnarray}
% &\lim_{ n \to\infty}
% \sup_{
% \substack{
% (t,x) \in(0,T) \times
% \left\{ y \in O \colonn
% \dist( y, \R^d\setminus O )
% < 1 / n
% \mbox{ or }
% \| y \| > n
% \right\}
% }
% }
% \frac{
% | u_1(t, x)
% -
% u_2(t, x) |
% }{
% e^{\rho t}V(x)
% }
% \\&
% \leq
% e^{|\rho|T}
% \lim_{ n \to\infty}
% \sup_{
% \substack{
% (t,x) \in(0,T) \times
% \left\{ y \in O \colonn
% \dist( y, \R^d\setminus O )
% < 1 / n
% \mbox{ or }
% \| y \| > n
% \right\}
% }
% }
% \frac{
% | u_1(t, x) |
% +
% | u_2(t, x) |
% }{
% V(x)
% }
% =0.
% \end{eqnarray} \end{equation}
Consequently,
Lemma~\ref{lcomparisonviscositysolution}
implies that $u_1\leq u_2$.
Repeating these arguments with $u_1$ and $u_2$ interchanged
finally shows that $u_2\leq u_1$ so that $u_1=u_2$.
This proves uniqueness and finishes the proof of
Corollary~\ref{coruniqueness2}.
\end{pf}

%%%%%%%%%%%%%%%%%%%%%%%%%%%%
%s4.4 #&#
\subsection{Viscosity solutions of
Kolmogorov equations}
\label{ssecViscositysolutionsofKolmogorovequations}

The main result of this subsection,
Theorem~\ref{thmKolmogorovviscosity}
below,
establishes that the transition
semigroup associated with a
suitable SDE
with locally Lipschitz continuous coefficients
is within a certain class of functions the unique viscosity solution of
the Kolmogorov
equation of the SDE.
To establish this result, we
first prove an
auxiliary result.

%le4.15 #&#
\begin{lemma}[(Existence of viscosity
solutions of Kolmogorov
equations with globally Lipschitz
continuous coefficients with
compact support)]
\label{lemKolmogorovviscosity2}
Let
$ d, m \in\mathbb{N} $,
let
$
( \Omega, \mathcal{F}, \P
)
$
be a probability space
with a normal filtration\break 
$
( \mathcal{F}_t )_{ t \in[0,\infty) }
$,
let
$
W \colonn
[0,\infty)
\times\Omega
\rightarrow\mathbb{R}^m
$
be a
standard
$
( \mathcal{F}_t )_{ t \in[0,\infty) }
$-Brownian motion,
let
$ O \subset\R^d $
be an open set,
let
$
\varphi\colonn O
\to\R
$
be a continuous function
and let
$ \mu\colonn O \to\R^d $
and
$
\sigma\colonn O \to
\R^{ d \times m }
$
be locally Lipschitz continuous
functions with compact support.
Then there exists
a family
$
X^x \colonn[0,\infty)
\times\Omega
\rightarrow O
$,
$ x \in O $,
of up to indistinguishability
unique adapted stochastic
processes
with continuous sample
paths
satisfying
%
%e4.95 #&#
\begin{equation}
\label{eqSDEviscosityexistence1} X^x(t) = x + \int_0^t
\mu\bigl( X^x(s) \bigr) \,ds + \int_0^t
\sigma\bigl( X^x(s) \bigr) \,dW(s)
\end{equation}
for all
$ t \in[0,\infty) $,
$\mathbb{P} $-a.s. and all
$ x \in O $
and the function
$
u \colonn(0,\infty) \times
O \to\R
$
given by
$
u(t,x) = \mathbb{E} [
\varphi( X^x(t) )
]
$
is a viscosity solution of
%
%e4.96 #&#
\begin{eqnarray}
\label{eqCDF2}&& \frac{ \partial}{ \partial t } u(t,x) - \bigl
\langle(\nabla_x u)
(t,x), \mu(x) \bigr\rangle- \frac{1}{2} \operatorname{tr} \bigl(
\sigma(x)
\bigl[ \sigma(x) \bigr]^{*} % \sigma(x)^{*}
(\operatorname{Hess}_x
u) (t,x) \bigr)
\nonumber
\\[-8pt]
\\[-8pt]
\nonumber
&&\qquad= 0
\end{eqnarray}
for $ (t,x) \in(0,\infty) \times O $.
\end{lemma}

\begin{pf}%{Proof
%of Lemma~\ref{lemKolmogorovviscosity2}}
First of all, observe that
since $ \mu$ and $ \sigma$
have compact supports,
they are globally Lipschitz continuous, so that \eqref
{eqSDEviscosityexistence1} has unique solutions.
It thus remains to show that
the function
$
u \colonn(0,\infty) \times O
\to\R
$
introduced above
is a viscosity solution
of \eqref{eqCDF2}.
Let $ U \subset O $
be a relatively compact open set
in $ O $
with the property that
$
\operatorname{supp}( \mu)
\cup
\operatorname{supp}( \sigma)
\subset U
$.
By assumption
$
\operatorname{supp}( \mu)
$
and
$
\operatorname{supp}( \sigma)
$
are compact sets,
and hence
such a set $U$ does indeed exist.
Next, let
$
\mu^{(n)}
\in\C^{ \infty}_{ \mathrm{cpt} }( O, \R^d )
$,
$ n \in\N$,
and
$
\sigma^{(n)}
\in\C^{ \infty}_{ \mathrm{cpt} }( O, \R)
$,
$ n \in\N$,
be sequences of smooth functions
satisfying
$
\lim_{ n \to\infty}
\sup_{ x \in U }
\|
\mu(x) -
\mu^{(n)}(x)
\|
=
\lim_{ n \to\infty}
\sup_{ x \in U }
\|
\sigma(x) -
\sigma^{(n)}(x)
\|_{
L( \R^m, \R^d )
}
= 0
$
and
$
\operatorname{supp} (
\mu^{ (n) }
)
\cup
\operatorname{supp} (
\sigma^{ (n) }
)
\subset U
$
for all $ n \in\N$
and denote by
$
X^{x,n}
\colonn[0,\infty) \times\Omega
\to O
$,
$ x \in O $,
$ n \in\N$,
the solutions to the corresponding SDEs.
Moreover,
let
$ \varphi_k \in\C^{ \infty}( O, \R) $,
$ k \in\N$,
be a sequence of smooth functions satisfying
$
\sup_{ x \in O_k }
\vert
\varphi(x) -
\varphi_k(x)
\vert
<\frac{1}{k}
$
for each $k\in\N$.
Now we
define functions
$
u^{ n, k } \colonn
(0,\infty) \times O
\to\R
$,
$ n, k \in\N$,
and
$
u^{ (k) } \colonn
(0,\infty) \times O
\to\R
$,
by
$
u^{n,k}(t,x)
:=
\mathbb{E} [
\varphi_k( X^{ x, n }(t) )
]
$
and
$
u^{(k)}(t,x)
:=
\mathbb{E} [
\varphi_k( X^{ x }(t) )
]
$.
For any fixed $n$ and $k$, the function
$
u^{ n, k } \colonn
(0,\infty) \times O
\to\R
$,
is smooth
and globally Lipschitz continuous
(see, e.g., Corollary 2.8.1 and Theorem 2.8.1 in~\cite{GihmanSkorohod1972}).
Theorem 4.3 in~\cite{PardouxPeng1992} then shows that
%
%e4.97 #&#
\begin{eqnarray}
\label{eqPDEnice}&& \biggl( \frac{ \partial}{ \partial t } u^{ n, k
} \biggr) (t,x) - \bigl
\langle\bigl( \nabla_x u^{ n, k } \bigr) (t,x),
\mu^{(n)}(x) \bigr\rangle\nonumber\\
&&\quad{}- \frac{1}{2} \operatorname{tr} \bigl(
\sigma^{(n)}(x) \bigl[ \sigma^{(n)}(x) \bigr]^{*}
\bigl( \operatorname{Hess}_x u^{ n, k } \bigr) (t,x) \bigr) \\
&&\qquad= 0\nonumber
\end{eqnarray}
for all
$
(t,x) \in(0,\infty)
\times O
$,
$ n, k \in\N$.
Remark~\ref{remclassicalsolution}
hence shows that
the functions
$
u^{ n, k }
$,
$ n, k \in\N$,
are also viscosity solutions
to these equations.
Furthermore, observe that
the smoothness of the functions
$
\varphi_k \in\C^{ \infty}( O, \R)
$,
$ k \in\N$,
and the global Lipschitz continuity
of the functions
$ ( \mu^{ (n) } )_{ n \in\N} $,
$ ( \sigma^{ (n) } )_{ n \in\N} $,
$ \mu$
and
$ \sigma$
imply that
%
%e4.98 #&#
%e4.99 #&#
%e4.100 #&#
\begin{eqnarray}
\label{eqlimitn}&& \lim_{ n \to\infty} \sup_{ t \in(0,T] } \sup
_{ x \in O } \bigl| u^{ (k) }(t, x) - u^{ n, k }(t, x) \bigr|\nonumber \\
&&\qquad=
\lim_{ n \to\infty} \sup_{ t \in(0,T] } \sup
_{ x \in\bar{U} } \bigl| \E\bigl[ \varphi_k\bigl(
X^{x,n}(t) \bigr) \bigr] - \E\bigl[ \varphi_k\bigl(
X^x(t) \bigr) \bigr] \bigr|
\nonumber
\\[-8pt]
\\[-8pt]
\nonumber
&&\qquad \leq\lim_{ n \to\infty} \sup_{ t \in(0,T] } \sup
_{ x \in\bar{U} } \E\bigl[\bigl | \varphi_k\bigl(
X^{x,n}(t) \bigr) - \varphi_k\bigl( X^x(t)
\bigr) \bigr| \bigr]
\\
&&\qquad \leq\Bigl( \sup_{ x \in\bar{U} } \bigl\| \varphi_k'(x)
\bigr\|_{
L( \R^d, \R)
} \Bigr) \cdot\Bigl( \lim_{ n \to\infty} \sup
_{ t \in(0,T] } \sup_{ x \in\bar{U} } \E\bigl[\bigl |
X^{x,n}(t) - X^x(t) \bigr| \bigr] \Bigr) = 0\nonumber
\end{eqnarray}
for all $ T \in(0,\infty) $
and all $ k \in\N$.
Combining this with
Lemma~\ref{thmlimitsofviscositysolutions}
shows
that for every $ k \in\N$
it holds that
$ u^{ (k) } $
is a viscosity solution
of \eqref{eqCDF2} with initial condition $\varphi_k$.
In addition, note that
%
%e4.101 #&#
%e4.102 #&#
\begin{eqnarray}
\label{eqapproximationvarphi} &&\lim_{ k \to\infty} \sup_{
(t,x) \in
(0,\infty) \times K
} \bigl|
u(t,x) - u^{ (k) }(t,x) \bigr| \nonumber\\
&&\qquad\leq\lim_{ k \to\infty} \sup
_{
(t,x) \in
(0,\infty) \times K
} \E\bigl[ \bigl\vert\varphi\bigl( X^{ x }(t)
\bigr) - \varphi_k\bigl( X^{ x }(t) \bigr) \bigr\vert
\bigr]
\\
&&\qquad \leq\lim_{ k \to\infty} \sup_{
y \in\overline{U}
\cup K
} \bigl\vert
\varphi( y ) - \varphi_k( y ) \bigr\vert= 0\nonumber
\end{eqnarray}
for all compact
sets $ K \subset O $.
Combining
this with
Lemma~\ref{thmlimitsofviscositysolutions}
eventually shows
that $ u $ is indeed a viscosity
solution of \eqref{eqCDF2} as claimed.
\end{pf}

The next result
is a
generalization and a
consequence of
Lemma~\ref{lemKolmogorovviscosity2}
above and constitutes the main result of this section.

%th4.16 #&#
\begin{theorem}[(Existence and uniqueness
of viscosity solutions
of
Kolmogorov equations)]
\label{thmKolmogorovviscosity}
Let
$ d, m \in\N$,
$ \rho\in\R$,
let
$ O \subset\R^d $
be an open set,
let
$
\varphi\colonn O \to\R
$
be a continuous function,
let
$
\mu\colonn O
\rightarrow\R^d
$
and
$
\sigma
\colonn
O \rightarrow
\R^{ d \times m }
$
be locally Lipschitz continuous
functions
and
let
$
V \in\C^2(O, (0,\infty))
$
be
such that
$
\lim_{n \to\infty}
\sup_{ x \in O_n^c }
\frac{
\vert \varphi(x) \vert
}{
1 + V(x)
}
= 0
$,
such that
%
%e4.103 #&#
\begin{equation}
\label{eqLyapunovcondition2} \bigl\langle( \nabla V ) (x), \mu(x)
\bigr\rangle+
\tfrac{ 1 }{ 2 } \operatorname{tr} \bigl( \sigma(x) \bigl[ \sigma(x)
\bigr]^{*} ( \operatorname{Hess} V ) (x) \bigr) \leq\rho\cdot V(x)
\end{equation}
for all $ x \in O $ and such that
$
\lim_{ n \to\infty}
\inf \{
V(x)
\colonn
x \in O_n^c \}
=
\infty
$.
Then
there exists a unique
continuous function
$
u \colonn[0,\infty)
\times O \to\R
$
which fulfills
$
u( 0, x ) = \varphi(x)
$
for all $ x \in O $,
which fulfills
$
\lim_{ n \to\infty}
\sup_{
(t,x) \in
[0,T] \times
O_n^c
}
\frac{
| u(t, x) |
}{
V(x)
}
= 0
$
for all $ T \in(0,\infty) $
and which is a viscosity
solution of
%
%e4.104 #&#
\begin{eqnarray}
\label{eqCDF}&& \frac{ \partial}{ \partial t } u(t,x) - \bigl
\langle(\nabla_x u)
(t,x), \mu(x) \bigr\rangle- \frac{1}{2} \operatorname{tr} \bigl(
\sigma(x)
\bigl[ \sigma(x) \bigr]^{*} % \sigma(x)^{*}
(\operatorname{Hess}_x
u) (t,x) \bigr)
\nonumber
\\[-8pt]
\\[-8pt]
\nonumber
&&\qquad= 0
\end{eqnarray}
for $ (t,x) \in(0,\infty) \times O $.
Moreover,
if
$
( \Omega, \mathcal{F}, \P
)
$
is a probability space
with a normal filtration
$
( \mathcal{F}_t )_{ t \in[0,\infty) }
$
and if
$
W \colonn
[0,\infty)
\times\Omega
\rightarrow\mathbb{R}^m
$
is a
standard
$
( \mathcal{F}_t )_{ t \in[0,\infty) }
$-Brownian motion,
then there exist
up to indistinguishability
unique global solutions
$
X^x \colonn[0,\infty)
\times\Omega
\rightarrow O
$,
$ x \in O $,
to
%
%e4.105 #&#
\begin{equation}
\label{eqSDEviscosityexistence3} X^x(t) = x + \int_0^t
\mu\bigl( X^x(s) \bigr) \,ds + \int_0^t
\sigma\bigl( X^x(s) \bigr) \,dW(s),
\end{equation}
$\P$-a.s. for all
$ t \in[0,\infty) $
and all
$ x \in O $.
In that case,
$u$ has the
probabilistic
representation
$
u(t,x) = \mathbb{E} [
\varphi( X^x(t) )
]
$
for
all $ (t,x) \in[0,\infty) \times O $.
\end{theorem}

\begin{pf}%{Proof
%of Theorem~\ref{thmKolmogorovviscosity}}
W.l.o.g. we assume throughout this proof that
$
( \Omega, \mathcal{F}, \P
)
$
is a probability space
with a normal filtration
$
( \mathcal{F}_t )_{ t \in[0,\infty) }
$
and that
$
W \colonn
[0,\infty)
\times\Omega
\rightarrow\mathbb{R}^m
$
is a
standard
$
( \mathcal{F}_t )_{ t \in[0,\infty) }
$-Brownian motion.
Then, since $ V $ is a Lyapunov function, \eqref
{eqSDEviscosityexistence3} does have global
solutions which furthermore (assuming without loss of generality that
$\rho\ge0$) have the property
that
%
%e4.106 #&#
\begin{equation}
\label{eboundSolLyap} \E\bigl[ V\bigl( X^x( t \wedge\tau) \bigr
) \bigr]
\le e^{\rho t} V(x)
\end{equation}
for any stopping time $ \tau\colonn\Omega\to[0,\infty) $.
As a consequence,
for every $ (t,x) \in[0,\infty) \times O $
it holds that
$
\E[ | \varphi( X^x(t) ) | ]
$
is finite
so that we can \textit{define}
$
u \colonn[0,\infty) \times O \to\R
$
by
$
u(t,x):= \mathbb{E} [ \varphi( X^x(t) ) ]
$
for all $ (t,x) \in[0,\infty) \times O $.
Note that as a consequence of our assumption on $ \varphi$,
for every $ \delta\in(0,\infty) $ there exists a constant $ C_{
\delta} \in(0,\infty) $
such that
%
%e4.107 #&#
\begin{equation}
\label{edecompPhi} \bigl|\varphi(x)\bigr| \le C_\delta+ \delta V(x)
\end{equation}
holds for all $x \in O$. The bound \eqref{eboundSolLyap}
immediately implies a similar bound on $u(t,\cdot)$, so that $u$ has
the required behaviour at infinity.
It therefore remains to show that $u$ is indeed a viscosity solution of
\eqref{eqCDF}, as uniqueness of such a solution
follows from Corollary~\ref{coruniqueness2}.
The proof for this goes again by approximation. Let $\mu^{(n)}$ and
$\sigma^{(n)}$ for $ n \in\N$ be any sequence of
Lipschitz continuous functions
such that for all $ x \in O $ it holds that
%
%e4.108 #&#
\begin{equation}
V(x) \le n \quad\Rightarrow\quad\mu^{(n)}(x) = \mu(x), \qquad \sigma^{(n)}(x)
= \sigma(x)
\end{equation}
and
% such that
%
%e4.109 #&#
\begin{equation}
V(x) \ge n+1 \quad\Rightarrow\quad\mu^{(n)}(x) = 0, \qquad \sigma^{(n)}(x) =
0.
\end{equation}
Denoting by $ X^{x,n} $, $ x \in O $, $ n \in\N$, the solutions to
the corresponding SDEs,
we set
$ u_n(t,x) = \E[ \varphi( X^{x,n}( t ) ) ] $
for all $ (t,x) \in[0,\infty) \times O $, $ n \in\N$.
It then follows from Lemma~\ref{lemKolmogorovviscosity2} that
$ u_n|_{ ( 0, \infty) \times O_n } $ is a viscosity solution to the
equation analogous to~\eqref{eqCDF}.
As a consequence of Lemma~\ref{thmlimitsofviscositysolutions},
it remains to show that
$ u_n \to u $, uniformly on compact subsets
of $ ( 0, \infty) \times O $.
For this, we introduce the stopping times
$
\tau_n^x:=
\inf (
\{ t \in(0,\infty) \colonn V(X^x(t)) \ge n\} \cup\{ \infty\}
)
$,
$ x \in O $,
$ n \in\N$.
As a consequence of \eqref{edecompPhi}, the fact that $ X^{x,n} $ and
$ X^x $ coincide
until time $ \tau_n^x $,
and the fact that
$
V ( X^{ x, n}( t ) ) \le n + 1
$, $\P$-a.s. provided that
$
V(x) \le n + 1
$,
we have
for all
$ n \in\N$
and all
$ (t,x) \in[0,\infty) \times O $
with $ V( x ) \leq n + 1 $
that
%
%e4.110 #&#
%e4.111 #&#
\begin{eqnarray}
\label{eboundDiffu}&& \bigl|u(t,x)- u_n(t,x)\bigr|\nonumber \\
&&\qquad \leq\E\bigl[
\mathbh{1}_{
\{
\tau^x_n \leq t
\}
} \bigl| \varphi\bigl( X^x(t) \bigr) \bigr| \bigr] + \E
\bigl[ \mathbh{1}_{
\{
\tau^x_n \leq t
\}
} \bigl| \varphi\bigl( X^{ x, n }(t) \bigr) \bigr|
\bigr]
\\
&&\qquad \le 2 C_\delta\P\bigl[ \tau_n^x \le t
\bigr] + \delta e^{ \rho t } V(x) + \delta( n + 1 ) \P\bigl[
\tau_n^x \le t \bigr].\nonumber
\end{eqnarray}
Using
% again
\eqref{eboundSolLyap}, we obtain from Chebychev's inequality that
for all $ (t,x) \in[0,\infty) \times O $ it holds that
%
%e4.112 #&#
\begin{equation}\qquad
\P\bigl[ \tau_n^x \le t \bigr] = \P\bigl[ V \bigl(
X^x\bigl( t \wedge\tau^x_n \bigr) \bigr) \ge
n \bigr] \leq{
\E[
V ( X^x( t \wedge\tau^x_n ) )
]
\over
n
} \le{e^{\rho t} V(x) \over n}.
\end{equation}
Inserting this into \eqref{eboundDiffu}, the required locally uniform
convergence follows at once.
\end{pf}

In the literature, there are many results proving an assertion
similar to
Theorem~\ref{thmKolmogorovviscosity}
and
Corollary~\ref{coruniqueness2},
respectively,
under
various assumptions on
the functions $ \mu$ and $ \sigma$.
Theorem 4.3 in
Pardoux and Peng~\cite{PardouxPeng1992}
%with $f\equiv0$
implies that the transition semigroup associated with the SDE~\eqref
{eqSDEviscosityexistence3}
is a viscosity solution
of~\eqref{eqCDF} if $\mu$ and $\sigma$ are globally Lipschitz continuous;
see also Peng~\cite{Peng1993}.
Theorem C.2.4 in Peng~\cite{Peng2010}
can be applied if $ \mu$ is
locally H\"{o}lder continuous
and if $ \sigma$ is constant
and then proves uniqueness
of an at most polynomially
growing
viscosity solution of~\eqref{eqCDF}.
Uniqueness of the viscosity solution
of~\eqref{eqCDF}
with given initial function
follows from
Theorem~8.2
in the User's guide
Crandall, Ishii and Lions~\cite{CrandallIshiiLions1992}
if $\mu$ is globally one-sided Lipschitz continuous,
that is,
if
there exists a constant $c\in\R$ such that
$
\langle x-y,\mu(x)-\mu(y)\rangle
\leq c \| x - y \|^2
$
for all $ x, y \in\R^d $,
and if
$\sigma$
is globally Lipschitz continuous.
Moreover,
Theorem~5.13
in Krylov~\cite{Krylov1999}
%with
%$
% f \equiv0 \equiv c
%$
implies that the
transition semigroup
solves the Kolmogorov
equation~\eqref{eqCDF}
in the sense of distributions
if $ \mu$ and $ \sigma$
are globally Lipschitz
continuous.
In addition,
Theorems~7.1.3 and~7.1.4 in Evans~\cite{Evans2010}
show that there exists a unique weak solution
of the PDE~\eqref{eqCDF}
if the coefficients $\mu$ and $\sigma$ are bounded and if
the PDE~\eqref{eqCDF} is
uniformly parabolic.

In many situations,
the open set $ O \subset\R^d $
and the Lyapunov-type
function $ V \colonn O \to\R$
in
Theorem~\ref{thmKolmogorovviscosity}
satisfy
$ O = \R^d $
and
$
V(x) = ( 1 + \| x \|^2 )^p
$
for all $ x \in\R^d $
where $ p \in[1,\infty) $ is
an arbitrary real number.
This is subject of the following Corollary~\ref{corKolmogorovviscosity}.
It is a direct
consequence of
Theorem~\ref{thmKolmogorovviscosity}
and its proof is therefore
omitted.

%co4.17 #&#
\begin{corollary}[(Existence and uniqueness
of at most polynomially growing
viscosity solutions of Kolmogorov
equations)]
\label{corKolmogorovviscosity}
Let
$ d, m \in\mathbb{N} $,
let
$
\varphi\colonn\R^d \to\R
$
be a continuous and
at most
polynomially growing function,\break 
let
$
\mu\colonn\R^d
\rightarrow\R^d
$
and
$
\sigma
\colonn
\R^d \rightarrow
\R^{ d \times m }
$
be locally Lipschitz continuous
functions
with
$
\sup_{ x \in\R^d }
\frac{
\anglel x, \mu(x) \angler
}{
(
1 + \| x \|^2
)
}
< \infty
$
and
$
\sup_{ x \in\R^d }
\frac{
\llVert \sigma(x) \rrVert
}{
(
1 + \| x \|
)
}
< \infty
$.
Then
there exists\break a unique
continuous function
$
u \colonn[0,\infty)
\times\R^d \to\R
$
which fulfills\break
$
\limsup_{ p \to\infty}
\sup_{ (t,x) \in[0,T] \times\R^d }
\frac{ | u(t,x) | }{
1 + \| x \|^p
}
< \infty
$
for all $ T \in( 0, \infty) $,
% which is at most polynomially growing as $ | x | \to\infty$ @@@,
% is such that
which fulfills\break 
$
u( 0, x ) = \varphi(x)
$
for all $ x \in\R^d $,
and which is a viscosity
solution of
%
%e4.113 #&#
\begin{eqnarray}\qquad
\label{eqCDF4} &&\frac{ \partial}{ \partial t } u(t,x) - \bigl
\langle(\nabla_x u)
(t,x), \mu(x) \bigr\rangle- \frac{1}{2} \operatorname{tr} \bigl(
\sigma(x)
\bigl[ \sigma(x) \bigr]^{*} % \sigma(x)^{*}
(\operatorname{Hess}_x
u) (t,x) \bigr)
\nonumber
\\[-8pt]
\\[-8pt]
\nonumber
&&\qquad = 0
\end{eqnarray}
for $ (t,x) \in(0,\infty) \times\R^d $.
Moreover,
if $ ( \Omega, \mathcal{F}, \P) $
is a probability space with a normal filtration
$ ( \mathcal{F}_t )_{ t \in[0,\infty) } $
and if $ W \colonn[0,\infty) \times\Omega\to\R^m $
is a standard $ ( \mathcal{F}_t )_{ t \in[0,\infty) } $-Brownian
motion, then
$ u $ has the
probabilistic
representation\break 
$
u(t,x) = \mathbb{E} [
\varphi( X^x(t) )
]
$
for all $ (t,x) \in[0,\infty) \times\R^d $,
where the stochastic processes $ X^x \colonn [0,\infty) \times\Omega
\to\R^d $, $ x \in\R^d $, are as before.
\end{corollary}

Note that all examples
in this article fulfill the
assumptions of
Corollary~\ref{corKolmogorovviscosity}.
In particular, observe that
$ \mu$ and $ \sigma$ from
the SDE~\eqref{eqSDEbsp1}
in Section~\ref{secex1},
%(see the
%PDE~\eqref{eqex1Kolmogorov}
%for details),
%
$ \mu$ and $ \sigma$ from
the SDE~\eqref{eqSDEbsp1b}
in Section~\ref{secex1},
$ \mu$ and $ \sigma$ from
the SDE~\eqref{eqSDEbsp1c}
in Section~\ref{secex1},
$ \mu$ and $ \sigma$ from
the SDE~\eqref{eqex2bSDE}
in Section~\ref{secex2},
$ \mu$ and $ \sigma$
from the
SDE~\eqref{eqSDESec3B}
%in Theorem~\ref{thmirr2}
in
Section~\ref{secex2}
as well as
$ \mu$ and $ \sigma$ from
the SDE~\eqref{eqex3SDE}
in Section~\ref{secex4}
%(see the
%PDE~\eqref{eqex3Kolmogorov}
%for details)
all fulfill the assumptions
of Corollary~\ref{corKolmogorovviscosity}.

%%%%%%%%%%%%%%%%%%%%%%%%%%%%
%s4.5 #&#
\subsection{Distributional solutions of
Kolmogorov equations}
\label{ssecSolutionsofKolmogorovequationsinthedistributionalsense}

In this section, we formulate a slight extension to Theorem~5.13 in
Krylov~\cite{Krylov1999}, which
states that the semigroup
associated to an SDE with
smooth coefficients solves
the corresponding
Kolmogorov equation in the distributional
sense, even if the coefficients are badly behaved near the boundary of
the domain
of definition~$O$.

%pr4.18 #&#
\begin{prop}
\label{propKolmogorov}
Let $ d, m \in\mathbb{N} $,
let $ O \subset\mathbb{R}^d $
be an open set,
let
$
\mu=
(
\mu_1, \ldots, \mu_d
)
\in\C^{ \infty}( O, \R^d )
$,
$
\sigma
=
(
\sigma_{ i, j }
)_{
i \in\{ 1, \ldots, d \},
j \in\{ 1, \ldots, m \}
}
\in
\C^{ \infty}( O, \R^{ d \times m } )
$,
let
$
\varphi\in
\C_b( O, \R)
$,
let
$ ( \Omega, \mathcal{F}, \P) $
be a probability space with a normal filitration
$
( \mathcal{F}_t )_{ t \in[0,\infty) }
$,
let
$ W \colonn[0,\infty) \times\Omega\to\R^m $
be a standard $ ( \mathcal{F}_t )_{ t \in[0,\infty) } $-Brownian motion
and let
$
X^x \colonn[0,\infty)
\times\Omega
\rightarrow O
$,
$ x \in O $,
be solutions to
%
%e4.114 #&#
\begin{equation}
X^x(t) = x + \int_0^t \mu\bigl(
X^x(s) \bigr) \,ds + \int_0^t
\sigma\bigl( X^x(s) \bigr) \,dW(s),
\end{equation}
$\P$-a.s. for all $ (t,x) \in[0,\infty) \times\Omega$.
% which is assumed to exist and take values on $O$ for all
% $t >0$ and all
% $x \in O $.
Then the function
$
u \colonn(0,\infty) \times O
\rightarrow\mathbb{R}
$
given by
$
u(t,x) = \mathbb{E} [
\varphi( X^x(t) )
]
$
for all $ (t,x) \in[0,\infty) \times O $
solves the Kolmogorov equation
%
%e4.115 #&#
\begin{equation}
\frac{
\partial u
}{
\partial t
} = \sum_{ i = 1 }^{ d }
\mu_i \frac{ \partial u}{
\partial x_{ i }
} + \frac{ 1 }{ 2 } \sum
_{ l = 1 }^{ m } \sum_{ i, j = 1 }^{ d }
\sigma_{ i, l } \sigma_{ j, l } \frac{ \partial^2 u}{
\partial x_i\, \partial x_j
}
\end{equation}
in the distributional sense.
\end{prop}

\begin{pf}%{Proof
%of Proposition~\ref{propKolmogorov}}
Let $O_n$ be as above, consider for every $ n \in\N$ smooth and
globally Lipschitz continuous
functions $ \mu^{(n)} $ and $ \sigma^{(n)} $ which agree with $\mu$
and $\sigma$ on $O_n$,
and denote by $ X^{x,n} $, $ x \in O $, solutions of the corresponding SDE.
Fix some final time $ T \in(0,\infty) $, denote by $P_x$ the law of $
X^x $ on $ C( [0,T], O ) $
and for every $ n \in\N$ by $ P_x^n $ the law of $ X^{x,n} $ on $ C(
[0,T], O ) $.
It then follows from the smoothness of the coefficients $ \mu$ and $
\sigma$
that $ O \ni x \mapsto P_x $ is weakly continuous; see Theorem~1.7 in
Krylov~\cite{Krylov1999}.
In particular, this implies that $u$ is continuous and that, for every compact
$K \subset O$, the set $ \{ P_x \colonn x \in K \} $ is tight.
Let now
$
u_n(x,t) = \E[ \varphi_n( X^{ x, n}( t ) ) ]
$
for all $ (t,x) \in(0,\infty) \times O $, $ n \in\N$,
where $ \varphi_n \colonn O \to\R$, $ n \in\N$, are smooth
approximations of
$ \varphi$ such that
$
\sup_{ x \in O_n }
| \varphi_n(x) - \varphi(x) |
\le
1 / n
$
and
$
\operatorname{supp}( \varphi_n ) \subset O_{ n + 1 }
$
for all $ n \in\N$
and such that
$
\sup_{ n \in\N}
\sup_{ x \in O }
| \varphi_n( x ) |
< \infty
$.
% and such that
% the functions $\phi_n$ are
% uniformly bounded, and
% $\phi_n$ is compactly supported on $O_{n+1}$.
%
Note now that
$
P_x|_{ \mathcal{B}( C( [0,T], O_n ) ) }
=
P_x^n|_{
\mathcal{B}( C( [0,T], O_n ) )
}
$
% $ P_x $ and $ P_x^n $ agree on the set $ C( [0,T], O_n ) $ of
%continuous paths with values in $ O_n $
and
that, locally uniformly in $ x $, the $ P_x $-measure of the set
$
C( [0,T], O_n )
$
converges to $ 1 $
as $ n \to\infty$. In particular, there exists a real number $ C \in
[0,\infty) $ such that
for all $ (t,x) \in(0,T] \times O $ it holds that
%
%e4.116 #&#
\begin{equation}
\bigl|u_n(x,t) - u(x,t)\bigr| \le{1\over n} + C \bigl[ 1 -
P_x \bigl( C\bigl( [0,T], O_n \bigr) \bigr) \bigr].
\end{equation}
As a consequence, one has $ u_n \to u $, locally uniformly in $ x $ and
$ t $. The claim now follows at once
from the fact that, by Theorem~5.13 in
Krylov~\cite{Krylov1999}, each of the $u_n$ solves the Kolmogorov
equation with $\mu^{(n)}$ and $\sigma^{(n)}$.
\end{pf}

%s5 #&#
\section{A counterexample
to the rate
of convergence
of the Euler--Maruyama method}
\label{secex4}

In this section, we use
the results
of Section~\ref{secex2}
to establish
the existence of
an SDE with smooth and
globally bounded coefficients
for which the Euler--Maruyama method
convergences
\textit{without
% slower than
any
arbitrarily small polynomial
%positive
rate of convergence}, thereby proving Theorem~\ref
{thmnonrateintroduction} of the \hyperref[sec1]{Introduction}.
Denote by $\hat C$ the constant
%
%e5.1 #&#
\begin{equation}
\hat C = \int_{ 0 }^1 e^{ -{1 }/{ ( 1 - u^2 ) } } \,du,
\end{equation}
and set
%
%e5.2 #&#
\begin{eqnarray}
\label{eqex3defmu} \mu(x) &= &\pmatrix{\displaystyle \mathbh{1}_{ (1,\infty) }(
x_4 ) \cdot\exp\biggl( -\frac{ 1 }{ x_4^2 - 1 } \biggr) \cdot\cos
\bigl( ( x_3 - \hat C ) \cdot\exp\bigl( x_2^3
\bigr) \bigr)
\vspace*{2pt}\cr
0
\vspace*{2pt}\cr
\displaystyle\mathbh{1}_{ (-1,1) }( x_4 ) \cdot\exp\biggl( -
\frac{1 }{ 1 - x_4^2 } \biggr)
\vspace*{2pt}\cr
1
}
, %\end{equation}
%and by
\nonumber
\\[-8pt]
\\[-8pt]
\nonumber
 B& =& \pmatrix{ 0 & 0 & 0 & 0
\vspace*{2pt}\cr
0 & 1 & 0 & 0
\vspace*{2pt}\cr
0 & 0 & 0 & 0
\vspace*{2pt}\cr
0 & 0 & 0 & 0}
\end{eqnarray}
for all
$
x = (x_1, x_2, x_3, x_4) \in
\mathbb{R}^4
$.
The function
$ \R\ni x\mapsto
\mathbh{1}_{ (-1,1) }( x )
\cdot
\exp(
- 1 / ( 1 - x^2 )
)
\in[0,1]
$
that appears in
$ \mu$ has been used
as a mollifier function
in Lemma~1.2.3
in
H\"{o}rmander~\cite{Hoermander1990}.
Note that
$
\mu\colonn\R^4 \to\R^4
$
is infinitely often
differentiable and globally bounded.
Moreover, let
$
( \Omega, \mathcal{F}, \P
)
$
be any probability space
supporting
a four-dimensional
standard Brownian
motion
$
W \colonn[0,\infty) \times\Omega
\to\R^4
$
with continuous
sample paths.
Then there exists
a unique stochastic process
$
X
\colonn[0,\infty) \times
\Omega\to\R^4
$
with
continuous sample paths
which fulfills
$X(t) =
\int_0^t \mu( X(s) ) \,ds
+
B W(t)$
for all $ t \in[0,\infty) $.
%see, e.g.,
%Corollary 2.6 in Gy\"ongy \citationand\ Krylov~
The stochastic
process
$
X = (X_1, X_2, X_3, X_4)
\colonn[0,\infty) \times
\Omega\to\R^4
$
is thus a solution
process of the SDE
%
%e5.3 #&#
\begin{eqnarray}
\label{eqex3SDE} d X_1(t) & =& \mathbh{1}_{
(1, \infty)
} \bigl(
X_4(t) \bigr) \cdot\exp\biggl( -\frac{1 }{
X_4(t)^2 - 1
} \biggr) \nonumber\\
&&{}\times
\cos\bigl( \bigl( X_3(t) - \hat C \bigr) \cdot\exp\bigl(
X_2(t) ^3 \bigr) \bigr) \,dt,
\nonumber\\
d X_2(t) & =& dW_2(t),
\\
d X_3(t) & =& \mathbh{1}_{
(-1,1)
} \bigl( X_4(t)
\bigr) \cdot\exp\biggl( -\frac{ 1 }{
1 - X_4(t)^2
} \biggr) \,dt,
\nonumber\\
d X_4(t) & =& 1 \,dt\nonumber
\end{eqnarray}
for $ t \in[0,\infty) $
satisfying
$ X(0) = 0 $.
%Moreover, observe that
% X_1(t) =
% \int_1^t
% \exp\left(
% \frac{ - 1 }{
% ( s^2 - 1 )
% }
% \right)
% ds
%for all $ t \in[1,\infty) $.
In the next step, we define the
Euler--Maruyama approximations
for the SDE~\eqref{eqex3SDE}
using
the following notation.
Let
$
\lfloor\cdot\rfloor_h
\colonn
[0,\infty) \to[0,\infty)
$,
$ h \in(0,\infty) $,
be a family of mappings
defined by
%
%e5.4 #&#
\begin{equation}
\label{eqfloorh} \lfloor t \rfloor_h:= \max\bigl\{ s \in\{ 0, h,
2h, \ldots\} \colonn s \leq t \bigr\}
\end{equation}
for all $ t \in[0,\infty) $
and all $ h \in(0,\infty) $.
Then let
$
Y^h =
( Y^h_1, Y^h_2, Y^h_3, Y^h_4 )
\colonn[0,\infty)
\times\Omega\to\R^4
$,
$ h \in(0,\infty) $,
be Euler--Maruyama
approximation processes
defined recursively
by
%
%e5.5 #&#
\begin{eqnarray}
\label{eqex3Euler} Y^h(0)&:=& 0\quad \mbox{and}
\nonumber
\\[-8pt]
\\[-8pt]
\nonumber
  Y^h( t )&:=&
Y^h\bigl( \lfloor t \rfloor_h \bigr) + \mu\bigl(
Y^h\bigl( \lfloor t \rfloor_h \bigr) \bigr) \cdot
\bigl( t - \lfloor t \rfloor_h \bigr) + B \bigl( W( t ) - W\bigl(
\lfloor t \rfloor_h \bigr) \bigr)
\end{eqnarray}
for all
$ t \in
( n h,
(n + 1) h
]
$,
$ n \in\{ 0, 1, \ldots\} $
and all $ h \in(0,\infty) $.
Observe that this definition
ensures that
%
%e5.6 #&#
\begin{eqnarray}
\label{eqexample3Y1rep}  Y^h_1(t) &=& \int_1^t
\mathbh{1}_{
(1, \infty)
} \bigl( \lfloor s \rfloor_h \bigr)
% \cdot
e^{
% \left(
-{ 1 }/{
(\lfloor s \rfloor_h^2 - 1)
}
% \right)
} %\\ &
% \cdot
\nonumber
\\[-8pt]
\\[-8pt]
\nonumber
&&\hspace*{13pt}{}\times\cos\biggl( \biggl(
\smallint_{ 0 }^{ \infty} \mathbh{1}_{
[0,1)
} \bigl( \lfloor
u \rfloor_h \bigr) % \cdot
e^{
-{ 1 }/{
(1 -
\lfloor u \rfloor_h^2)
}
} \,du - \hat C \biggr)
% \cdot
e^{
W_2( \lfloor s \rfloor_h )^3
} \biggr) \,ds
\end{eqnarray}
for all $ t \in[1,\infty) $
and all $ h \in(0,\infty) $.
The following
Theorem~\ref{thmnonrate}
proves that the
Euler--Maruyama method~\eqref{eqex3Euler}
for the SDE~\eqref{eqex3SDE}
convergences
\textit{%
% slower than
without any
arbitrarily small polynomial
%positive
rate of convergence}.
Theorem~\ref{thmnonrate}
together with an elementary
transformation argument
[dealing with general
$ x_0 \in\R^4 $
and general $ T \in( 0, \infty) $]
then implies
Theorem~\ref{thmnonrateintroduction}.

%%%%%%%%%%%%%%%%%%%%%%%%%%%%%%%%%%%%%%%
%
%th5.1 #&#
\begin{theorem}[(A counterexample
to the rate of convergence
of the Euler--Maruyama method)]
\label{thmnonrate}
Let
$
X
= ( X_1, X_2, X_3, X_4 )
\colonn[0,\infty) \times\Omega
\to\R^4
$
be
%the up to indistinguishability
%unique
a solution process
of the SDE~\eqref{eqex3SDE}
with continuous sample paths and
with $ X(0) = 0 $.
Then
%
%e5.7 #&#
\begin{equation}
\E\bigl[ X_1(t) \bigr] - \E\bigl[ Y^h_1(t)
\bigr] \geq\exp\bigl( - 14 \bigl| \ln( h ) \bigr|^{ 2 / 3 } \bigr)
\end{equation}
for all
$
h \in( 0,
\frac{ 1 }{ 22 }
]
$
and all
$ t \in[2, \infty) $
and, therefore, we obtain
%
%e5.8 #&#
\begin{eqnarray}
\lim_{ h \searrow0 } \biggl( \frac{
\E[
\|
X(t)
-
Y^h(t)
\|
]
}{
h^{ \alpha}
} \biggr)& =& \lim
_{ h \searrow0 } \biggl( \frac{
\llVert
\E[
X(t)
]
-
\E[
Y^h(t)
]
\rrVert
}{
h^{ \alpha}
} \biggr)
\nonumber
\\[-8pt]
\\[-8pt]
\nonumber
&=& \cases{ 0, &\quad $
\alpha= 0,$ \vspace*{2pt}
\cr
\infty,& \quad $\alpha> 0,$ }
\end{eqnarray}
for all $ \alpha\in[0,\infty) $
and all $ t \in[2,\infty) $.
In particular, for every
$ t \in[2,\infty) $
and
every $ \alpha, C, h_0 \in(0,\infty) $
there exists a real number
$ h \in( 0, h_0 ) $ such that
$
\llVert
\E[
X(t)
]
-
\E [
Y^h(t)
]
\rrVert
>
C
\cdot h^{ \alpha}
$.
\end{theorem}

The proof of
Theorem~\ref{thmnonrate} is deferred
to the end of this section.
To the best of our
knowledge,
the SDE~\eqref{eqex3SDE}
is the first SDE with smooth coefficients
in the literature
for which it has been
established that the Euler--Maruyama
scheme converges in the strong
and numerical weak sense
\textit{without
% slower than
any arbitrarily
small rate of convergence}.
Using the results of
Section~\ref{secex2},
one can show
that the SDE~\eqref{eqex3SDE}
is not locally H\"{o}lder
continuous with respect to
the initial value.
This is summarized in the next corollary
of Lemma~\ref{lemex2lowerbound}
in Section~\ref{secex2}.
%Its proof
%is a straightforward
%consequence of
%Lemma~\ref{lemex2lowerbound}
%in Section~\ref{secex2}
%and is therefore omitted.

%co5.2 #&#
\begin{corollary}
\label{corex4initial}
Let
$
X^x \colonn[0,\infty) \times\Omega
\to\R^4
$,
$ x \in\R^4 $,
be solution processes
of the SDE~\eqref{eqex3SDE}
with continuous
sample paths
and with
$
X^x(0)
=
x
$
for all
$ x \in\R^4 $.
Then
for every $ t \in(0,\infty) $
the function
$
\R^4 \ni x \mapsto
\E [
X^x(t)
]
\in\R^4
$
is not locally
H\"{o}lder
continuous.
\end{corollary}

\begin{pf}%{Proof
%of Corollary~\ref{corex4initial}}
Note that
%
%e5.9 #&#
%e5.10 #&#
%e5.11 #&#
%e5.12 #&#
\begin{eqnarray}
&&\bigl\llVert\E\bigl[ X^{
(0,0, \hat C,
2 )
}(t) \bigr] - \E\bigl[
X^{
(0,0,
h + \hat C,
2 )
}(t) \bigr] \bigr\rrVert
\nonumber\\
&&\qquad \geq\bigl\vert\E\bigl[
X_1^{
(0,0,
\hat C,
2 )
}(t) - X_1^{
(0,0,
h + \hat C,
2 )
}(t) \bigr]
\bigr\vert
\nonumber
\\[-8pt]
\\[-8pt]
\nonumber
&&\qquad = \biggl\vert\int_0^t \exp\biggl(
\frac{ - 1 }{
(
( 2 + s )^2 - 1
)
} \biggr) \E\bigl[ 1 - \cos\bigl( h \cdot\exp\bigl( \bigl[
W_2(s) \bigr]^3 \bigr) \bigr) \bigr] \,ds \biggr\vert
\\
&&\qquad \geq\exp\biggl( -\frac{1 }{
3
} \biggr) \int_0^t
\bigl( 1 - \E\bigl[ \cos\bigl( h \cdot\exp\bigl( \bigl[ W_2(s)
\bigr]^3 \bigr) \bigr) \bigr] \bigr) \,ds\nonumber
\end{eqnarray}
for all $ h, t \in( 0, \infty) $.
Combining this with
Lemma~\ref{lemex2lowerbound}
in Section~\ref{secex2}
completes
the proof of
Corollary~\ref{corex4initial}.
\end{pf}

%Note also that the Kolmogorov
%PDE associated with~\eqref{eqex3SDE}
%reads as
% \frac{ \partial}{ \partial t }
% u(t,x)
%=
% \frac{ 1 }{ 2 }
% \frac{ \partial^2 }{ \partial x_2^2 }
% u(t,x)
% +
% \mu_1(x)
% \cdot
% \frac{ \partial}{ \partial x_1 }
% u(t,x)
% +
% \mu_3(x)
% \cdot
% \frac{ \partial}{ \partial x_3 }
% u(t,x)
% +
% \frac{ \partial}{ \partial x_4 }
% u(t,x)
%for $ t \in(0,\infty) $
%and
%$ x = (x_1, x_2, x_3, x_4) \in\R^4 $.
%Corollary~\ref{corKolmogorovviscosity}
%above proves
%existence and uniqueness
%of viscosity solutions
%of the
%PDE~\eqref{eqex3Kolmogorov}
%with a given initial condition
%(see
%Corollary~\ref{corKolmogorovviscosity}
%for details).
%
In the following, the size
of the quantity
$
\|
\E[ X(T) ]
-
\E[ Y^h(T) ]
\|
\in[0,\infty)
$
is analyzed
for sufficiently small
$ h \in(0,\infty) $
and thereby
Theorem~\ref{thmnonrate}
is established.
To do so, we first establish
a few auxiliary results.
We begin with an elementary
estimate for the numerical
integration of concave
functions.

%le5.3 #&#
\begin{lemma}[(Numerical
integration of concave functions)]
\label{lemnumintegrate}
Let
$
\lfloor\cdot\rfloor_h
\colonn\break 
[0, \infty) \to[0,\infty)
$,
$ h \in(0,\infty) $,
be given by \eqref{eqfloorh},
let $ b \in(0,\infty) $
be a real number and
let
$
\psi\colonn[0,b]
\rightarrow\R
$
be a continuously differentiable
function with a nonincreasing
derivative.
Then
%
%e5.13 #&#
\begin{eqnarray}\qquad
&&\int_0^b \bigl( \psi( s ) - \psi\bigl(
\lfloor s \rfloor_h \bigr) \bigr) \,ds
\nonumber
\\[-8pt]
\\[-8pt]
\nonumber
&&\qquad\leq\frac{ 1 }{ 2 } \bigl[
\psi'( 0 ) \cdot h^2 + \bigl( \psi\bigl( \lfloor b
\rfloor_h - h \bigr) - \psi( 0 ) \bigr) \cdot h + \psi'
\bigl( \lfloor b \rfloor_h \bigr) \cdot\bigl( b - \lfloor b
\rfloor_h \bigr)^2 \bigr]
\end{eqnarray}
for all $ h \in(0,b] $.
\end{lemma}

\begin{pf}%{Proof
%of Lemma~\ref{lemnumintegrate}}
The fundamental theorem
of calculus and monotonicity
of $ \psi' $ imply
%
%e5.14 #&#
%e5.15 #&#
%e5.16 #&#
%e5.17 #&#
%e5.18 #&#
\begin{eqnarray}
 &&\int_0^b \bigl( \psi( s ) - \psi\bigl(
\lfloor s \rfloor_h \bigr) \bigr) \,ds\nonumber\\
&&\qquad = \int_0^b
\int_{ \lfloor s \rfloor_h }^s \psi'( u ) \,du \,ds \leq
\int_0^b \int_{ \lfloor s \rfloor_h }^s
\psi'\bigl( \lfloor s \rfloor_h \bigr) \,du \,ds
\\
&&\qquad  =\int_0^{
h
} \int_{ \lfloor s \rfloor_h }^s
\psi'\bigl( \lfloor s \rfloor_h \bigr) \,du \,ds + \int
_h^{
\lfloor b \rfloor_h
} \int_{ \lfloor s \rfloor_h }^s
\psi'\bigl( \lfloor s \rfloor_h \bigr) \,du \,ds\nonumber\\
&&\qquad\quad{} + \int
_{
\lfloor b \rfloor_h
}^b \int_{ \lfloor s \rfloor_h }^s
\psi'\bigl( \lfloor s \rfloor_h \bigr) \,du \,ds\nonumber
\nonumber\\
&&\qquad = \psi'( 0 ) \cdot\frac{ h^2 }{ 2 } + \frac{ h^2 }{ 2 }
\biggl( \sum_{
n \in\N,
n h < \lfloor b \rfloor_h
%}
} \psi'( n h ) \biggr)
+ \psi'\bigl( \lfloor b \rfloor_h \bigr) \cdot
\frac{
( b - \lfloor b \rfloor_h )^2
}{
2
} %\\ & =
% \psi'( 0 )
% \cdot
% \frac{ h^2 }{ 2 }
% +
% \frac{ h }{ 2 }
% \left(
% \sum_{
% \substack{
% n \in\N
% \\
% n h < \lfloor b \rfloor_h
% }
% }
% \int_{ (n - 1) h }^{ n h }
% \psi'( n h )
% ds
% \right)
% +
% \psi'( \lfloor b \rfloor_h )
% \cdot
% \frac{
% \left( b - \lfloor b \rfloor_h \right)^2
% }{
% 2
% }
\nonumber\\
&&\qquad \leq\psi'( 0 ) \cdot\frac{ h^2 }{ 2 } + \frac{ h }{ 2 }
\biggl( \sum_{
n \in\N,
n h < \lfloor b \rfloor_h
%}
} \int_{ (n - 1) h }^{ n h }
\psi'( s ) \,ds \biggr) \nonumber\\
&&\qquad\quad{}+ \psi'\bigl( \lfloor b
\rfloor_h \bigr) \cdot\frac{
( b - \lfloor b \rfloor_h )^2
}{
2
}
\nonumber\\
&&\qquad = \psi'( 0 ) \cdot\frac{ h^2 }{ 2 } + \bigl( \psi\bigl(
\lfloor b \rfloor_h - h \bigr) - \psi( 0 ) \bigr) \cdot
\frac{ h }{ 2 }\nonumber\\
&&\qquad\quad{} + \psi'\bigl( \lfloor b \rfloor_h
\bigr) \cdot\frac{
( b - \lfloor b \rfloor_h )^2
}{
2
}\nonumber
\end{eqnarray}
for all $ h \in(0,b] $.
This finishes the proof
of
Lemma~\ref{lemnumintegrate}.
\end{pf}

Using Lemma~\ref{lemnumintegrate},
we establish in the next lemma
a simple lower bound
for the numerical integration
of
the function
$
\mathbh{1}_{ (-1,1) }( x )
\cdot
\exp(
- 1 / ( 1 - x^2 )
)
$,
$ x \in\R$,
in the third component of
$ \mu\colonn\R^4 \to\R^4 $.

%le5.4 #&#
\begin{lemma}[{[Numerical integration
of the function
$
\mathbh{1}_{ (-1,1) }( x )
\cdot
\exp(
- 1 / ( 1 - x^2 )
)
$,
$ x \in\R$]}]
\label{lemmollifier}
Let
$
\lfloor\cdot\rfloor_h
\colonn
[0,\infty) \to[0,\infty)
$,
$ h \in(0,\infty) $,
be given by \eqref{eqfloorh}.
Then
%
%e5.19 #&#
\begin{equation}
\frac{ h }{ 20 } \leq\int_0^{ \infty}
\mathbh{1}_{
[ 0, 1 )
} \bigl( \lfloor s \rfloor_h \bigr) \cdot
\exp\biggl( -\frac{1 }{
1 - \lfloor s \rfloor_h ^2
} \biggr) \,ds - \hat C \leq2 h
\end{equation}
for all $ h \in(0,\frac{1}{8}] $.
\end{lemma}

\begin{pf}%{Proof
%of Lemma~\ref{lemmollifier}}
First of all, observe that
%
%e5.20 #&#
\begin{eqnarray}
\label{eqex3mollifiersecondderivative} \frac{
d
}{
d x
} \bigl( e^{ - 1 / ( 1 - x^2 ) } \bigr) &=&
\frac{
- 2 x \cdot
e^{ - 1 / ( 1 - x^2 ) }
}{
( 1 - x^2 )^2
} \quad\mbox{and}
\nonumber
\\[-8pt]
\\[-8pt]
\nonumber
 \frac{
d^2
}{
d x^2
} \bigl( e^{ - 1 / ( 1 - x^2 ) } \bigr)
&=& \frac{
6 \cdot e^{ - 1 / ( 1 - x^2 ) }
}{
( 1 - x^2 )^4
} \biggl( x^4 - \frac{ 1 }{ 3 } \biggr)
\end{eqnarray}
for all $ x \in(-1,1) $.
We hence obtain
that the function
$
[0,3^{-1/4}]
\ni
s \mapsto\break 
e^{ - 1 / ( 1 - s^2 ) }
\in\R
$
has a nonincreasing derivative.
Applying Lemma~\ref{lemnumintegrate}
and using that the function
$
[0, \infty)
\ni
s \mapsto
\mathbh{1}_{ [0,1) }(s) \cdot
e^{ - 1 / ( 1 - s^2 ) }
\in\R
$
is nonincreasing therefore
results in
%
%e5.21 #&#
%e5.22 #&#
%e5.23 #&#
%e5.24 #&#
%e5.25 #&#
\begin{eqnarray}
\label{eqex3mollifierlowerbound} && \int_0^{ \infty}
\mathbh{1}_{
[ 0, 1 )
} \bigl( \lfloor s \rfloor_h \bigr) \cdot
\exp\biggl( \frac{ - 1 }{
(
1 - \vert \lfloor s \rfloor_h \vert^2
)
} \biggr) \,ds - \int_0^{ 1 }
\exp\biggl( \frac{ - 1 }{
(
1 - s^2
)
} \biggr) \,ds
\nonumber\\
&&\qquad = \int_0^{ \infty} \underbrace{
\mathbh{1}_{
[ 0, 1 )
} \bigl( \lfloor s \rfloor_h \bigr) \cdot
\exp\biggl( \frac{ - 1 }{
(
1 - \vert \lfloor s \rfloor_h \vert^2
)
} \biggr) - \mathbh{1}_{
[ 0, 1 )
} ( s ) \cdot
\exp\biggl( \frac{ - 1 }{
(
1 - s^2
)
} \biggr) }_{ \geq0 } \,ds
\nonumber\\
&&\qquad \geq\int_0^{
3^{ - 1 / 4 }
} \exp\biggl(
\frac{ - 1 }{
(
1 - \vert \lfloor s \rfloor_h \vert^2
)
} \biggr) - \exp\biggl( \frac{ - 1 }{
(
1 - s^2
)
} \biggr) \,ds
\nonumber\\
&&\qquad \geq\frac{ h }{ 2 } \cdot\biggl( \exp\biggl( \frac{ - 1 }{
(
1 - 0^2
)
}
\biggr) - \exp\biggl( \frac{ - 1 }{
(
1 -
\vert
\lfloor3^{ - 1 / 4 } \rfloor_h - h
\vert^2
)
} \biggr) \biggr)\\
&&\qquad\quad{} + \frac{
2 \cdot\lfloor3^{ - 1 / 4 } \rfloor_h
\cdot
e^{
- 1 /
( 1 - | \lfloor3^{ - 1 / 4 } \rfloor_h |^2 )
}
}{
[ 1 -
\vert \lfloor3^{ - 1 / 4 } \rfloor_h
\vert^2
]^2
}
\cdot\frac{
(
3^{ - 1 / 4 } -
\lfloor3^{ - 1 / 4 } \rfloor_h
)^2
}{ 2
}
\nonumber\\
&&\qquad \geq\frac{ h }{ 2 } \cdot\biggl( e^{ - 1 } - \exp\biggl(
\frac{ - 1 }{
(
1 -
[
3^{ - 1 / 4 } - 2 h
]^2
)
} \biggr) \biggr)\nonumber\\
&&\qquad \geq\frac{ h }{ 2 } \cdot\biggl(
e^{ - 1 } - \exp\biggl( \frac{ - 1 }{
(
1 -
[
{ 1 }/{ 2 }
]^2
)
} \biggr) \biggr) = h \cdot
\frac{
(
e^{ - 1 }
-
e^{ - 4 / 3 }
)
}{ 2 } > \frac{ h }{ 20 }\nonumber
\end{eqnarray}
for all $ h \in(0,\frac{1}{8}] $
where we used the estimate
$
3^{ - 1 / 4 } - 2 h
\geq
\frac{ 1 }{ 3^{ 1 / 4 } }
- \frac{ 1 }{ 4 }
\geq
\frac{ 1 }{ 2 }
$
for all $ h \in( 0, \frac{ 1 }{ 8 } ] $
in the penultimate inequality in
\eqref{eqex3mollifierlowerbound}.
Moreover, note that
\eqref{eqex3mollifiersecondderivative}
implies that
%
%e5.26 #&#
%e5.27 #&#
%e5.28 #&#
\begin{eqnarray}
\label{eqex3mollifierupperbound} && \int_0^{ \infty}
\mathbh{1}_{
[ 0, 1 )
} \bigl( \lfloor s \rfloor_h \bigr) \cdot
\exp\biggl( \frac{ - 1 }{
(
1 - \vert \lfloor s \rfloor_h \vert^2
)
} \biggr) \,ds \nonumber\\
&&\quad{}- \int_0^{ 1 }
\exp\biggl( \frac{ - 1 }{
(
1 - s^2
)
} \biggr) \,ds
\nonumber\\
&&\qquad \leq h + \int_0^{ 1 } \biggl\vert\exp
\biggl( \frac{ - 1 }{
(
1 - \vert \lfloor s \rfloor_h \vert^2
)
} \biggr) - \exp\biggl( \frac{ - 1 }{
(
1 - s^2
)
} \biggr)
\biggr\vert \,ds
\nonumber
\\[-8pt]
\\[-8pt]
\nonumber
&&\qquad \leq h + \sup_{ x \in(0,1) } \biggl[
\frac{
2 x \cdot
e^{ - 1 / ( 1 - x^2 ) }
}{
( 1 - x^2 )^2
} \biggr] \cdot h
\\
&&\qquad = h + \biggl[ \frac{
2 \cdot3^{ -1/4 } \cdot
e^{ - 1 / ( 1 - 3^{ - 1 / 2 } ) }
}{
( 1 - 3^{ - 1 / 2 } )^2
} \biggr] \cdot h\nonumber \\
&&\qquad= h + % \underbrace{
\biggl[
\frac{
6
}{
3^{ 1 / 4 } \cdot
( \sqrt{ 3 } - 1 )^2
\cdot
e^{ \sqrt{ 3 } / ( \sqrt{ 3 } - 1 ) }
} \biggr] % }_{ \leq1 }
\cdot h \leq2 h\nonumber
\end{eqnarray}
for all $ h \in(0,\infty) $.
Combining
\eqref{eqex3mollifierlowerbound}
and
\eqref{eqex3mollifierupperbound}
completes the proof
of Lem\-ma~\ref{lemmollifier}.
\end{pf}

We are now ready to prove Theorem~\ref{thmnonrate}.
Its proof uses
Lemma~\ref{lemmollifier}
as well as
Lemma~\ref{lemex2lowerbound}
in Section~\ref{secex2} above.

\begin{pf*}{Proof
of
Theorem~\ref{thmnonrate}}
First of all, note that
$X_1(t)
=
\int_1^t
\exp (
\frac{ - 1 }{
( s^2 - 1 )
}
)
\,ds$,
$\P$-a.s. for
all $ t \in[1,\infty) $.
Combining this
with \eqref{eqexample3Y1rep}
then shows that
\begin{eqnarray*}
 &&\E\bigl[ X_1(t) \bigr] - \E\bigl[ Y^h_1(t)
\bigr]\\
&&\qquad = \underbrace{ \int_1^t \exp\biggl( -
\frac{ 1 }{
s^2 - 1
} \biggr) - \mathbh{1}_{
(1, \infty)
} \bigl( \lfloor s
\rfloor_h \bigr) \cdot\exp\biggl( -\frac{ 1 }{
\lfloor s \rfloor_h^2 - 1
} \biggr) \,ds
}_{
\geq0
}
\\
&&\qquad\quad{} + \int_1^t \mathbh{1}_{
(1, \infty)
}
\bigl( \lfloor s \rfloor_h \bigr) e^{
-{ 1 }/{
(\lfloor s \rfloor_h^2 - 1)
}
} %\\ &
% \cdot
\\
&&\hspace*{22pt}\qquad\quad{}\times \E\biggl[ 1 -
\cos\biggl( \biggl( \smallint_{ 0 }^{ \infty}
\mathbh{1}_{
[0,1)
} \bigl( \lfloor u \rfloor_h \bigr)
% \cdot
e^{
-{ 1 }/{
(1 - \lfloor u \rfloor_h^2)
}
} \,du\\
&&\hspace*{181pt}{} - \smallint_{ 0 }^{ 1 }
e^{ -{ 1 }/{ (1 - u^2) } } \,du \biggr) % \cdot
e^{
W_2( \lfloor s \rfloor_h )^3
} \biggr) \biggr] \,ds
\\
&&\qquad \geq\int_{ 3 / 2 }^t \mathbh{1}_{
(1, \infty)
}
\bigl( \lfloor s \rfloor_h \bigr) % \cdot
e^{
-{1 }/{
(\lfloor s \rfloor_h ^2 - 1)
}
}
% \cdot
\\
&&\hspace*{51pt}{}\times\E\biggl[ 1 - \cos\biggl( \biggl( \smallint_{ 0 }^{ \infty}
\mathbh{1}_{
[0,1)
} \bigl( \lfloor u \rfloor_h \bigr)
% \cdot
e^{
-{ 1 }/{
(1 - \lfloor u \rfloor_h^2)
}
} \,du \\
&&\hspace*{176pt}{}- \smallint_{ 0 }^{ 1 }
e^{ -{ 1 }/{ (1 - u^2) } } \,du \biggr) % \cdot
e^{
W_2( \lfloor s \rfloor_h ) ^3
} \biggr) \biggr] \,ds
\end{eqnarray*}
for all $ t \in[\frac{3}{2}, \infty) $
and all $ h \in(0,\infty) $.
The estimate
$
\lfloor s \rfloor_h
\geq
\lfloor\frac{3}{2} \rfloor_h
\geq
\frac{3}{2} - h
\geq\frac{ 11 }{ 8 }
$
for all
$ s \in[\frac{3}{2},\infty) $,
$ h \in(0,\frac{1}{8}] $
and
Lemmas~\ref{lemmollifier}
and~\ref{lemex2lowerbound}
therefore show that
{\fontsize{10}{12}{\selectfont{
\begin{eqnarray*}
&& \E\bigl[ X_1(t) \bigr] - \E\bigl[ Y^h_1(t)
\bigr]
\\
& &\qquad\geq\exp\biggl( -\frac{ 1 }{
{ 121 }/{ 64 } - 1
} \biggr) %\\ &
% \cdot
\\
&&\qquad\quad{}\times\int
_{ { 3 }/{ 2 } }^v \E\biggl[ 1 - \cos\biggl( \underbrace{
\biggl( \smallint_{ 0 }^{ \infty} \mathbh{1}_{
[0,1)
} \bigl(
\lfloor u \rfloor_h \bigr) e^{
- 1 /
(
1 -
\vert \lfloor u \rfloor_h \vert^2
)
} \,du -
\smallint_{ 0 }^{ 1 } e^{ - 1 / ( 1 - u^2 ) } \,du \biggr)
}_{z
{ h }/{ 20 } \leq\cdots
\leq2 h\
\mathrm{due\ to\ Lemma~\scriptsize{\ref{lemmollifier}}}
} % \cdot
\\
&&\hspace*{282pt}{}\times e^{
W_2( \lfloor s \rfloor_h )^3
} \biggr) \biggr] \,ds
\\
&&\qquad \geq e^{
-{ 64 }/{
57
}
}\\
&&\qquad\quad{}\times \int_{ { 3 }/{ 2 } }^v \exp
\biggl( \frac{ - 8 }{ \lfloor s \rfloor_h }\\
&&\hspace*{77pt}{} \times\biggl\vert\ln\biggl
( { \pi}\Big/\biggl(
2
\biggl(
\smallint_{ 0 }^{ \infty}
\mathbh{1}_{
[0,1)
} (
\lfloor u \rfloor_h
)
\cdot
e^{
- 1 /
(
1 -
\vert \lfloor u \rfloor_h \vert^2
)
}
\,du
\\
&&\hspace*{202pt}{}-
\smallint_{ 0 }^{ 1 }
e^{ - 1 / ( 1 - u^2 ) }
\,du
\biggr)
\biggr) \biggr)
\biggr\vert^{ 2 / 3 } \biggr) \,ds
\\
&&\qquad \geq\frac{
( v - { 3 }/{ 2 } )
}{ 4 } \cdot\exp\biggl( -\frac{ 64 }{ 11 } \biggl\vert
\ln\biggl( \frac{ 10 \pi}{
h
} \biggr) \biggr\vert^{ 2 / 3 } \biggr)
\end{eqnarray*}}}}
\hspace*{-2pt}for all
$
h \in ( 0,
\min\{ \frac{1}{8},
\frac{ \pi}{ 4 }
\exp( - v^{ 3 / 2 } )
\}
]
$,
$ t \in[v,\infty) $
and all
$ v \in[\frac{3}{2}, \infty) $.
Hence, we finally obtain that
%
%e5.29 #&#
\begin{eqnarray}
%&
&&\E\bigl[ X_1(t) \bigr] - \E\bigl[
Y^h_1(t) \bigr]
\nonumber
\\[-8pt]
\\[-8pt]
\nonumber
&&\qquad\geq\exp\biggl( - \ln(8) -
\frac{ 64 }{ 11 } \bigl\vert\ln( 10 \pi) \bigr\vert^{ 2 / 3 } -
\frac{ 64 }{ 11 } \bigl\vert\ln( h ) \bigr\vert^{ 2 / 3 } \biggr)
\end{eqnarray}
for all
$
h \in( 0,
\frac{ 1 }{ 22 }
]
$
and all
$ t \in[2,\infty) $.
This completes the proof
of Theorem~\ref{thmnonrate}.
\end{pf*}

%It is somewhat surprising
%that there exist SDEs with globally bounded and
%smooth coefficients for which the Euler-Maruyama
%approximations converge
%but not with any
%arbitrarily small positive
%rate of convergence.
%Next we address the impact of
%this theoretical result on applications.
%

\begin{figure}

\includegraphics{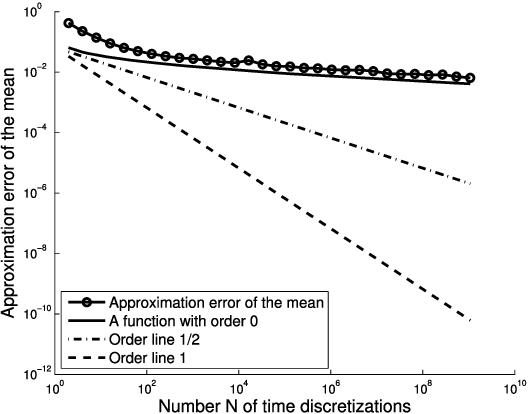}

\caption{The norm
$
\|
\E[ X(T) ] - \E[ Y^{{T}/{N}}(T) ]
\|
$
of the difference between the mean of
the solution of the SDE~\protect\eqref{eqex3SDE}
and the mean
of the Euler--Maruyama approximations~\protect\eqref{eqex3Euler}
for $ T = 2 $
and
$
N \in
\{
2^1,
2^2,
\ldots,
2^{ 29 },
2^{ 30 }
\}
$.
The function with convergence order $0$ is given by~\protect\eqref
{eqorder0function}.}

\label{fnonrate}
\end{figure}

In the next step, we illustrate
the lower bound
on the weak approximation error
in Theorem~\ref{thmnonrate}
by a numerical
simulation.
% on a computer.
More precisely,
we ran Monte Carlo simulations
and
approximatively calculated
the quantity
$
\|
\E[ X(T) ] - \E[ Y^{{T}/{N}}(T) ]
\|
$
for $ T = 2 $
and
$
N \in
\{
2^1,
2^2,
\ldots,
2^{ 29 },
2^{ 30 }
\}
$.
We approximated these
differences of expectations
with an average
over $ 100\mbox{,}000 $
independent Monte
Carlo realizations.
Moreover,
we discretized the integrals
$
X_1(2)
=
\int_1^2
\exp (
\frac{ - 1 }{
( s^2 - 1 )
}
)
\,ds
$
and
$
X_3(2)
=
\int_0^1
\exp (
\frac{ - 1 }{
( 1 - s^2 )
}
)
\,ds
$
in the exact solution with
a uniform grid and mesh size
$
\frac{2}{2^{31}}
=
2^{ - 30 }
$.
Figure~\ref{fnonrate}
depicts the resulting graph.

%f1 #&#

In addition to
the weak approximation
error
$
\|
\E[ X(T) ] - \E[ Y^{{T}/{N}}(T) ]
\|
$
for $ T = 2 $
and
$
N \in
\{
2^1,
2^2,
\ldots,
2^{ 29 },
2^{ 30 }
\}
$,
we also plotted
the function
%
%e5.30 #&#
\begin{eqnarray}
\label{eqorder0function}\quad &&\bigl\{2^{1},2^{2},\ldots,2^{30}
\bigr\}\ni N
\nonumber
\\[-8pt]
\\[-8pt]
\nonumber
&&\qquad\mapsto\frac{1}{15\cdot(\ln(N))^{{1}/{3}}} \exp
\biggl(-\frac{1}{2T} \biggl(
\ln(N) -\frac{1}{2T} \bigl(\ln(N) \bigr)^{{2}/{3}}
\biggr)^{{2}/{3}} \biggr)\in(0,1]
\end{eqnarray}
(a function with order $ 0 $),
the function
$
\{ 2^{1}, 2^{2}, \ldots, 2^{30} \}
\ni N
\mapsto\frac{ 1 }{ 15 \cdot\sqrt{ N } }
\in
(0,1]
$
(order line $ \frac{ 1 }{ 2 } $)
and
the function
$
\{ 2^{1}, 2^{2}, \ldots, 2^{30} \}
\ni N
\mapsto\frac{ 1 }{ 15 \cdot N }
\in
(0,1]
$
(order line $ 1 $)
in Figure~\ref{fnonrate}.
In the standard literature
in computational stochastics
(see, e.g., Kloeden
and Platen~\cite{kp92})
the Euler--Maruyama scheme
is shown to converge
in the numerically weak sense with
order $ 1 $ if the coefficients
of the SDE are smooth and
globally Lipschitz
continuous
(see Chapter~8 in
Kloeden and Platen~\cite{kp92} for the
precise assumptions)
and, therefore, the order line~$ 1 $
is plotted in Figure~\ref{fnonrate}.
Moreover, the function with
order~$ 0 $
is included in Figure~\ref{fnonrate}
so that one can compare the graph visually with a function
which has convergence order $0$.
According to our simulations,
the approximation error for the mean
$
\E[X(2)]
$
does not drop far below
$ \frac{ 1 }{ 100 } $
even for
$
N = 2^{ 30 } > 10^9
$
time discretizations.
This indicates that calculating the mean
$
\E[X(T)]
$
with the Euler--Maruyama method
up to a high precision requires
a huge computational effort.
In particular, this suggests for
applications that an approximation
cannot, in general, be assumed to
be very close to the
exact value even after
a very high computational effort.

\section*{Acknowledgements}

We gratefully acknowledge
Verena B\"ogelein,
Sonja Cox,
Weinan E,
Alessandra Lunardi,
{\'E}tienne Pardoux,
Michael R\"{o}ckner
and
Tobias Weth
for helpful remarks
and
for pointing
out useful
references to us.
Special thanks are due to
Shige Peng
for fruitful discussions about
questions on uniqueness of
viscosity solutions,
in particular,
for pointing out his
quite instructive
book~\cite{Peng2010} to us.

%

% imsref loaded by akundreckaite, 2014-01-10 08:58:40
% imsref loaded by akundreckaite, 2014-01-10 09:47:19
%


%suskaldyti doi

\printaddresses

\end{document}